\documentclass[11pt,leqno]{amsart}
\usepackage{euscript, amssymb, amsmath, amsthm}
\usepackage{graphicx}
\usepackage{float}
\usepackage{pdfsync}
\usepackage{subcaption}
\usepackage{color}

\newcommand{\ep}{\varepsilon}

\date{\today}

\title[Weighted Essentially Non-Oscillatory Schemes]{\Small{A New Adaptive Weighted Essentially Non-Oscillatory WENO-$\theta$ Scheme for Hyperbolic Conservation Laws}}

 \theoremstyle{definition}
 
 \theoremstyle{remark}
 \newtheorem{rem}{Remark}[section]
 \theoremstyle{remark}
 \newtheorem{exam}{Example}[section]
 \numberwithin{equation}{section}

\begin{document}

\maketitle

\centerline{\scshape Chang-Yeol Jung and Thien Binh Nguyen}
\centerline{\small{Department of Mathematical Sciences, School of
Natural Science,}} \centerline{\small{Ulsan National Institute of
Science and Technology,}} \centerline{\small{UNIST-gil 50, Ulsan
689-798, Republic of Korea}} \centerline{\small{cjung@unist.ac.kr,\
thienbinh84@unist.ac.kr}}

\begin{abstract}
A new adaptive weighted essentially non-oscillatory WENO-$\theta$
scheme in the context of finite difference is proposed. Depending on
the smoothness of the large stencil used in the reconstruction of
the numerical flux, a parameter $\theta$ is set adaptively to switch
the scheme between a 5th-order upwind and 6th-order central
discretization. A new indicator $\tau^{\theta}$ measuring the
smoothness of the large stencil is chosen among two candidates which
are devised based on the possible highest-order variations of the
reconstruction polynomials in $L^2$ sense. In addition, a new set of
smoothness indicators $\tilde{\beta}_k$'s of the sub-stencils is
introduced. These are constructed in a central sense with respect to
the Taylor expansions around the point $x_{j}$.

Numerical results show that the new scheme combines good properties
of both 5th-order upwind schemes, e.g., WENO-JS (\cite{JS}), WENO-Z
(\cite{BCCD}), and 6th-order central schemes, e.g., WENO-NW6
(\cite{YC}), WENO-CU6 (\cite{HWA}). In particular, the new scheme
captures discontinuities and resolves small-scaled structures much
better than the 5th-order schemes; overcomes the loss of accuracy
near some critical regions and is able to maintain symmetry which
are drawbacks detected in the 6th-order ones.
\end{abstract}


{\small{\textbf{Keywords.} Hyperbolic conservation laws, Euler
equations, shock-capturing methods, Weighted essentially
non-oscillatory (WENO) schemes, Adaptive upwind-central schemes,
Smoothness indicators.\\}}


\section{Introduction}

In this work, we consider the following one-dimensional hyperbolic
conservation law
\begin{align}\label{mod_con_law}
\begin{cases}
&\mathbf{u}_t+\mathbf{f}(\mathbf{u})_x=0,\quad x\in\mathbb{R},\ t>0,\\
&\mathbf{u}(x,0)=\mathbf{u}_0(x),
\end{cases}
\end{align}
where $\mathbf{u}=(u_1,\ldots,u_m)^T$ is an $m$-dimensional vector
of conserved quantities and its flux $\mathbf{f}(\mathbf{u})$ is a
vector-valued function with $m$ components, $x$ and $t$ denote space
and time, respectively. Eq. (\ref{mod_con_law}) is called hyperbolic
assuming that all eigenvalues $\lambda_s$'s of the Jacobian
$A(\mathbf{u})=\partial{\mathbf{f}}/\partial{\mathbf{u}}$ are real
and the set of all eigenvectors $\mathbf{r}_s$'s is complete.

It is well-known that shocks and discontinuities may develop in the
solution of Eq. (\ref{mod_con_law}) even if the initial condition is
smooth. Thus classical numerical methods which depend on Taylor
expansions in general do not work in this case. As a result, there
exist spurious oscillations near these discontinuities.

In order to overcome this difficulty, in \cite{Ha83}, \cite{Ha84}
Harten introduced the Total Variation Diminishing (TVD) schemes
which are of high-order resolutions as well as oscillations free.
The schemes are constructed based on the principle that the total
variation of the numerical approximation must be non-increasing in
time. A drawback is that TVD schemes are only at most first-order
near smooth extrema (see \cite{OC}). Later on, Harten et al. in
\cite{HOEC}, \cite{HO}, and \cite{HEOC} tried to tackle this
disadvantage by relaxing the TVD condition and allowing spurious
oscillations in the order of the truncation error to occur but the
$\mathcal{O}(1)$ Gibbs-like ones are essentially prevented. Thus
these new schemes were named essentially non-oscillatory (ENO). For
an $r$th-order ENO scheme, only the smoothest stencil is chosen
among $r$ candidates to approximate the numerical flux. The
smoothness of the solution on each stencil is determined by an
indicator of smoothness. Later on, Liu, Osher, and Chan (\cite{LOC})
upgraded ENO schemes and introduced the Weighted ENO (WENO) by
combining all stencil candidates (hereafter sub-stencils) in the
numerical flux approximation. Here, a nonlinear weight is assigned
to each sub-stencil to control its contribution in the procedure.
WENO schemes maintain the essentially non-oscillatory property of
the ENO near discontinuities and outperform the latter in smooth
regions where the accuracy order is increased to $(r+1)$th-order if
$r$ sub-stencils are used. Consequently, Jiang and Shu (see
\cite{JS}, also \cite{Sh03}, and the review \cite{Sh09}) constructed
WENO schemes in the framework of finite difference and further
improved the order to $(2r-1)$th in smooth regions by introducing a
new class of smoothness indicators. Hereafter, we denote WENO-JS for
the 5th-order finite difference WENO developed in \cite{JS}. In
\cite{BS}, \cite{SZ} higher order than 5th-order WENO schemes are
given.

Since the introduction of WENO, many improvements and derivatives of
the schemes have been developed and introduced. Henrick et al. in
\cite{HAP} carefully analyzed the necessary and sufficient
conditions of the nonlinear weights and found that WENO-JS does not
achieve the designed 5th-order but reduces to only 3rd-order in
cases where the first and third derivatives of the flux do not
simultaneously vanish (e.g., $f^{\prime}(x_j)=0$ but
$f^{\prime\prime\prime}(x_j)\ne0$ for the scalar case of Eq.
(\ref{mod_con_law})). They then suggested an improved version which
is called mapped WENO, abbreviated by WENO-M. By using a mapping on
the nonlinear weights, WENO-M satisfies the sufficient condition on
which WENO-JS fails and obtains optimal order near simple smooth
extrema. In a different approach on the construction of the
nonlinear weights, in \cite{BCCD} Borges et al. introduced the
5th-order WENO-Z scheme. Here, the authors also measured the
smoothness of the large stencil which comprises all sub-stencils and
incorporated this in devising the new smoothness indicators and
nonlinear weights. It was proven numerically that WENO-Z is less
dissipative than WENO-JS and more efficient than WENO-M,
respectively. It was also checked that WENO-Z attains 4th-order near
simple smooth extrema comparing with 3rd-order of WENO-JS. For
higher order WENO-Z schemes, we refer readers to \cite{CCD}. Another
approach to improve WENO schemes is the new designs of the
smoothness indicators. In \cite{HKLY}, $L^1$-norm based smoothness
indicators are suggested, and the ones devised from Lagrange
interpolation polynomials are given in \cite{Fa}, and \cite{FSTY}.
See also \cite{FHW} for a new mapped WENO scheme.

We notice that for a general flux where the signs of the eigenvalues
of the Jacobian $A(\mathbf{u})$ are not uniform throughout the
domain, a flux splitting technique, for example, the global or local
Lax-Friedrichs or the Roe with entropy fix (see \cite{JS} and the
references therein) is needed. This increases the number of grid
points in the numerical flux approximating procedure by one. We take
the 5th-order WENO-JS scheme for example, the total number of grid
points used in the reconstruction for both positive and negative
fluxes will be six instead of five. We also note that with these six
points, one can indeed improve the scheme up to 6th-order in smooth
regions. The difficulty of this approach lies in the dispersive
nature of a central scheme if six points are employed. In this case,
oscillations are expected to occur near discontinuities. In
\cite{YC} (see also \cite{CFY} for the boundary condition
treatment), Yamaleev and Carpenter for the first time introduced a
6th-order WENO scheme by adding one more sub-stencil into the
numerical flux approximation. We denote this scheme WENO-NW6. For
this most downwind sub-stencil, an \emph{ad hoc} treatment on the
smoothness indicator $\beta_3$ was suggested. The idea is originated
from that of Mart\'{i}n et al. in \cite{MTW}. In order that
oscillations do not happen, $\beta_3$ is computed using the
information of $\mathbf{f}(\mathbf{u})$ on all grid points of the
large stencil, i.e, six points. Hence, the sub-stencil only plays
roles in case the solution is smooth over this large stencil. In a
similar manner, recently Hu, Wang, and Adams in \cite{HWA} proposed
an adaptive central-upwind WENO-CU6 scheme which switches between a
5th-order upwind and 6th-order central WENO scheme automatically.
The difference of their work from that given in \cite{YC} is that
$\beta_3$ is defined via a Lagrange interpolating polynomial of
degree five over the large stencil. In \cite{HA11}, the authors
successfully applied WENO-CU6 in the LES simulation of scale
separation. Other hybrid WENO schemes can be found in, for examples,
\cite{CD}, or \cite{LQ}, \cite{HP}, etc.

A drawback of the presented 6th-order WENO schemes (i.e., WENO-NW6,
WENO-CU6) is that they suffer from a loss of accuracy near the
smooth critical region which is just behind another one where the
first derivative of the flux is undefined. To illustrate this, we
consider Eq. (\ref{mod_con_law}) in a scalar case where $f(u)=u$ in
the following example.

\begin{exam}\label{exam_1}
\begin{align}\label{ex}
\begin{cases}
&u_t+u_x=0,\quad x\in(-1,1),\\
&u_0(x)=\max(-\sin(\pi x),0),
\end{cases}
\end{align}
subject to periodic boundary conditions.

We approximate the solution of (\ref{ex}) by the WENO-JS, WENO-Z,
WENO-NW6, and WENO-CU6 schemes. The results at time $t=2.4$ with
$200$ grid intervals are plotted in Fig. \ref{fig_exam_1} with the
critical region zoomed in. It is clearly shown the above mentioned
defect of the WENO-NW6 and WENO-CU6 schemes. Near the smooth
critical region, we note that these schemes are worse than both
WENO-JS and WENO-Z. Since there are many problems whose solution
often exhibits the same behavior as mentioned above, we notice that
this loss of accuracy is an important issue.

\begin{figure}
  \begin{center}
    \begin{tabular}{cc}
      \resizebox{70mm}{!}{\includegraphics{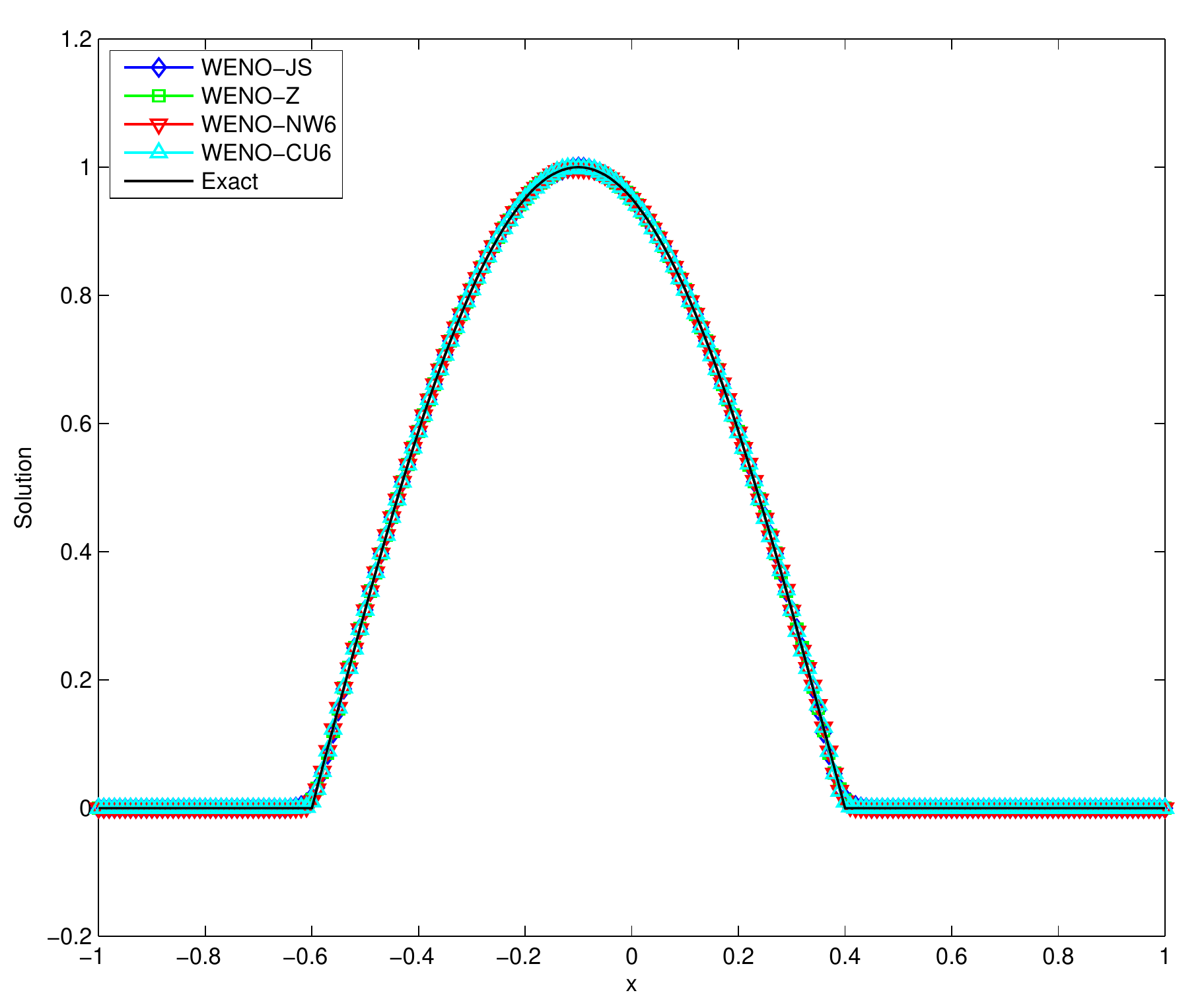}} &
      \resizebox{70mm}{!}{\includegraphics{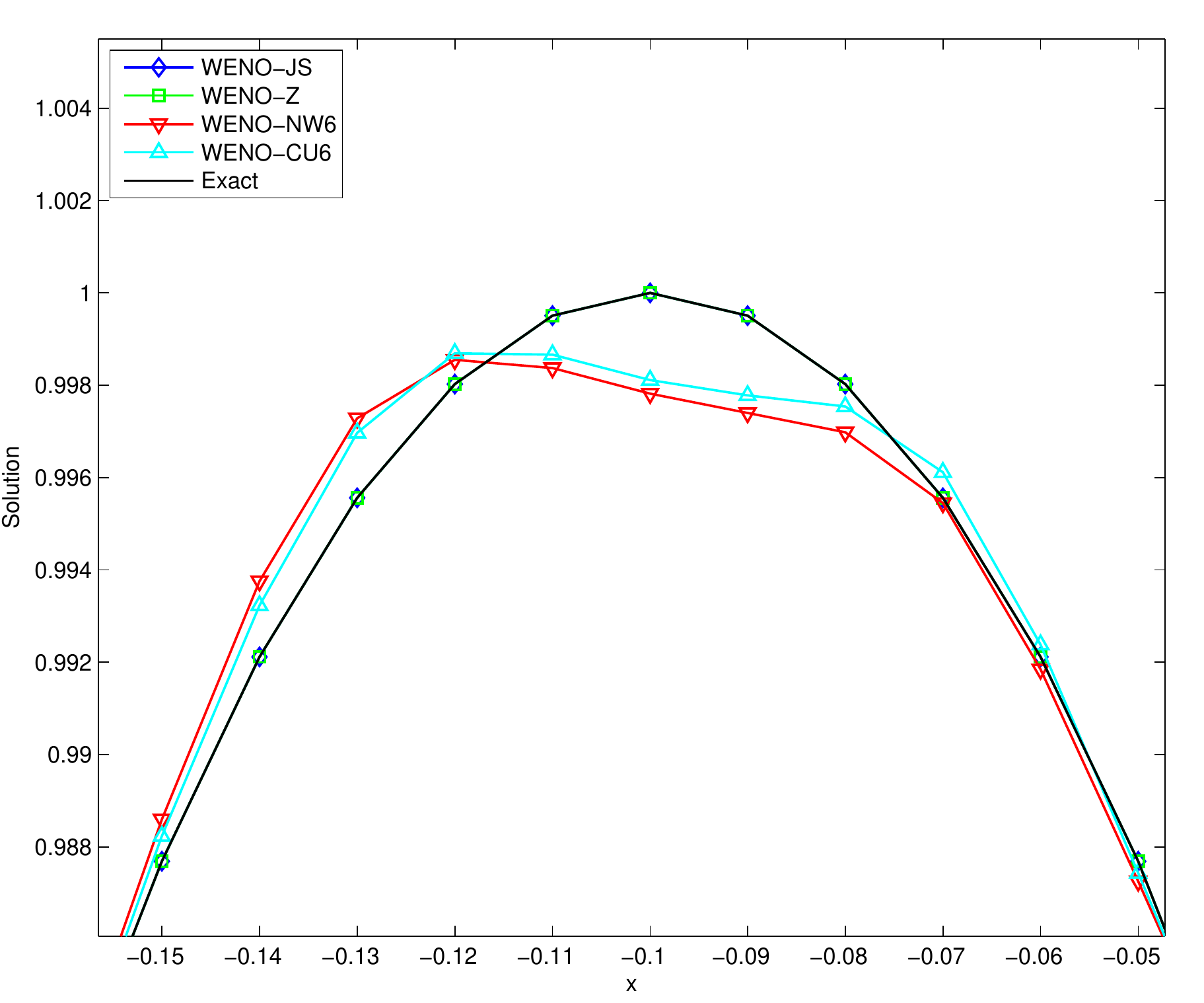}}
    \end{tabular}
    \caption{\small{Left: Numerical solutions of Eq. (\ref{ex}) at time $t=2.4$ obtained from different WENO schemes. Right: Zoom near the critical region.}} \label{fig_exam_1}
  \end{center}
\end{figure}

\end{exam}

Our goal in this work is to construct a new WENO scheme which
overcomes the drawback of WENO-NW6 and WENO-CU6 presented in the
previous example. For this, we introduce a different switching
mechanism between a 5th-order upwind and 6th-order central scheme.
Unlike the WENO-NW6 or WENO-CU6 scheme in which the change depends
on the smoothness indicator of the most downwind sub-stencil, in our
scheme, whether the scheme is upwind or central is due to the
smoothness indicator of the large stencil. Moreover, instead of
using all six points for the indicator $\beta_3$, we reduce the
number of points down to only four. The reason for this is explained
in the below section. We also introduce a new set of smoothness
indicators which are constructed in a central sense in Taylor
expansions with respect to the point $x_{j}$. The feature of our new
scheme is that it automatically switches between a 5th-order upwind
scheme near discontinuities to prevent spurious oscillations, and a
6th-order central scheme in smooth regions which improves the loss
of accuracy of the WENO-NW6 and WENO-CU6 schemes. Moreover, it is
shown in the numerical results below that the new scheme maintains
symmetry in the solutions much better than the 6th-order ones.

We organize our paper as follows. In section 2, we summarize the
mentioned above WENO schemes which relate to our work. From there,
we construct our new scheme in section 3. In this section, we first
start with the new definition of the central smoothness indicators.
We then introduce a new switching mechanism for 5th-order upwind and
6th-order central scheme. Numerical results comparing the
performances of all schemes are presented in section 4. Finally, we
close our discussions with a conclusion section.

\section{Summary on Finite Difference WENO Schemes}

For simplicity, we consider Eq. (\ref{mod_con_law}) in a scalar
case. We rewrite the equation as follows,
\begin{align}\label{mod_con_scalar}
\begin{cases}
&u_t+f(u)_x=0,\quad x\in[x_l,x_r],\\
&u(x,0)=u_0(x).
\end{cases}
\end{align}

We first define a uniformly spatial grid $x_j=x_l+j\Delta x$,
$j=0,\ldots,N$, where $\Delta x$ is the grid size. We denote the
interval $I_j=[x_{j-\frac{1}{2}},x_{j+\frac{1}{2}}]$ where
$x_{j\pm\frac{1}{2}}=x_j\pm\frac{\Delta x}{2}$ are the interfaces of
$I_j$. We also denote all quantities with a subscript $(\cdot)_j$
their grid values at $x_j$, for examples, $u_j=u(x_j,\cdot)$,
$f_j=f(u_j)$, etc.; and so as with a subscript
$(\cdot)_{j+\frac{1}{2}}$ for the quantities at the interface
$x_{j+\frac{1}{2}}$. Whether these quantities are exact or
approximate depends on particular circumstances.

We denote $h(x)$ the numerical flux function defined as follows,
\begin{align}\label{f_ave}
f(u(x,\cdot))=\bar{h}(x):=\dfrac{1}{\Delta x}\int_{x-\frac{\Delta
x}{2}}^{x+\frac{\Delta x}{2}}h(y)dy.
\end{align}
Evaluating Eq. (\ref{mod_con_scalar}) at grid point $x_j$, we obtain
the semi-discretized form as follows,
\begin{align}\label{mod_semi}
\dfrac{du_j}{dt}=-\dfrac{\partial f}{\partial
x}\bigg|_{x=x_j}=-\dfrac{h(x_{j+\frac{1}{2}})-h(x_{j-\frac{1}{2}})}{\Delta
x}=:\mathcal{L}(u).
\end{align}
It is noticed that Eq. (\ref{mod_semi}) is exact since there are no
approximating errors in the formula.

\subsection{Time Integration}

We first mention about time advancing for Eq. (\ref{mod_semi}).
Following \cite{GS} and the references therein, for all below WENO
schemes, we employ the 3rd-order TVD Runge-Kutta method as below.
For TVD, we mean that the time integrator follows the same property
mentioned above. Here, $\Delta t$ is the time step satisfying some
proper CFL condition.
\begin{align}\label{3rd_RK}
\begin{split}
&u^{(1)}=u^n+\Delta t\hat{\mathcal{L}}(u^n),\\
&u^{(2)}=\dfrac{3}{4}u^n+\dfrac{1}{4}u^{(1)}+\dfrac{1}{4}\Delta t
\hat{\mathcal{L}}(u^{(1)}),\\
&u^{n+1}=\dfrac{1}{3}u^n+\dfrac{2}{3}u^{(2)}+\dfrac{2}{3}\Delta
t\hat{\mathcal{L}}(u^{(2)}),
\end{split}
\end{align}
where $\hat{\mathcal{L}}(u)$ obtained from some method is an
approximation of the spatial operator $\mathcal{L}(u)$. In
particular, see the below WENO discretizations
$\hat{\mathcal{L}}^5(u)$ in Eq. (\ref{approx_dis}) where
$\hat{f}^5_{j+\frac{1}{2}}=\hat{f}_{j+\frac{1}{2}}$ follows Eq.
(\ref{weno}) and $\hat{\mathcal{L}}^6(u)$ in Eq.
(\ref{approx_dis_6}) where $\hat{f}^6_{j+\frac{1}{2}}$ is defined in
Eq. (\ref{fh6_lin}) with $\gamma_k$'s replaced by $\omega_k$'s.

We now proceed to the discussions on the spatial discretizations.

\subsection{5th-order Upwind WENO Reconstruction}\label{recon}

We notice that for simplicity, we can assume that
$f^{\prime}(u)\ge0$ over the whole computational domain. In case
there is a change in signs of $f^{\prime}(u)$, a flux splitting
technique is invoked. We discuss this in Remark \ref{rmk_flux_split}
below.

Originally, WENO schemes were constructed in the context of finite
volume (see \cite{LOC}). Thanks to Lemma $3.1$ given in \cite{Sh09},
the schemes can be transformed into finite difference through
relation (\ref{f_ave}). The $\bar{h}_j$ is called an average value
of the numerical flux $h(x)$ over the interval $I_j$. We then seek
for an approximating polynomial $\hat{f}^5(x)$ of degree four of
$h(x)$ as below
\begin{align}\label{fh_poly}
h(x)\approx\hat{f}^5(x)=a_0+a_1x+a_2x^2+a_3x^3+a_4x^4,
\end{align}
over the large stencil
$S^5=\{x_{j-2},x_{j-1},x_j,x_{j+1},x_{j+2}\}$. We note that $S^5$ is
chosen biased to the left with respect to the point
$x_{j+\frac{1}{2}}$ for the stability purpose. Hence, the scheme is
in an upwind sense.

Replacing the integrand $h(x)$ in Eq. (\ref{f_ave}) by its
approximation $\hat{f}^5(x)$ in Eq. (\ref{fh_poly}) and evaluating
at $x_k$, $k=j-2,\ldots,j+2$, we can uniquely determine the
coefficients $a_k$'s, $k=0,\ldots,4$. Since the procedure is via the
average $\bar{h}_k=\bar{h}(x_k)=f(u(x_k,\cdot))$ in Eq.
(\ref{f_ave}), it is called reconstruction and $\hat{f}^5(x)$ is the
reconstruction polynomial. Evaluating $\hat{f}^5(x)$ at
$x_{j+\frac{1}{2}}$, we obtain the approximation of
$h_{j+\frac{1}{2}}$ as follows,
\begin{align}\label{fh_5}
\hat{f}^5_{j+\frac{1}{2}}=\dfrac{2}{60}f_{j-2}-\dfrac{13}{60}f_{j-1}
+\dfrac{47}{60}f_j+\dfrac{27}{60}f_{j+1}-\dfrac{3}{60}f_{j+2}.
\end{align}
Here, we recall $f_k=f(u_k)=f(u(x_k,\cdot))$.

To justify the approximation error, we denote the polynomial $H(x)$
such that
\begin{align}\label{H_h}
H^{\prime}(x)=h(x).
\end{align}

We then deduce from Eq. (\ref{f_ave}) that
\begin{align}\label{f_jH}
f(u(x,\cdot))=\dfrac{1}{\Delta x}\int_{x-\frac{\Delta
x}{2}}^{x+\frac{\Delta x}{2}}H^{\prime}(y)dy=\dfrac{H(x+\frac{\Delta
x}{2})-H(x-\frac{\Delta x}{2})}{\Delta x}.
\end{align}

Substituting Eq. (\ref{f_jH}) into the approximation (\ref{fh_5}),
evaluating at $x_k$, $k=j-2,\ldots,j+2$, and applying Taylor
expansions of $H(x)$ at $x=x_{j+\frac{1}{2}}$, we obtain the
following truncation error,
\begin{align}\label{fh_5_tay}
\begin{split}
\hat{f}^5_{j+\frac{1}{2}}&=\dfrac{1}{60\Delta
x}(-2H_{j-\frac{5}{2}}+15H_{j-\frac{3}{2}}-60H_{j-\frac{1}{2}}+20H_{j+\frac{1}{2}}+30H_{j+\frac{3}{2}}-3H_{j+\frac{5}{2}})\\
&=H^{\prime}_{j+\frac{1}{2}}-\dfrac{1}{60}\dfrac{d^6H}{dx^6}\bigg|_{x=x_{j+\frac{1}{2}}}\Delta
x^5+\mathcal{O}(\Delta
x^6)
=h_{j+\frac{1}{2}}-\dfrac{1}{60}\dfrac{\partial^5f}{\partial
x^5}\bigg|_{x=x_j}\Delta x^5+\mathcal{O}(\Delta x^6).
\end{split}
\end{align}

The last equality in Eq. (\ref{fh_5_tay}) is justified as follows.
Thanks to the relation in (\ref{f_jH}), by a Taylor expansion around
$x_{j+\frac{1}{2}}$ we have that
{\small
\begin{align}\label{fp5_jH}
\begin{split}
\dfrac{\partial^5 f}{\partial x^5}\bigg|_{x=x_j}&=\dfrac{1}{\Delta
x}\left[\dfrac{d^5H}{dx^5}\bigg|_{x=x_{j+\frac{1}{2}}}-\dfrac{d^5H}{dx^5}\bigg|_{x=x_{j+\frac{1}{2}}-\Delta
x}\right]
=\dfrac{d^6H}{dx^6}\bigg|_{x=x_{j+\frac{1}{2}}}+\mathcal{O}(\Delta
x).
\end{split}
\end{align}
}
Together with Eq. (\ref{H_h}), we deduce the last equality in Eq.
(\ref{fh_5_tay}).

Similarly, we have
\begin{align}
\hat{f}^5_{j-\frac{1}{2}}&=\dfrac{2}{60}f_{j-3}-\dfrac{13}{60}f_{j-2}+\dfrac{47}{60}f_{j-1}+\dfrac{27}{60}f_{j}-\dfrac{3}{60}f_{j+1}\\
&=h_{j-\frac{1}{2}}-\dfrac{1}{60}\dfrac{\partial^5f}{\partial
x^5}\bigg|_{x=x_j}\Delta x^5+\mathcal{O}(\Delta x^6).\nonumber
\end{align}

Hence, we have
\begin{align}\label{approx_dis}
\dfrac{du_j}{dt}\approx-\dfrac{\hat{f}^5_{j+\frac{1}{2}}-\hat{f}^5_{j-\frac{1}{2}}}{\Delta
x}=:\hat{\mathcal{L}}^5(u).
\end{align}
The scheme is 5th-order of accuracy in space.

\begin{figure}
  \begin{center}
    \begin{tabular}{c}
      \resizebox{105mm}{!}{\includegraphics{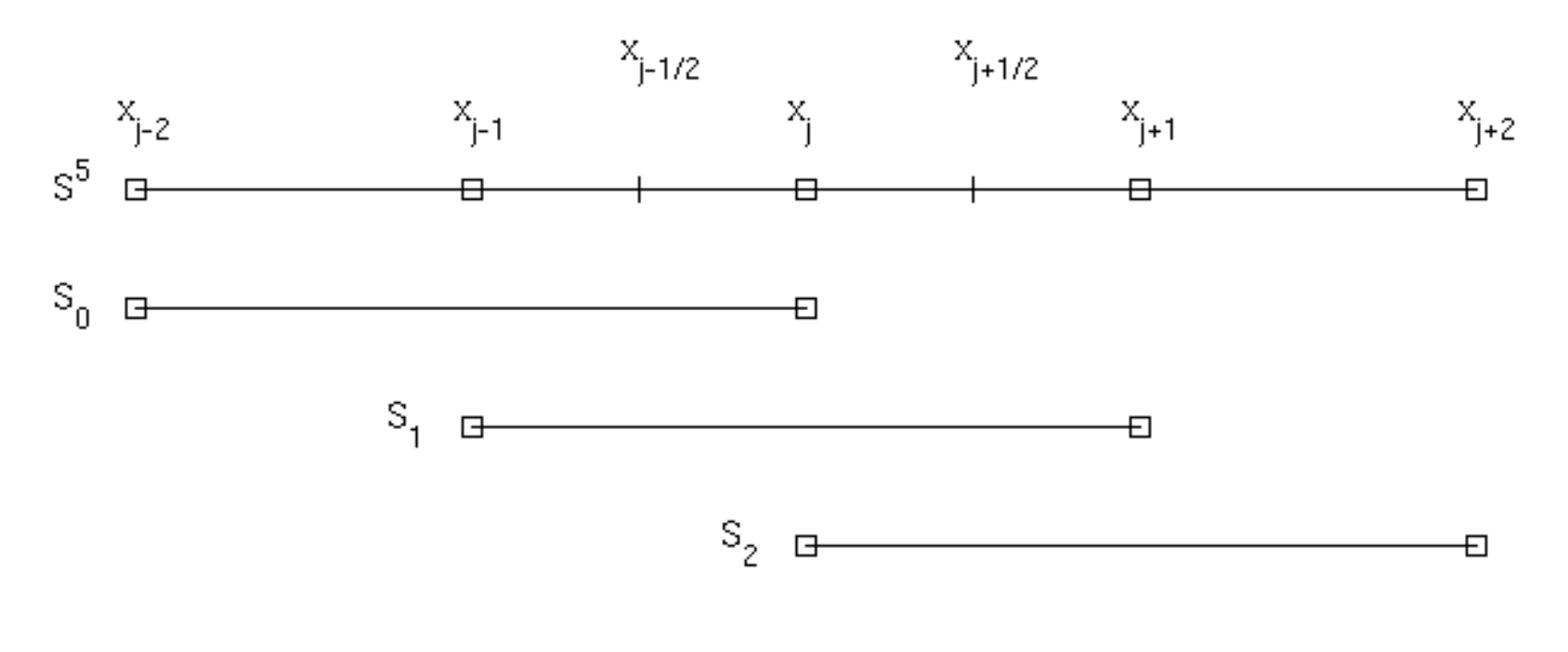}}
    \end{tabular}
    \caption{\small{Stencils for 5th-order WENO schemes.}}
\label{fig_sten_5}
  \end{center}
\end{figure}

For non-smooth solutions, we employ WENO reconstruction. The idea of
WENO schemes is that, instead of the 5-point stencil $S^5$, a convex
combination of three 3-point sub-stencils are facilitated for an
adaptive choice of candidates for the reconstruction. That is,
\begin{align}\label{weno}
\hat{f}_{j+\frac{1}{2}}=\sum_{k=0}^2\omega_k\hat{f}^k_{j+\frac{1}{2}},
\end{align}
where $\hat{f}^k_{j+\frac{1}{2}}$'s are defined below and  $\omega_k$ is the non-linear weight satisfying $\omega_k\ge0,\
\forall k$ and
\begin{align}\label{om_norm}
\sum_{k=0}^2\omega_k=1.
\end{align}
The necessity of non-negative nonlinear weights is discussed in
\cite{LSZ} and in \cite{SHS02} for practical implementations. And
$\hat{f}^k_{j+\frac{1}{2}}$ is the approximation of
$h_{j+\frac{1}{2}}$ by the reconstruction polynomial
$\hat{f}^k(x)=b_0+b_1x+b_2x^2$ over the sub-stencil $S_k$,
$k=0,1,2$. Here, $S_0=\{x_{j-2},x_{j-1},x_j\}$,
$S_1=\{x_{j-1},x_{j},x_{j+1}\}$, and $S_2=\{x_{j},x_{j+1},x_{j+2}\}$
(see Fig. \ref{fig_sten_5}). Carrying a similar process as for the
large stencil $S^5$, we find that around $x_j=0$, {\small
\begin{align}\label{f0_up}
\hat{f}^0(x)=\dfrac{-f_{j-2}+2f_{j-1}+23f_{j}}{24}+\left(\dfrac{f_{j-2}-4f_{j-1}+3f_j}{2\Delta
x}\right)x+\left(\dfrac{f_{j-2}-2f_{j-1}+f_j}{2\Delta
x^2}\right)x^2,
\end{align}
\begin{align}\label{f1_up}
\hat{f}^1(x)=\dfrac{-f_{j-1}+26f_{j}-f_{j+1}}{24}+\left(\dfrac{f_{j+1}-f_{j-1}}{2\Delta
x}\right)x+\left(\dfrac{f_{j-1}-2f_{j}+f_{j+1}}{2\Delta
x^2}\right)x^2,
\end{align}
\begin{align}\label{f2_up}
\hat{f}^2(x)=\dfrac{23f_{j}+2f_{j+1}-f_{j+2}}{24}+\left(\dfrac{-3f_{j}+4f_{j+1}-f_{j+2}}{2\Delta
x}\right)x+\left(\dfrac{f_{j}-2f_{j+1}+f_{j+2}}{2\Delta
x^2}\right)x^2.
\end{align}
}

Evaluating each of these $\hat{f}^k(x)$'s at $x_{j+\frac{1}{2}}$, we
obtain that
\begin{align}\label{fh0}
\hat{f}^0_{j+\frac{1}{2}}=\dfrac{2}{6}f_{j-2}-\dfrac{7}{6}f_{j-1}+\dfrac{11}{6}f_{j},
\end{align}
\begin{align}\label{fh1}
\hat{f}^1_{j+\frac{1}{2}}=-\dfrac{1}{6}f_{j-1}+\dfrac{5}{6}f_{j}+\dfrac{2}{6}f_{j+1},
\end{align}
\begin{align}\label{fh2}
\hat{f}^2_{j+\frac{1}{2}}=\dfrac{2}{6}f_{j}+\dfrac{5}{6}f_{j+1}-\dfrac{1}{6}f_{j+2}.
\end{align}

Carrying a similar process as given in Eqs. (\ref{fh_5_tay}) -
(\ref{fp5_jH}) with $\hat{f}^k_{j+\frac{1}{2}}$ replacing
$\hat{f}^5_{j+\frac{1}{2}}$, we obtain that
\begin{align}\label{fh_3_tay}
h_{j+\frac{1}{2}}=\hat{f}^k_{j+\frac{1}{2}}+\mathcal{O}(\Delta x^3).
\end{align}

Comparing between $\hat{f}^5_{j+\frac{1}{2}}$ given in Eq.
(\ref{fh_5}) and the ones in Eqs. (\ref{fh0}) - (\ref{fh2}), we
deduce the following linear relation
\begin{align}\label{fh_linear}
\hat{f}^5_{j+\frac{1}{2}}=\sum_{k=0}^2\gamma_k\hat{f}^k_{j+\frac{1}{2}},
\end{align}
where
\begin{align}\label{gam}
\gamma_0=\frac{1}{10},\quad\gamma_1=\frac{6}{10},\quad\gamma_2=\frac{3}{10},
\end{align}
are called the linear (optimal) weights. We note that
\begin{align}\label{gam_norm}
\sum_{k=0}^2\gamma_k=1.
\end{align}

Adding and subtracting
$\sum_{k=0}^2\gamma_k\hat{f}^k_{j+\frac{1}{2}}$ into and from Eq.
(\ref{weno}), thanks to the truncation errors in Eqs.
(\ref{fh_5_tay}), (\ref{fh_3_tay}), the normalization in Eqs.
(\ref{om_norm}), (\ref{gam_norm}), and the linear relation
(\ref{fh_linear}) we obtain that
\begin{align}
\begin{split}
\hat{f}_{j\pm\frac{1}{2}}&=\sum_{k=0}^2(\omega^{\pm}_k-\gamma^{\pm}_k)\hat{f}^k_{j^{\pm}\frac{1}{2}}
+\sum_{k=0}^2\gamma^{\pm}_k\hat{f}^k_{j^{\pm}\frac{1}{2}}\\
&=\sum_{k=0}^2(\omega^{\pm}_k-\gamma^{\pm}_k)(h_{j^{\pm}\frac{1}{2}}+\mathcal{O}(\Delta
x^3))+\left(h_{j^{\pm}\frac{1}{2}}-\dfrac{1}{60}\dfrac{\partial^5
f}{\partial
x^5}\bigg|_{x=x_j}\Delta x^5+\mathcal{O}(\Delta x^6)\right)\\
&=h_{j^{\pm}\frac{1}{2}}-\dfrac{1}{60}\dfrac{\partial^5 f}{\partial
x^5}\bigg|_{x=x_j}\Delta
x^5+\sum_{k=0}^2(\omega^{\pm}_k-\gamma^{\pm}_k)\mathcal{O}(\Delta
x^3)+\mathcal{O}(\Delta x^6),
\end{split}
\end{align}
where $\gamma^{\pm}_k$, $\omega^{\pm}_k$ are the linear and
non-linear weights of the sub-stencils $S^{j\pm\frac{1}{2}}_k$
corresponding to the interfaces $x_{j\pm\frac{1}{2}}$, respectively.

Hence, in order that the discretization in Eq. (\ref{approx_dis})
where $\hat{f}^5_{j\pm\frac{1}{2}}=\hat{f}_{j\pm\frac{1}{2}}$ follow
the nonlinear relation (\ref{weno}) to be 5th-order, we deduce the
sufficient condition for the nonlinear weights as follows,
\begin{align}\label{cond_2}
\omega^{\pm}_k-\gamma^{\pm}_k=\mathcal{O}(\Delta x^3),\quad\forall
k.
\end{align}

Different WENO schemes depends on how these nonlinear weights and
the smoothness indicators are defined. The latter ones are
introduced in the below section. In the following subsections, we
summarize the 5th-order upwind and 6th-order central WENO schemes
discussed previously. Since WENO-Z is a good replacement for WENO-M,
we omit the latter one in our comparison.

\subsubsection{WENO-JS}

In \cite{JS}, Jiang and Shu defined the nonlinear weights as
follows,
\begin{align}\label{om_js}
\omega^{JS}_k=\dfrac{\alpha^{JS}_k}{\sum_{l=0}^2\alpha^{JS}_l},\quad
\alpha^{JS}_{k}=\dfrac{\gamma_k}{(\ep+\beta_k)^p},
\end{align}
where $\gamma_k$ is defined in Eq. (\ref{gam}), $\beta_k$ is called
the smoothness indicator of $S_k$ which measures how smooth the
solution is over this sub-stencil. The authors defined these
$\beta_k$'s through the normalized $L^2$-norm of high-order
variations of the reconstruction polynomials given in Eqs.
(\ref{fh0}) - (\ref{fh2}). Explicitly, for a 5th-order scheme, we
have
\begin{align}\label{be}
\beta_k=\Delta
x\int_{x_{j-\frac{1}{2}}}^{x_{j+\frac{1}{2}}}\left(\dfrac{d}{dx}\hat{f}^k(x)\right)^2dx+
\Delta
x^3\int_{x_{j-\frac{1}{2}}}^{x_{j+\frac{1}{2}}}\left(\dfrac{d^2}{dx^2}\hat{f}^k(x)\right)^2dx,
\end{align}
where $\hat{f}^k(x)$'s are as in Eqs. (\ref{f0_up}) - (\ref{f2_up})
and the sub-stencils $S_0$, $S_1$, $S_2$ are centered around
$x_j=0$.

Evaluating for each $k$, we obtain that
\begin{align}\label{be_0}
\beta_0&=\frac{13}{12}(f_{j-2}-2f_{j-1}+f_j)^2+\frac{1}{4}(f_{j-2}-4f_{j-1}+3f_{j})^2\\
&=f^{\prime2}\Delta
x^2+\left(\frac{13}{12}f^{\prime\prime2}-\frac{2}{3}f^{\prime}f^{\prime\prime\prime}\right)\Delta
x^4+\mathcal{O}(\Delta x^5),\nonumber
\end{align}
\begin{align}\label{be_1}
\beta_1&=\frac{13}{12}(f_{j-1}-2f_{j}+f_{j+1})^2+\frac{1}{4}(f_{j+1}-f_{j-1})^2\\
&=f^{\prime2}\Delta
x^2+\left(\frac{13}{12}f^{\prime\prime2}+\frac{1}{3}f^{\prime}f^{\prime\prime\prime}\right)\Delta
x^4+\mathcal{O}(\Delta x^6),\nonumber
\end{align}
\begin{align}\label{be_2}
\beta_2&=\frac{13}{12}(f_{j}-2f_{j+1}+f_{j+2})^2+\frac{1}{4}(3f_{j}-4f_{j+1}+f_{j+2})^2\\
&=f^{\prime2}\Delta
x^2+\left(\frac{13}{12}f^{\prime\prime2}-\frac{2}{3}f^{\prime}f^{\prime\prime\prime}\right)\Delta
x^4+\mathcal{O}(\Delta x^5),\nonumber
\end{align}
where the derivatives are evaluated at $x=x_j$.

In formula (\ref{om_js}), $\ep$ is a small parameter to prevent
division by zero. In most cases, WENO-JS works well with
$\ep=10^{-6}$. A thorough analysis of the role of $\ep$ can be found
in \cite{HAP}. The parameter $p$ is to increase the dissipation of
the scheme. For WENO-JS, $p=2$ is chosen.

If $f^{\prime}_j = f'(x_j) \ne0$, $\forall k$, $\beta_k$ can be written in the
form
\begin{align}\label{be_js}
\beta_k=(f^{\prime}_j\Delta x)^2(1+\mathcal{O}(\Delta x^2)).
\end{align}
Substituting these into Eq. (\ref{om_js}) with the removal of $\ep$,
since $(1+y)^{-2}=1+\mathcal{O}(y)$, by Eq. (\ref{gam_norm}) we
obtain that
\begin{align}\label{om_js_tay}
\begin{split}
\omega^{JS}_k=\dfrac{\gamma_k(f^{\prime}_j\Delta
x)^{-2}(1+\mathcal{O}(\Delta x^2))}{(f^{\prime}_j\Delta
x)^{-2}\sum_{l=0}^2\gamma_l(1+\mathcal{O}(\Delta
x^2))}=\gamma_k+\mathcal{O}(\Delta x^2),
\end{split}
\end{align}
which is a relaxed form of (\ref{cond_2}). We notice that for
WENO-JS, $\omega^{JS}_k$ cannot satisfy condition (\ref{cond_2})
directly. Moreover, if $f^{\prime}_j=0$, it is observed from Eqs.
(\ref{be_0}) - (\ref{be_2}) that
$\beta_k=\frac{13}{12}(f^{\prime\prime}_j)^2\Delta
x^4(1+\mathcal{O}(\Delta x))$ for $k=0,2$ and
$\beta_1=\frac{13}{12}(f^{\prime\prime}_j)^2\Delta
x^4(1+\mathcal{O}(\Delta x^2))$, in which the condition
(\ref{be_js}) is not satisfied for all $k$'s. Similarly, we find
that
\begin{align}
\omega^{JS}_k=\gamma_k+\mathcal{O}(\Delta x),
\end{align}
which is a loss of accuracy near this critical point. This accuracy
loss is improved by the WENO-Z scheme which is summarized in the next subsection.

\subsubsection{WENO-Z}

Borges et al. in \cite{BCCD} proposed a new WENO-Z. In their scheme,
the nonlinear weights are defined differently from those of WENO-JS.
They are as follows,
\begin{align}\label{om_z}
\omega^Z_k=\dfrac{\alpha^Z_k}{\sum_{l=0}^2\alpha^Z_l},\quad
\alpha^Z_{k}=\gamma_k\left(1+\left(\dfrac{\tau^Z}{\ep+\beta_k}\right)^q\right),
\end{align}
where the smoothness indicators $\beta_k$'s are the same as those
given in Eqs. (\ref{be_0}) - (\ref{be_2}), $\ep=10^{-40}$, and
$\tau^Z$ is the smoothness indicator of the large stencil $S^5$. The
power $q$ is used to tune the relation between the dispersive and
dissipative property of the scheme. It is checked numerically in
\cite{BCCD} that the scheme becomes more dissipative when $q$ is
increased. For WENO-Z, $\tau^Z$ is defined as follows,
\begin{align}\label{tau_z}
\tau^Z=|\beta_0-\beta_2|=\frac{13}{3}|f^{\prime\prime}f^{\prime\prime\prime}|\Delta
x^5+\mathcal{O}(\Delta x^6).
\end{align}

We note that if $f^{\prime}_j\ne0$ and $q=1$,
$\beta_k=\mathcal{O}(\Delta x^2)$, $\forall k$. Then
\begin{align}\label{tau_be_z}
\dfrac{\tau^Z}{\beta_k}=\mathcal{O}(\Delta x^3),\quad\forall k.
\end{align}

Similarly to Eq. (\ref{om_js_tay}), we obtain that
\begin{align}\label{om_z_tay}
\omega^Z_k=\gamma_k+\mathcal{O}(\Delta x^3),
\end{align}
directly without using the relaxed version as the WENO-JS scheme. It
was also proven in \cite{BCCD} that WENO-Z is 4th-order near simple
smooth critical points (i.e. where $f'_j=0$) for $q=1$ and attains
the designed 5th-order for $q=2$. The tradeoff for the latter case
is that the scheme is more dissipative. Throughout this work, we
choose $q=1$. One more advantage of WENO-Z over WENO-JS is that the
former is more central in a sense that the stencil over which the
solution is discontinuous plays more roles in the approximation of
the numerical flux. This assessment is checked as follows. We
suppose that $S_2$ contains a discontinuity whereas the solution is
smooth over the other two sub-stencils. Hence,
$\tau^Z=\mathcal{O}(1)$, $\beta_2=\mathcal{O}(1)$, and
$\beta_k=\mathcal{O}(\Delta x^2)$, $k=0,1$. Then,
\begin{align}
\begin{split}
\dfrac{\alpha^Z_2}{\alpha^Z_k}=\dfrac{\gamma_2\left(1+\frac{\tau^Z}{\beta_2}\right)}
{\gamma_k\left(1+\frac{\tau^Z}{\beta_k}\right)}
=\dfrac{\frac{\gamma_2}{\beta_2}(\beta_2+\tau^Z)}{\frac{\gamma_k}{\beta_k}(\beta_k+\tau^Z)}
=\dfrac{\alpha^{JS}_2(\beta_2+\tau^Z)}{\alpha^{JS}_k(\beta_k+\tau^Z)}
\ge\dfrac{\alpha^{JS}_2}{\alpha^{JS}_k},
\end{split}
\end{align}
since
\begin{align}
\dfrac{\beta_2+\tau^Z}{\beta_k+\tau^Z}\approx\dfrac{\beta_2+\tau^Z}{\tau^Z}\ge1.
\end{align}

Hence, WENO-Z has a sharper capturing of discontinuities than
WENO-JS.

\begin{rem}\label{rmk_flux_split}
\begin{enumerate}
\item[]
\item[i.] In case the condition $f^{\prime}(u)\ge0$ is not satisfied,
which is general in real applications, we apply a flux splitting
technique to decompose $f(u)$ into positive and negative components.
In most applications, the global Lax-Friedrichs flux splitting is
used (see \cite{JS} and the references therein),
\begin{align}\label{flux_split_lf}
f^{\pm}(u)=\dfrac{1}{2}(f(u)\pm\alpha u),
\end{align}
where $\alpha=\max|f^{\prime}(u)|$ over the whole computational
domain. Then,
\begin{align}
\hat{f}_{j+\frac{1}{2}}=\hat{f}^{+}_{j+\frac{1}{2}}+\hat{f}^{-}_{j+\frac{1}{2}}.
\end{align}
The negative flux $\hat{f}^{-}_{j+\frac{1}{2}}$ is symmetric to
$\hat{f}^{+}_{j+\frac{1}{2}}$ with respect to $x_{j+\frac{1}{2}}$.
\item[ii.] If the flux splitting is employed, the overall number of grid
points used in the reconstruction of the numerical flux is increased
by one. That is, let $S^{5+}$ and $S^{5-}$ be the stencils over
which $\hat{f}^{+}_{j+\frac{1}{2}}$ and
$\hat{f}^{-}_{j+\frac{1}{2}}$ are determined, then
\begin{align}
S^6:=S^{5+}\bigcup
S^{5-}=\{x_{j-2},x_{j-1},x_j,x_{j+1},x_{j+2},x_{j+3}\},
\end{align}
which consists of six points. We note that with this $S^6$, there
exists a polynomial of degree five which reconstructs $h(x)$ over
the stencil. Therefore, the accuracy of WENO schemes can be
increased up to sixth. These schemes are discussed in the subsection
below.
\end{enumerate}
\end{rem}

\subsection{6th-order Central WENO Reconstruction}

Carrying a similar procedure as described in the previous subsection
for the 6-point large stencil $S^6$, we can deduce that
\begin{align}\label{fh_6}
\begin{split}
\hat{f}^6_{j+\frac{1}{2}}&=\dfrac{1}{60}f_{j-2}-\dfrac{8}{60}f_{j-1}+\dfrac{37}{60}f_j+\dfrac{37}{60}f_{j+1}-\dfrac{8}{60}f_{j+2}+\dfrac{1}{60}f_{j+3}\\
&=h_{j+\frac{1}{2}}+\dfrac{1}{140}\dfrac{\partial^6f}{\partial
x^6}\bigg|_{x=x_j}\Delta x^6+\mathcal{O}(\Delta x^7).
\end{split}
\end{align}

Similarly,
\begin{align}
\begin{split}
\hat{f}^6_{j-\frac{1}{2}}&=\dfrac{1}{60}f_{j-3}-\dfrac{8}{60}f_{j-2}+\dfrac{37}{60}f_{j-1}
+\dfrac{37}{60}f_{j}-\dfrac{8}{60}f_{j+1}+\dfrac{1}{60}f_{j+2}\\
&=h_{j-\frac{1}{2}}+\dfrac{1}{140}\dfrac{\partial ^6f}{\partial
x^6}\bigg|_{x=x_j}\Delta x^6+\mathcal{O}(\Delta x^7).
\end{split}
\end{align}

Hence we obtain that
\begin{align}\label{approx_dis_6}
\dfrac{du_j}{dt}\approx-\dfrac{\hat{f}^6_{j+\frac{1}{2}}-\hat{f}^6_{j-\frac{1}{2}}}{\Delta
x}=:\hat{\mathcal{L}}^6(u).
\end{align}
The scheme is increased to 6th-order of accuracy in space.

\begin{figure}
  \begin{center}
    \begin{tabular}{c}
      \resizebox{105mm}{!}{\includegraphics{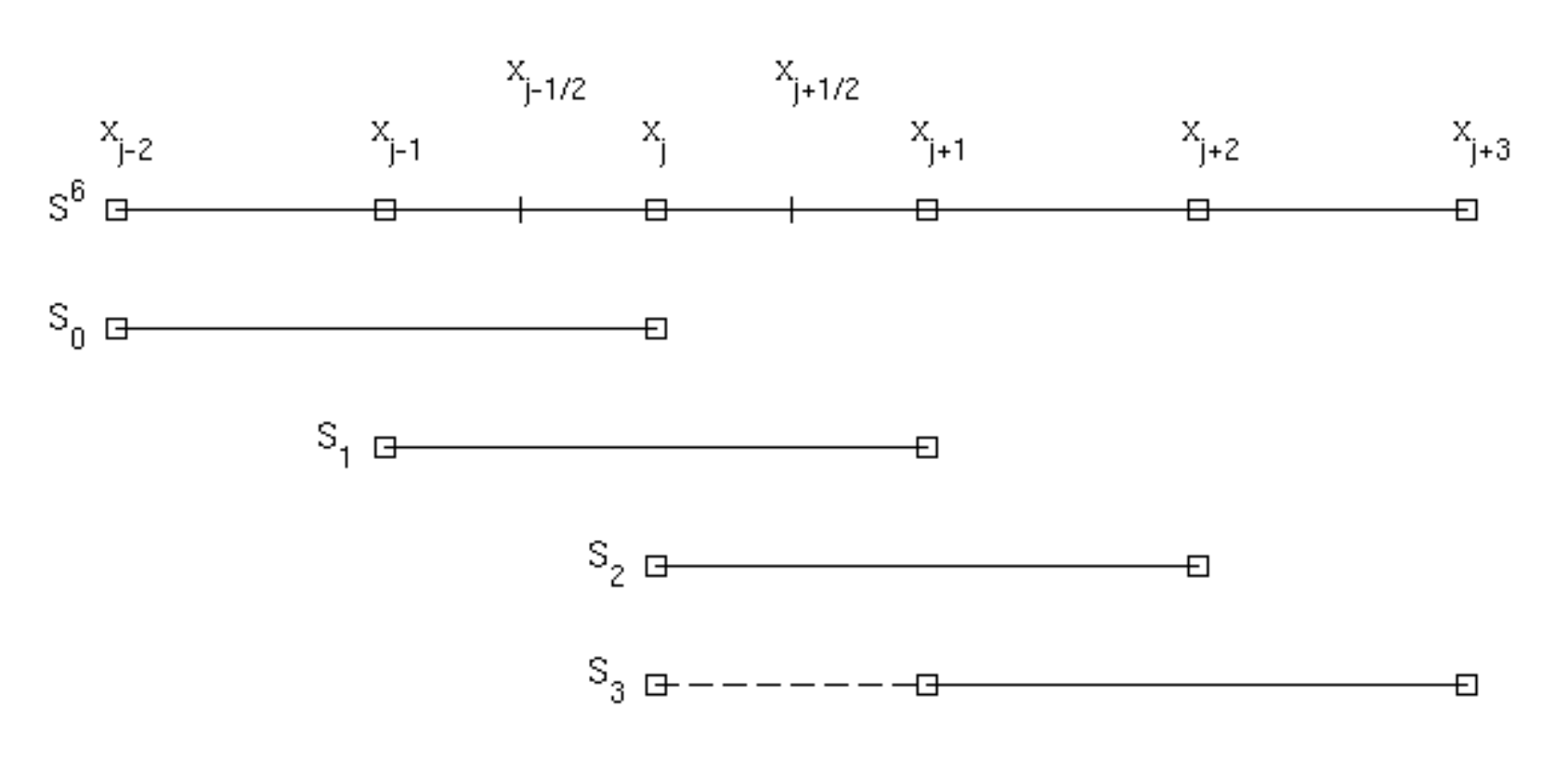}}
    \end{tabular}
    \caption{\small{Stencils for 6th-order WENO schemes.}}
\label{fig_sten_6}
  \end{center}
\end{figure}

Adding one more sub-stencil $S_3=\{x_{j+1},x_{j+2},x_{j+3}\}$ into
the approximation of the interfaced value $\hat{f}_{j+\frac{1}{2}}$
(see Fig. \ref{fig_sten_6}), we deduce a similar linear relation
with Eq. (\ref{fh_linear}) as follows
\begin{align}\label{fh6_lin}
\hat{f}^6_{j+\frac{1}{2}}=\sum_{k=0}^3\gamma_k\hat{f}^k_{j+\frac{1}{2}},
\end{align}
where
\begin{align}\label{gam_6}
\gamma_0=\gamma_3=\dfrac{1}{20},\quad\gamma_1=\gamma_2=\dfrac{9}{20}.
\end{align}

Here, $\hat{f}^3_{j+\frac{1}{2}}$ is the 3rd-order approximation of
the numerical flux $h(x)$ at the interface $x_{j+\frac{1}{2}}$ from
the reconstruction polynomial $\hat{f}^3(x)$ over the sub-stencil
$S_3$. Explicitly, we have
\begin{align}\label{fh3}
\hat{f}^3_{j+\frac{1}{2}}=\dfrac{11}{6}f_{j+1}-\dfrac{7}{6}f_{j+2}+\dfrac{2}{6}f_{j+3}.
\end{align}

The other approximations $\hat{f}^k_{j+\frac{1}{2}}$'s, $k=0,1,2$
follow Eqs. (\ref{fh0}) - (\ref{fh2}).

The nonlinear combination using the nonlinear weights is similar to
that of the linear case (\ref{fh6_lin}), except for the linear
weight $\gamma_k$ replaced by its nonlinear version $\omega_k$,
$k=0,\ldots,3$.

\begin{rem}\label{rmk_cen}
\begin{enumerate}
\item[]
\item[i.] Since the accuracy order is increased to sixth, the
sufficient condition for the nonlinear weights given in
(\ref{cond_2}) is also increased by one. That is,
\begin{align}\label{cond6_2}
\omega^{\pm}_k-\gamma^{\pm}_k=\mathcal{O}(\Delta x^4),\ \forall k.
\end{align}
\item[ii. ] Observing from Eq. (\ref{fh6_lin}) that the
approximations $\hat{f}_{j+\frac{1}{2}}$'s are now symmetric with
respect to $x_{j+\frac{1}{2}}$. It means that the scheme now is
central. Hence, spurious oscillations are expected to occur near
discontinuities. A treatment on the most downwind sub-stencil $S_3$
is needed to sustain the essentially non-oscillatory property of the
scheme. We now overview the 6th-order WENO schemes for this case.
\end{enumerate}
\end{rem}

\subsubsection{WENO-NW6}

In \cite{YC}, Yamaleev and Carpenter proposed a 6th-order
energy-stable WENO scheme. They introduced an artificial dissipative
term and proved that this makes the new scheme be stable in $L^2$
sense. In this work, we only discuss their treatment on the
nonlinear weights and omit this artificial dissipative term (see
\cite{YC} for a detailed discussion on the term). Hence, the scheme
here is denoted by WENO-NW6, not ESWENO as in their paper.

The nonlinear weights (denoted $\omega_k^{NW}$, $\alpha_k^{NW}$)
follow those defined in the WENO-Z scheme which are given in Eq.
(\ref{om_z}), for $k=0,\ldots,3$. The differences lie on the
smoothness indicator of the most downwind sub-stencil $\beta_3$ and
the one for the large stencil $S^6$. For the former, in order that
there are no oscillations occurring near discontinuities, all grid
values of the flux over $S^6$ are accounted for the computation of
$\beta_3$. It is as follows,
\begin{align}\label{be_3_es}
\beta_3=\dfrac{1}{4}(\beta_0^4+\beta_1^4+\beta_2^4+\tilde{\beta}_3^4)^{1/4},
\end{align}
where $\tilde{\beta}_3$ is computed using the formula given in Eq.
(\ref{be}). That is,
\begin{align}\label{be_3ti}
\begin{split}
\tilde{\beta}_3&=\dfrac{13}{12}(f_{j+1}-2f_{j+2}+f_{j+3})^2+\dfrac{1}{4}(-5f_{j+1}+8f_{j+2}-3f_{j+3})^2\\
&=f^{\prime2}\Delta
x^2+\left(\frac{13}{12}f^{\prime\prime2}-\frac{11}{3}f^{\prime}f^{\prime\prime\prime}\right)\Delta
x^4+\mathcal{O}(\Delta x^5),
\end{split}
\end{align}
where the derivatives are evaluated at $x=x_j$. The other indicators
$\beta_k$'s, $k=0,1,2$ follow Eqs. (\ref{be_0}) - (\ref{be_2}).

The smoothness indicator of the large stencil $S^6$ is defined as
the highest, i.e., fifth-degree, undivided difference as follows
\begin{align}\label{tau_es}
\begin{split}
\tau^{NW}&=(f_{j-2}-5f_{j-1}+10f_{j}-10f_{j+1}+5f_{j+2}-f_{j+3})^2\\
&=(f^{(5)})^2\Delta x^{10}+\mathcal{O}(\Delta x^{11}).
\end{split}
\end{align}

Carrying a similar procedure as in Eqs. (\ref{tau_be_z}) -
(\ref{om_z_tay}) with a change in $\tau$, thanks to Eq.
(\ref{tau_es}), we obtain that
\begin{align}\label{tau_be_nw}
\dfrac{\tau^{NW}}{\beta_k}=
\begin{cases}
&\mathcal{O}(\Delta x^8),\quad\text{if }f^{\prime}_j\ne0,\\
&\mathcal{O}(\Delta x^6),\quad\text{if }f^{\prime}_j=0,
\end{cases}
\quad\forall k,
\end{align}
thus condition (\ref{cond6_2}) is satisfied. Hence, the scheme is
6th-order in smooth regions. The case where the derivatives of $f$
vanish will be checked numerically in the below section.

\subsubsection{WENO-CU6}

In \cite{HWA}, Hu et al. developed the adaptive central-upwind
scheme WENO-CU6 based on the principle that the most downwind
sub-stencil only plays roles in smooth regions and is suppressed
near discontinuities. Hence the scheme is central in smooth regions
and upwind near discontinuities. The scheme is different from
WENO-NW6 in defining the smoothness indicators $\beta_3$ and $\tau$.
In particular, they are as below, {\small
\begin{align}\label{be_3_cu}
\begin{split}
\beta_3&=\dfrac{1}{120960}[271779f_{j-2}^2+f_{j-2}(-2380800f_{j-1}+4086352f_j-3462252f_{j+1}+1458762f_{j+2}-245620f_{j+3})\\
&+f_{j-1}(5653317f_{j-1}-20427884f_j+17905032f_{j+1}-7727988f_{j+2}+1325006f_{j+3})+f_{j}(19510972f_j\\
&-35817664f_{j+1}+15929912f_{j+2}-2792660f_{j+3})+f_{j+1}(17195652f_{j+1}-15880404f_{j+2}\\
&+2863984f_{j+3})+f_{j+2}(3824847f_{j+2}-1429976f_{j+3})+139633f_{j+3}^2]\\
&=f^{\prime2}\Delta x^2+\dfrac{13}{12}f^{\prime\prime2}\Delta
x^4+\mathcal{O}(\Delta x^6).
\end{split}
\end{align}
\emph{It is noticed that there is a typo in Eq. (25) in \cite{HWA},
and we have corrected it in Eq. (\ref{be_3_cu}).}

From there, the smoothness indicator of the large stencil $S^6$ is
defined as follows,
\begin{align}\label{tau_cu}
\tau^{CU}=\beta_3-\dfrac{1}{6}(\beta_0+4\beta_1+\beta_2)=\mathcal{O}(\Delta
x^6).
\end{align}

Hence, we have for $k=0,\ldots,3,$,
\begin{align}
\dfrac{\tau^{CU}}{\beta_k}=
\begin{cases}
&\mathcal{O}(\Delta x^4),\quad\text{if }f^{\prime}_j\ne0,\\
&\mathcal{O}(\Delta x^2),\quad\text{if }f^{\prime}_j=0,
\end{cases}
\end{align}
which satisfies the condition (\ref{cond6_2}).

It is also noteworthy that $\alpha^{CU}_k$ in Eq. (\ref{om_z}) has a
change as below,
\begin{align}
\alpha^{CU}_k=\gamma_k\left(C+\dfrac{\tau^{CU}}{\ep+\beta_k}\right),
\end{align}
where $C\gg1$ is to increase the contribution of the linear weights
when the smoothness indicators have comparable magnitudes (see
\cite{TWM}). Following \cite{HWA}, we choose $C=20$.

As indicated in example \ref{exam_1}, the 6th-order WENO-NW6 and
WENO-CU6 schemes suffer from the loss of accuracy near the smooth
critical points just right behind, with respect to the
characteristic direction, a critical point where the first
derivative of the solution is undefined, that is, the solution is
just $C^0$ at that point. The explanation for this defect is given
in the below section. In the next section, we propose a new scheme
which automatically switches between a 6th-order central and
5th-order upwind scheme and overcomes the defect occurred in the
mentioned schemes.

\section{The New Scheme}

We first observe that the 5th- and 6th-order linear approximations
given in Eqs. (\ref{fh_5}) and (\ref{fh_6}), respectively, can be
combined linearly in the following manner,
\begin{align}\label{fh_we6}
\hat{f}_{j+\frac{1}{2}}=\gamma^{\theta}_0\hat{f}^0_{j+\frac{1}{2}}+\gamma^{\theta}_1\hat{f}^1_{j+\frac{1}{2}}
+\gamma^{\theta}_2\hat{f}^2_{j+\frac{1}{2}}+\gamma^{\theta}_3\hat{f}^3_{j+\frac{1}{2}},
\end{align}
where
\begin{align}\label{gam_we6}
\gamma^{\theta}_0=\frac{1}{20}(1+\theta),\quad\gamma^{\theta}_1=\frac{3}{20}(3+\theta),
\quad\gamma^{\theta}_2=\frac{3}{20}(3-\theta),\quad\gamma^{\theta}_3=\frac{1}{20}(1-\theta),
\end{align}
and $\hat{f}^k_{j+\frac{1}{2}}$, $k=0,1,2,3$ is given in Eqs.
(\ref{fh0}) - (\ref{fh2}) and (\ref{fh3}), respectively.

We deduce that
\begin{align}
\hat{f}_{j+\frac{1}{2}}=
\begin{cases}
& \hat{f}^5_{j+\frac{1}{2}}\quad\text{if }\theta=1,\\
& \hat{f}^6_{j+\frac{1}{2}}\quad\text{if }\theta=0.
\end{cases}
\end{align}

We then propose a new scheme in which $\hat{f}_{j+\frac{1}{2}}$ is
chosen between $\hat{f}^5_{j+\frac{1}{2}}$ and
$\hat{f}^6_{j+\frac{1}{2}}$ adaptively. Hence, the scheme is
$5$th-order upwind or $6$th-order central depending on the
smoothness of the stencils $S^5$ and $S^6$. We expect that this will
get over the drawback of accuracy degeneration of the above
mentioned central 6th-order schemes. To proceed, we first rewrite
Eq. (\ref{fh_we6}) using instead the non-linear weights
$\omega^{\theta}_k$'s as follows,
\begin{align}\label{fh_nonlinear}
\hat{f}_{j+\frac{1}{2}}=\omega^{\theta}_0\hat{f}^0_{j+\frac{1}{2}}+\omega^{\theta}_1\hat{f}^1_{j+\frac{1}{2}}
+\omega^{\theta}_2\hat{f}^2_{j+\frac{1}{2}}+\omega^{\theta}_3\hat{f}^3_{j+\frac{1}{2}},
\end{align}
where, for $k=0,\ldots,3$, and
\begin{align}\label{om}
\omega^{\theta}_k=\dfrac{\alpha^{\theta}_k}{\sum_{l=0}^3\alpha^{\theta}_l},\quad
\alpha^{\theta}_k=\gamma^{\theta}_k\left(1+\dfrac{\tau^{\theta}}{\ep+\tilde{\beta}_k}\right).
\end{align}
Here, $\tau^{\theta}$ is the smoothness indicator of the large
stencil, and $\tilde{\beta}_k$ is the smoothness indicator of the
sub-stencil $S_k$. We define these indicators in the following
subsection.

\subsection{The Central Smoothness Indicators $\tilde{\beta}_k$ and New $\tau^{\theta}$}

For a 6th-order central scheme over the large stencil $S^6$,
spurious oscillations are expected to occur near discontinuities. To
overcome this, WENO-NW6 and WENO-CU6 choose to construct their
$\beta_3$ over all points of $S^6$. A more careful observation
reveals that this cost can be reduced in the following way. We
remind the principle of WENO schemes is that there is at least one
smoothest sub-stencil is used in the reconstruction of the numerical
flux. We suppose that $\beta_3$ follows Eq. (\ref{be_3ti}), that is,
it measures the smoothness of the most downwind $S_3$ only. We
further assume that the grid size $\Delta x$ is so small that a
discontinuity does not spread over two neighboring grid points, then
for a 6th-order WENO scheme, the only case where oscillations occur
is when a discontinuity is in between $x_j$ and $x_{j+1}$. In this
case, both $\hat{f}^0_{j+\frac{1}{2}}$ and
$\hat{f}^3_{j+\frac{1}{2}}$ play main roles in the combination
(\ref{fh_nonlinear}) since $\beta_0$ and $\beta_3$ are much smaller
than the other two. This leads to oscillations since the downwind
$\hat{f}^3_{j+\frac{1}{2}}$ is wrongly chosen. To prevent this from
happening, we choose $\beta_3$ to measure the smoothness of an
extended sub-stencil $\tilde{S}_3:=\{x_j,x_{j+1},x_{j+2},x_{j+3}\}$
(see Fig. \ref{fig_sten_6} with $S_3$ extended by the dashed line).
It is observed that $S_2$ is now a subset of $\tilde{S}_3$ and all
sub-stencils share the point $x_j$. Hence, the case where
oscillations occur is essentially eliminated. Moreover, the cost of
computing $\beta_3$ is now much reduced since the computation
involves only four grid points defined in $\tilde{S}_3$ instead of
six in $S^6$ as in WENO-NW6 or WENO-CU6.

For 5th-order schemes, all sub-stencils are symmetric with respect
to $x_j$. As a result, the smoothness indicators $\beta_k$'s are
also symmetric with respect to $x_j$ (see Eqs. (\ref{be_0}) -
(\ref{be_2})). That is, $\beta_0$ and $\beta_2$ are equal to each
other up to order $\mathcal{O}(\Delta x^4)$ in Taylor expansions. We
recall that WENO discretizations choose the sub-stencils depending
on the non-linear weights $\omega_k$'s which are very sensitive to
the smoothness indicators $\beta_k$'s due to the latter's smallness
in smooth regions (see Eq. (\ref{om_js}) for WENO-JS, Eq.
(\ref{om_z}) for WENO-Z, WENO-NW6, and WENO-CU6). For the
sensitivity, we mean that a small change in any $\beta_k$ leads to a
large difference among $\alpha_k$'s, thus $\omega_k$'s. In that
sense, the symmetry in terms of Taylor expansions of $\beta_k$'s
reduces the effects of this sensitivity, especially in transition
regions where the solution is smooth and discontinuous. We refer to
Figs. \ref{fig_soln_ini3} and \ref{fig_soln_lin} below for numerical
evidences for this assessment in which the schemes with symmetric
$\beta_k$'s (i.e., WENO-Z and WENO-$\theta$) show better results
than the ones without this property. Unfortunately, the 6th-order
methods lack of this (comparing $\beta_3$ in Eq. (\ref{be_3_es}) for
WENO-NW6, and Eq. (\ref{be_3_cu}) for WENO-CU6 with the other
$\beta_k$'s, $k=0,1,2$, defined Eqs. (\ref{be_0}) - (\ref{be_2})).
In our new scheme, we try to recover the property. We devise our new
smoothness indicators in a central sense. That is, they are
constructed based on the reconstruction polynomials which are
symmetric with respect to $x_{j+\frac{1}{2}}$. In addition, it is
shown below that the new indicators are symmetric in terms of Taylor
expansions with respect to $x_j$. We notice that Taylor expansions
about $x_j$ are natural since the approximation of $f(u)_x$ is at
the interval center $x_j$ (see Eq. (\ref{mod_semi})) although the
reconstruction of the numerical flux function $h(x)$ is at the
interface $x_{j+\frac{1}{2}}$ (see Eqs. (\ref{fh_5}) and
(\ref{fh_6}) for 5th-order and 6th-order schemes, respectively).

Proceeding the reconstruction procedure as given in subsection
\ref{recon}, but instead around $x_{j+\frac{1}{2}}=0$ and with
$\tilde{S}_3$ replacing $S_3$, we obtain that {\small
\begin{align}\label{f0_cen}
\tilde{f}^0(x)=\dfrac{2f_{j-2}-7f_{j-1}+11f_{j}}{6}+\left(\dfrac{f_{j-2}-3f_{j-1}+2f_j}{\Delta
x}\right)x+\left(\dfrac{f_{j-2}-2f_{j-1}+f_j}{2\Delta
x^2}\right)x^2,
\end{align}
\begin{align}\label{f1_cen}
\tilde{f}^1(x)=\dfrac{-f_{j-1}+5f_{j}+2f_{j+1}}{6}+\left(\dfrac{f_{j+1}-f_{j}}{\Delta
x}\right)x+\left(\dfrac{f_{j-1}-2f_{j}+f_{j+1}}{2\Delta
x^2}\right)x^2,
\end{align}
\begin{align}\label{f2_cen}
\tilde{f}^2(x)=\dfrac{2f_{j}+5f_{j+1}-f_{j+2}}{6}+\left(\dfrac{f_{j+1}-f_{j}}{\Delta
x}\right)x+\left(\dfrac{f_{j}-2f_{j+1}+f_{j+2}}{2\Delta
x^2}\right)x^2,
\end{align}
\begin{align}\label{f3_cen}
\begin{split}
\tilde{f}^3(x)&=\dfrac{3f_{j}+13f_{j+1}-5f_{j+2}+f_{j+3}}{12}+\left(\dfrac{-11f_j+9f_{j+1}+3f_{j+2}-f_{j+3}}{12\Delta
x}\right)x\\
&+\left(\dfrac{3f_j-7f_{j+1}+5f_{j+2}-f_{j+3}}{4\Delta
x^2}\right)x^2+\left(\dfrac{-f_j+3f_{j+1}-3f_{j+2}+f_{j+3}}{6\Delta
x^3}\right)x^3.
\end{split}
\end{align}
}

We notice that these reconstruction polynomials are different from
those given in Eqs. (\ref{f0_up}) - (\ref{f2_up}) which are
constructed symmetrically with respect to $x_j=0$. Substituting
these into Eq. (\ref{be}) with $\tilde{f}^k(x)$ replacing
$\hat{f}^k(x)$, $k=0,1,2,3$, we deduce the new central smoothness
indicators as follows,
\begin{align}\label{be0_cen}
\begin{split}
\tilde{\beta}_0&=\frac{13}{12}(f_{j-2}-2f_{j-1}+f_j)^2+(f_{j-2}-3f_{j-1}+2f_{j})^2\\
&=f^{\prime2}\Delta x^2+f^{\prime}f^{\prime\prime}\Delta
x^3+\left(\frac{4}{3}f^{\prime\prime2}-\frac{5}{3}f^{\prime}f^{\prime\prime\prime}\right)\Delta
x^4+\mathcal{O}(\Delta x^5),
\end{split}
\end{align}
\begin{align}\label{be1_cen}
\begin{split}
\tilde{\beta}_1&=\frac{13}{12}(f_{j-1}-2f_{j}+f_{j+1})^2+(f_{j+1}-f_{j})^2\\
&=f^{\prime2}\Delta x^2+f^{\prime}f^{\prime\prime}\Delta
x^3+\left(\frac{4}{3}f^{\prime\prime2}+\frac{1}{3}f^{\prime}f^{\prime\prime\prime}\right)\Delta
x^4+\mathcal{O}(\Delta x^5),
\end{split}
\end{align}
\begin{align}\label{be2_cen}
\begin{split}
\tilde{\beta}_2&=\frac{13}{12}(f_{j}-2f_{j+1}+f_{j+2})^2+(f_{j}-f_{j+1})^2\\
&=f^{\prime2}\Delta x^2+f^{\prime}f^{\prime\prime}\Delta
x^3+\left(\frac{4}{3}f^{\prime\prime2}+\frac{1}{3}f^{\prime}f^{\prime\prime\prime}\right)\Delta
x^4+\mathcal{O}(\Delta x^5),
\end{split}
\end{align}
and
\begin{align}\label{be3_cen}
\begin{split}
\tilde{\beta}_3&=\frac{13}{48}(3f_{j}-7f_{j+1}+5f_{j+2}-f_{j+3})^2+(2f_{j+1}-3f_{j+2}+f_{j+3})^2\\
&=f^{\prime2}\Delta x^2+f^{\prime}f^{\prime\prime}\Delta
x^3+\left(\frac{4}{3}f^{\prime\prime2}-\frac{5}{3}f^{\prime}f^{\prime\prime\prime}\right)\Delta
x^4+\mathcal{O}(\Delta x^5),
\end{split}
\end{align}
where the derivatives are evaluated at $x=x_j$. We note that for the
most downwind $\tilde{\beta}_3$, we treated $\tilde{f}_3(x)$ in Eq.
(\ref{f3_cen}) as a 2nd-degree polynomial by ignoring the 3rd-degree
term when substituting it into Eq. (\ref{be}) so that
$\int_{x_{j-\frac{1}{2}}}^{x_{j+\frac{1}{2}}}\left(\frac{d^2}{dx^2}\tilde{f}^3(x)\right)^2dx$
is the highest-order variation. This is for the consistency with the
other reconstruction polynomials $\tilde{f}^k$'s, $k=0,1,2$, which
are of lower-degree. We also modified the second term in the
obtained indicator so that its Taylor expansion agrees with that of
$\tilde{\beta}_0$ up to order $\mathcal{O}(\Delta x^4)$; thus all
$\tilde{\beta}_k$'s are now symmetric with respect to $x_j$ in
Taylor expansions, which is our goal in designing these new
smoothness indicators. The original smoothness indicator $\hat{\beta}_3$ for
$\tilde{\beta}_3$ obtained from Eq. (\ref{be}) was as below,
\begin{align}
\hat{\beta}_3=\frac{13}{48}(3f_{j}-7f_{j+1}+5f_{j+2}-f_{j+3})^2
+\frac{1}{144}(-11f_j+9f_{j+1}+3f_{j+2}-f_{j+3})^2.
\end{align}

In order to enhance the dispersion of the scheme, following the
approach by Taylor et al. (\cite{TWM}), we set a restriction on the
smoothness indicators as below, for $k=0,\ldots,3$,
\begin{align}\label{be_res}
\tilde{\beta}_k=
\begin{cases}&0,\quad\text{if }R(\tilde{\beta})\le\alpha_R,\\
&\tilde{\beta}_k,\quad\text{otherwise};
\end{cases}
\end{align}
where
\begin{align}
R(\tilde{\beta})=\dfrac{\max_k(\tilde{\beta}_k)}{\ep+\min_k(\tilde{\beta}_k)}.
\end{align}
Here, $\alpha_R$ is a threshold value depending on the
configurations of flows. $\alpha_R$ is taken small for flows with
the presence of shocks. For a detailed discussion, consult
\cite{TWM}.

We next devise the smoothness indicator of the large stencil $S^6$.
Since the one proposed by Yamaleev and Carpenter in \cite{YC} (see
Eq. (\ref{tau_es})) is too dispersive and may lead to oscillations
(see the evidence in \cite{YC}), we introduce a new smoothness
indicator $\tau^{\theta}$ which is based on Eq. (\ref{be}) but at a
much higher order variations for the large stencil $S^6$. We
consider the following $\tau$'s,
\begin{align}\label{tau_5}
\begin{split}
\tau_5&=\Delta
x^5\int_{x_{j-\frac{1}{2}}}^{x_{j+\frac{1}{2}}}\left(\dfrac{d^3}{dx^3}\tilde{f}^5(x)\right)^2dx+\Delta
x^7\int_{x_{j-\frac{1}{2}}}^{x_{j+\frac{1}{2}}}\left(\dfrac{d^4}{dx^4}\tilde{f}^5(x)\right)^2dx\\
&=\dfrac{13}{12}(f_{j-2}-4f_{j-1}+6f_{j}-4f_{j+1}+f_{j+2})^2+(-f_{j-1}+3f_{j}-3f_{j+1}+f_{j+2})^2\\
&=f^{\prime\prime\prime2}\Delta
x^6+f^{\prime\prime\prime}f^{(4)}\Delta
x^7+\left(\dfrac{1}{2}f^{\prime\prime\prime}f^{(5)}+\dfrac{4}{3}(f^{(4)})^2\right)\Delta
x^8+\mathcal{O}(\Delta x^9),
\end{split}
\end{align}
and
\begin{align}\label{tau_6}
\begin{split}
\tau_6&=\Delta
x^7\int_{x_{j-\frac{1}{2}}}^{x_{j+\frac{1}{2}}}\left(\dfrac{d^4}{dx^4}\tilde{f}^6(x)\right)^2dx+\Delta
x^9\int_{x_{j-\frac{1}{2}}}^{x_{j+\frac{1}{2}}}\left(\dfrac{d^5}{dx^5}\tilde{f}^6(x)\right)^2dx\\
&=\dfrac{13}{12}(-f_{j-2}+5f_{j-1}-10f_{j}+10f_{j+1}-5f_{j+2}+f_{j+3})^2\\
&+\dfrac{1}{4}(f_{j-2}-3f_{j-1}+2f_{j}+2f_{j+1}-3f_{j+2}+f_{j+3})^2\\
&=(f^{(4)})^2\Delta x^8+f^{(4)}f^{(5)}\Delta
x^9+\left(\dfrac{5}{6}f^{(4)}f^{(5)}+\dfrac{4}{3}(f^{(5)})^2\right)\Delta
x^{10}+\mathcal{O}(\Delta x^{11}),
\end{split}
\end{align}
where $\tilde{f}^5(x)$ and $\tilde{f}^6(x)$ are the central
reconstruction polynomials around $x_{j+\frac{1}{2}}=0$ constructed
in a similar way as with $\tilde{f}^k(x)$'s in Eqs. (\ref{f0_cen}) -
(\ref{f3_cen}) but for the large stencils $S^5$ and $S^6$,
respectively; and the derivatives are evaluated at $x_j$.

We then choose our $\tau^{\theta}$ and set $\theta$ in Eq.
(\ref{gam_we6}) as follows,
\begin{align}\label{tau_we6}
(\tau^{\theta},\theta)=
\begin{cases}
&(\tau_6,0)\quad\text{if }\tau_6<\tau_5,\\
&(\tau_5,1)\quad\text{if }\tau_6\ge\tau_5.
\end{cases}
\end{align}

It is noted that by choosing such $\tau^{\theta}$ and $\theta$ as in
Eq. (\ref{tau_we6}), the scheme achieves a 6th-order in smooth
regions since $\tau_6\ll\tau_5$; whereas it adaptively chooses the
smoother large stencil between $S^5$ and $S^6$ in the WENO
reconstruction near discontinuities or unresolved regions. The new
scheme now chooses the smoothest not only sub-stencils but also
large one in the reconstruction procedure. The non-linear weights
follow Eq. (\ref{om}) above. Since our new method depends on
$\theta$ to switch between a 5th-order upwind and 6th-order central
scheme, we name it WENO-$\theta$. In the numerical simulations
below, we use the name WENO-$\theta6$ for the compatibility with the
other 6th-order schemes.

\begin{rem}
\begin{enumerate}
\item[]
\item[i.] Although the definition of $\tau^{\theta}$ has a switching
mechanism by an \emph{if} statement, it does not ruin the
methodology of WENO schemes. This is because the switching applies
to the smoothness indicator of the large stencil, not to the choice
of smoother sub-stencils. 
\item[ii.] Although the switching is discontinuous in nature,
WENO-$\theta$ is robust for problems with highly unstable fluid
flows. We illustrate this by conducting a numerical simulation of
the Rayleigh-Taylor instability problem in the below section.
\item[iii.] The role of $\ep$ is well investigated in \cite{HAP}. We
note here that for cases with increasing number of vanishing
derivatives, since both $\tau$ and $\beta_k$ are very small at
critical points, $\ep$ does play roles to sustain the designed
formal accuracy order. For this reason, except for WENO-JS where
$\ep=10^{-6}$, we choose $\ep=10^{-10}$ for other schemes. See the
accuracy tests in the below section.
\end{enumerate}
\end{rem}

In the next step, we test the accuracy, and resolutions of the new
scheme.

\subsection{Accuracy Tests}

We note that either $\tau_5$ or $\tau_6$ is chosen in Eq.
(\ref{tau_we6}), the sufficient condition (\ref{cond6_2}) is always
satisfied. Hence the new scheme is 6th-order in smooth regions.

For the tests of accuracy, we choose the linear scalar conservation
law,
\begin{align}\label{lin_eq}
\begin{cases}
&u_t+u_x=0,\quad x\in(-1,1),\\
&u(x,0)=u_0(x),
\end{cases}
\end{align}
subject to periodic boundary conditions. The following initial data
are considered:

$\bullet$ \emph{Initial condition 1:}
\begin{align}\label{sca_ini_1}
u_0(x)=\sin(\pi x);
\end{align}
and

$\bullet$ \emph{Initial condition 2:}
\begin{align}\label{sca_ini_2}
u_0(x)=\left(x+\dfrac{1}{2}\right)^k\exp\left(-100\left(x+\frac{1}{2}\right)^2\right),
\end{align}
where $k-1$ is the number of vanishing spatial derivatives at
$x=-\frac{1}{2}$, that is, $0=\frac{\partial f}{\partial
x}\big|_{x=0}=\ldots=\frac{\partial^{(k-1)}f}{\partial
x^{(k-1)}}\big|_{x=0}\ne \frac{\partial^{(k)}f}{\partial
x^{(k)}}\big|_{x=0}$.

{\footnotesize
\begin{table}[t]
\centering \caption{\small{Convergence of $u_t+u_x=0$ with initial
conditions (\ref{sca_ini_1}) and (\ref{sca_ini_2}), at time $t=1$.}}
\begin{tabular}{@{} p{1.3cm} p{0.5cm} | p{2.0cm} p{2.0cm} | p{2.0cm} p{2.0cm} | p{2.0cm} p{2.0cm} @{}}
\hline\hline
                &           &     \emph{Eq. (\ref{sca_ini_1})}    &       &    \emph{Eq. (\ref{sca_ini_2}), $k=2$}  &   &   \emph{Eq. (\ref{sca_ini_2}),  $k=3$}  \\
\cline{3-8}
     &      $N$      &     $L^1$ error      &    $L^{\infty}$ error    &   $L^1$ error      &    $L^{\infty}$ error    &   $L^1$ error      &    $L^{\infty}$ error\\
\cline{2-8}
WENO-           &   40   &   4.5E-07\ (-)      &   3.4E-07\ (-)     &   4.5E-04\ (-)      &   2.2E-03\ (-)        &   4.2E-05\ (-)      &     1.3E-04\ (-)       \\
CU6             &   80   &   6.9E-09\ (6.0)    &   5.4E-09\ (6.0)   &   3.8E-05\ (3.6)    &   2.0E-04\ (3.5)      &   8.9E-06\ (2.2)    &     4.5E-05\ (1.5)    \\
                &   160  &   1.1E-10\ (6.0)    &   8.4E-11\ (6.0)   &   6.5E-07\ (5.9)    &   3.6E-06\ (5.8)      &   1.3E-07\ (6.1)    &     8.1E-07\ (5.8)    \\
                &   320  &   4.1E-13\ (8.0)    &   3.8E-13\ (7.8)   &   1.1E-08\ (5.9)    &   6.0E-08\ (5.9)      &   1.8E-09\ (6.1)    &     1.1E-08\ (6.2)    \\
\cline{2-8}
WENO-           &   40   &   4.5E-07\ (-)      &   3.4E-07\ (-)     &   5.0E-04\ (-)      &   2.2E-03\ (-)        &   4.3E-05\ (-)      &     1.4E-04\ (-)       \\
NW6             &   80   &   6.9E-09\ (6.0)    &   5.3E-09\ (6.0)   &   4.1E-05\ (3.6)    &   2.2E-04\ (3.4)      &   8.7E-06\ (2.3)    &     4.7E-05\ (1.6)    \\
                &   160  &   1.1E-10\ (6.0)    &   8.4E-11\ (6.0)   &   6.4E-07\ (6.0)    &   3.6E-06\ (5.9)      &   1.4E-07\ (5.9)    &     9.6E-07\ (5.6)    \\
                &   320  &   4.1E-13\ (7.8)    &   3.6E-13\ (7.9)   &   1.0E-08\ (5.9)    &   6.0E-08\ (5.9)      &   1.8E-09\ (6.3)    &     1.1E-08\ (6.5)    \\
\cline{2-8}
WENO-           &   40   &   4.5E-07\ (-)      &   3.4E-07\ (-)     &   3.7E-04\ (-)      &   1.8E-03\ (-)        &   3.2E-05\ (-)      &     1.2E-04\ (-)        \\
$\theta$6       &   80   &   6.9E-09\ (6.0)    &   5.3E-09\ (6.0)   &   4.2E-05\ (3.1)    &   2.3E-04\ (2.9)      &   4.8E-06\ (2.7)    &     2.2E-05\ (2.4)     \\
                &   160  &   1.1E-10\ (6.0)    &   8.4E-11\ (6.0)   &   7.6E-07\ (5.8)    &   3.9E-06\ (5.9)      &   1.2E-07\ (5.3)    &     6.3E-07\ (5.1)     \\
                &   320  &   4.1E-13\ (8.0)    &   3.7E-13\ (7.8)   &   1.3E-08\ (5.9)    &   8.0E-08\ (5.6)      &   2.1E-09\ (5.9)    &     1.1E-08\ (5.9)     \\
\hline\hline
\end{tabular}
\label{tab_err}
\end{table}
}

$L^1$ and $L^{\infty}$ errors of 6th-order schemes at time $t=1$ are
measured and listed in Table \ref{tab_err} together with the order
of accuracy (in brackets), and are plotted in Fig. \ref{fig_err}. We
choose the time step $\Delta t=\Delta x^{6/3}$ so that the numerical
errors in time do not contribute to the results. In the figures, we
also show the errors of the 5th-order schemes for comparison. It is
observed that for all initial conditions, all schemes converge to
the designed order of accuracy. Moreover, the errors of the new
WENO-$\theta6$ schemes are almost indistinguishable with those of
other 6th-order ones, or even better at some grid sizes.

\begin{figure}[h]
  \begin{center}
    \begin{tabular}{ccc}
      \resizebox{52mm}{!}{\includegraphics{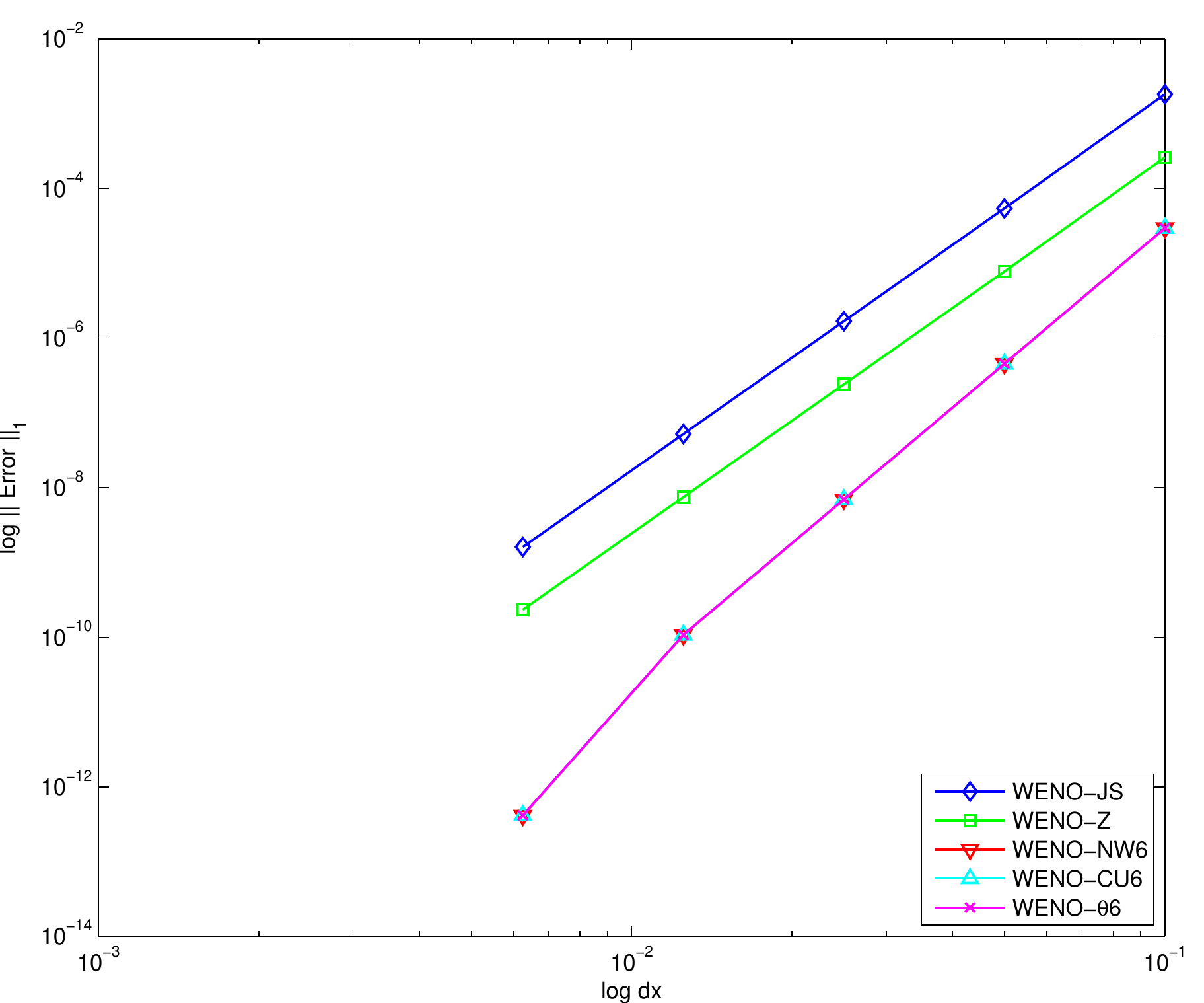}} &
      \resizebox{52mm}{!}{\includegraphics{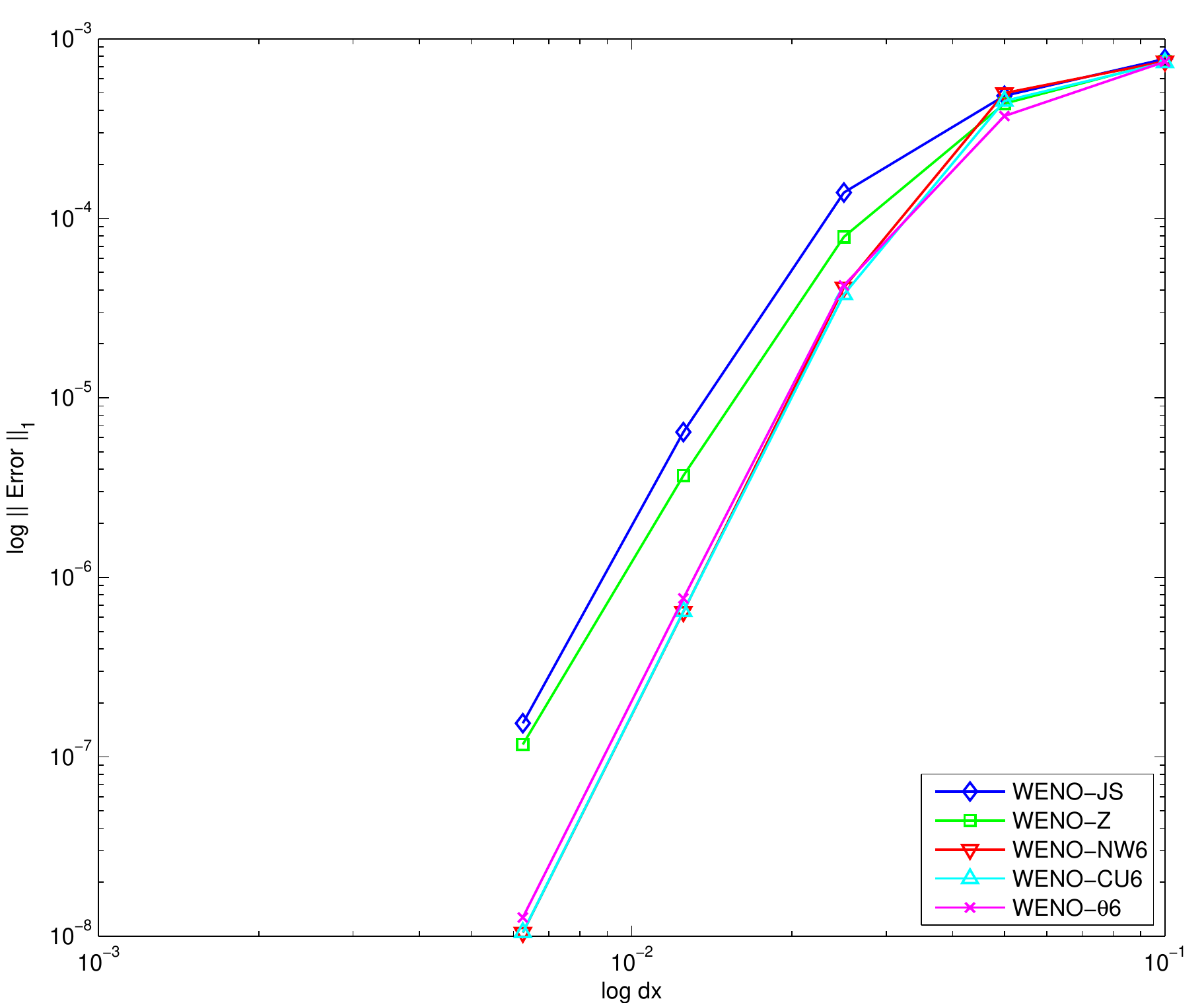}} &
      \resizebox{52mm}{!}{\includegraphics{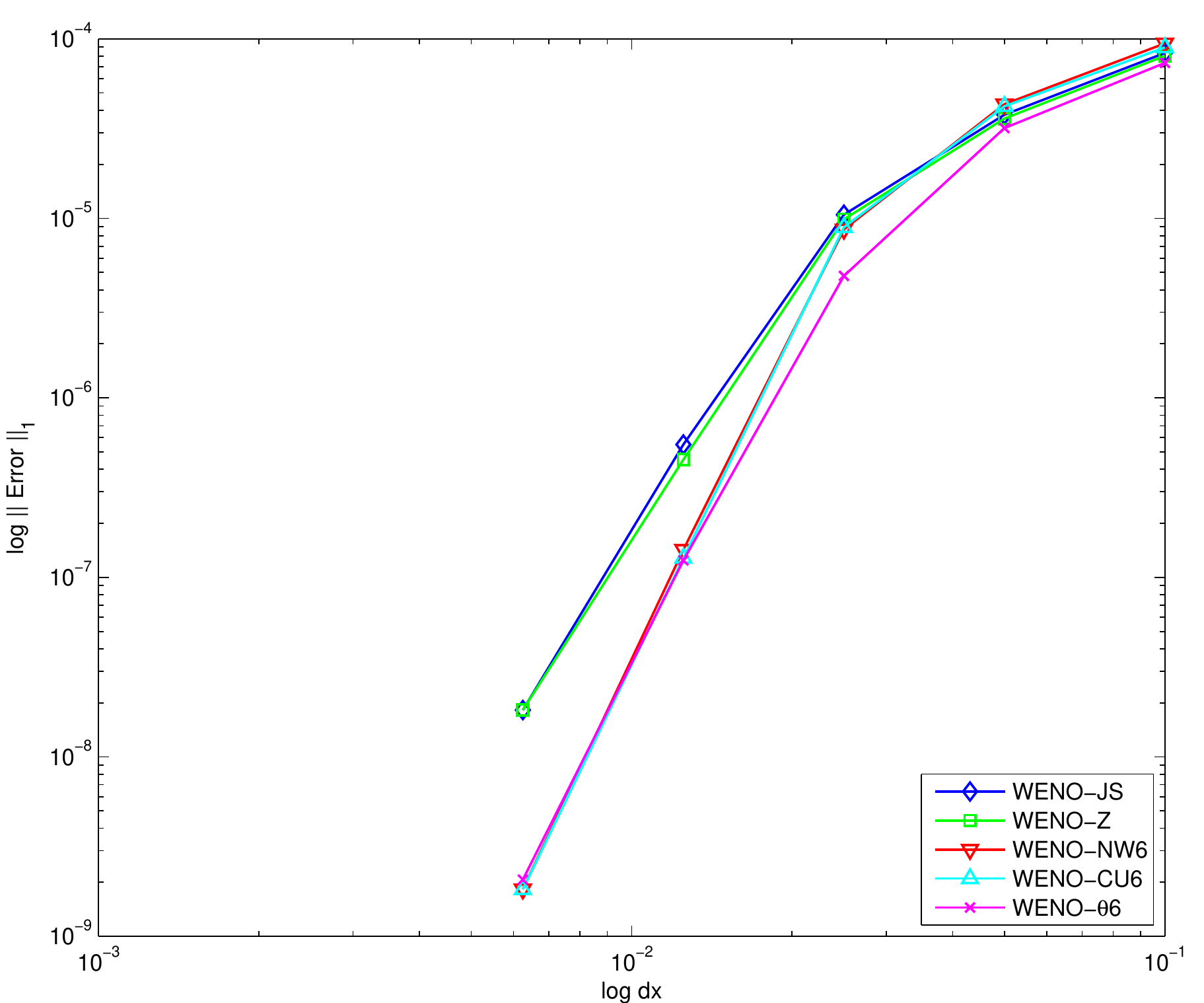}}\\
      \resizebox{52mm}{!}{\includegraphics{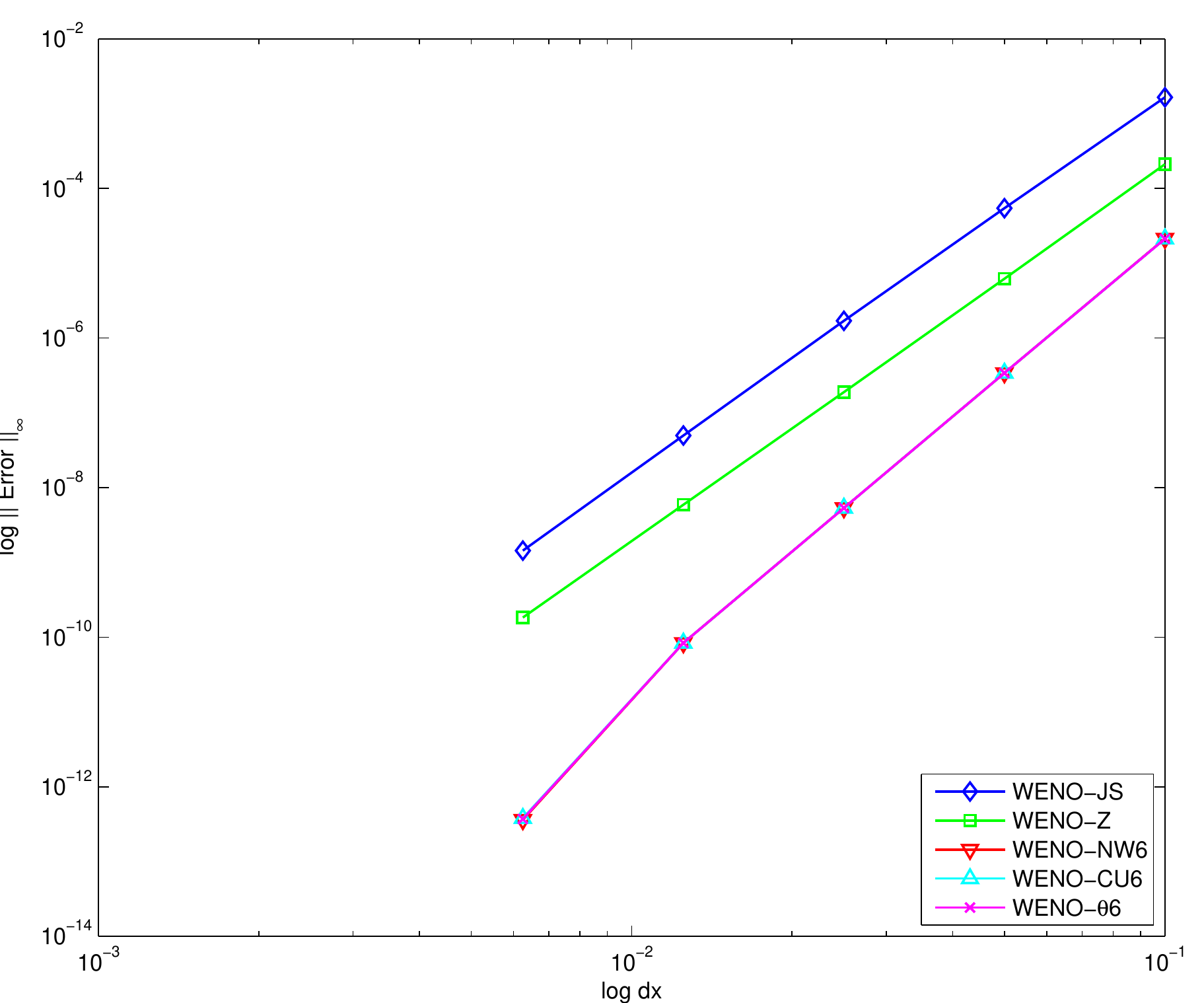}}&
      \resizebox{52mm}{!}{\includegraphics{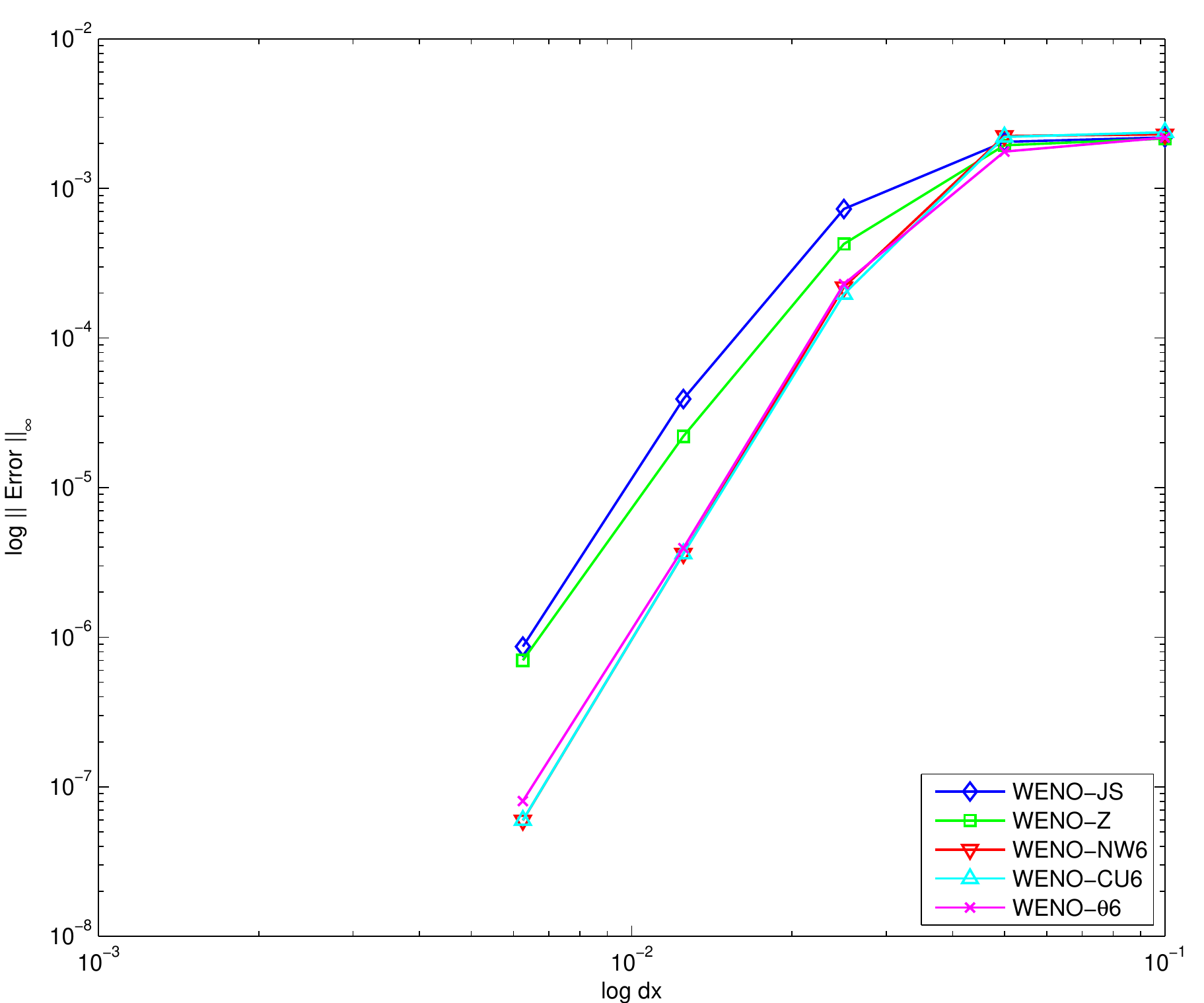}}&
      !\resizebox{52mm}{!}{\includegraphics{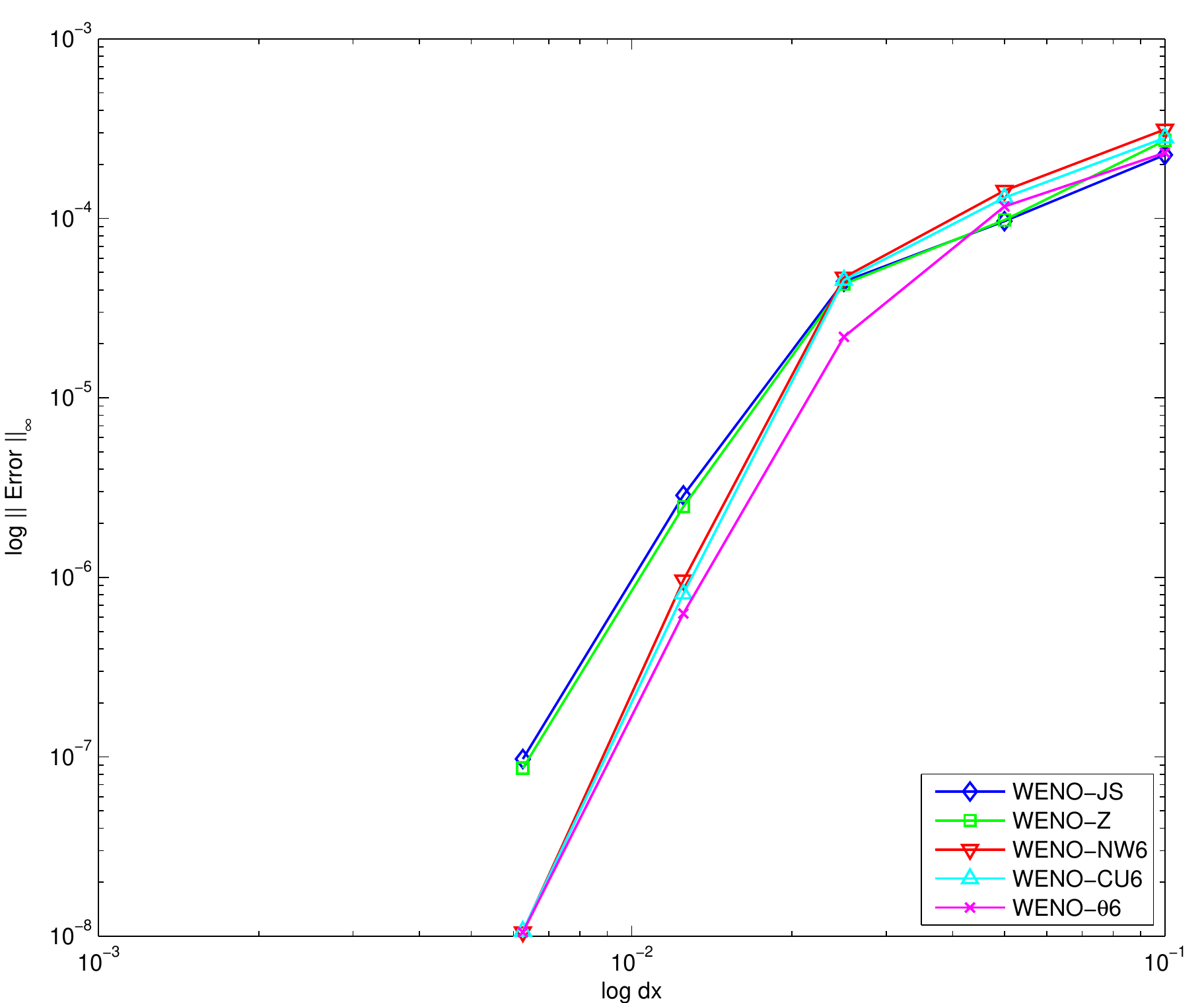}}
    \end{tabular}
    \caption{\small{Convergence of Eq. (\ref{lin_eq}) at time $t=1$.
    Left: Initial condition (\ref{sca_ini_1}); Middle: Initial condition (\ref{sca_ini_2}) with $k=2$.
    Right: Initial condition (\ref{sca_ini_2}) with $k=3$. Top: in $L^1$ norm; Bottom: in $L^{\infty}$ norm.}}
  \label{fig_err}
  \end{center}
\end{figure}

\subsection{Resolution Tests}

We now test if our new scheme overcomes the loss of accuracy of
WENO-CU6 and WENO-NW6. We revisit the initial condition given in
example \ref{exam_1} which is as follows,
\begin{align}\label{sca_ini_3}
u_0(x)=\max(\sin(\pi x),0),\quad x\in(-1,1).
\end{align}

The numerical solution obtained from our new scheme is added and
shown in Fig. \ref{fig_soln_ini3}, together with those given in
example \ref{exam_1}. We also plot the pointwise errors in
$L^{\infty}$ norm in the same figure. It is shown that
WENO-$\theta6$ approximates the critical region around $x=-0.1$ much
better than WENO-CU6 and WENO-NW6. Indeed, the pointwise errors of
the former around this region is comparable to those of WENO-JS and
WENO-Z, where WENO-NW6 and WENO-CU6 show a loss of accuracy, which
in turn causes problems for approximating solutions where symmetry
is required. See below tests for numerical evidence.

\begin{figure}
  \begin{center}
    \begin{tabular}{ccc}
      \resizebox{52mm}{!}{\includegraphics{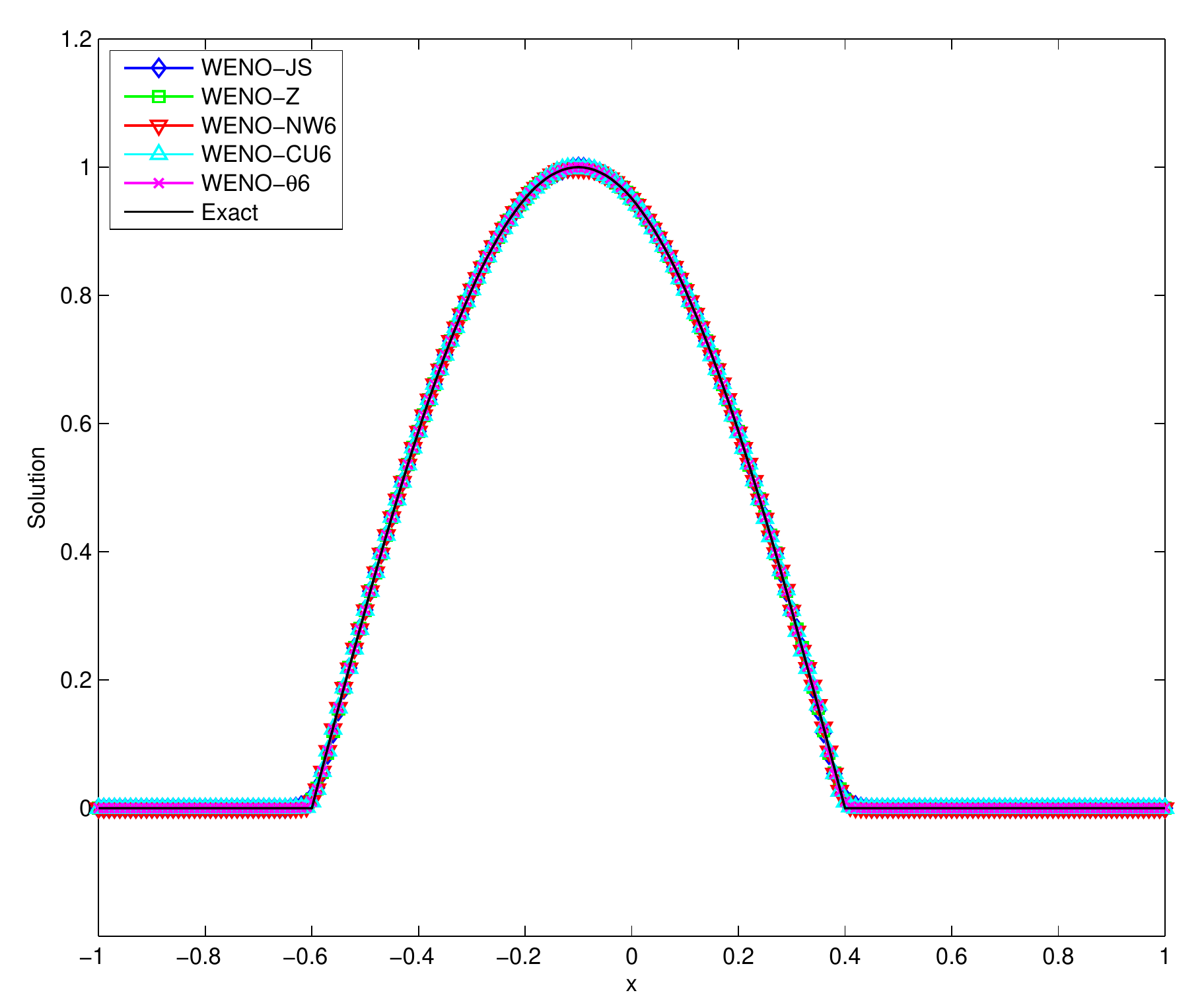}} &
      \resizebox{52mm}{!}{\includegraphics{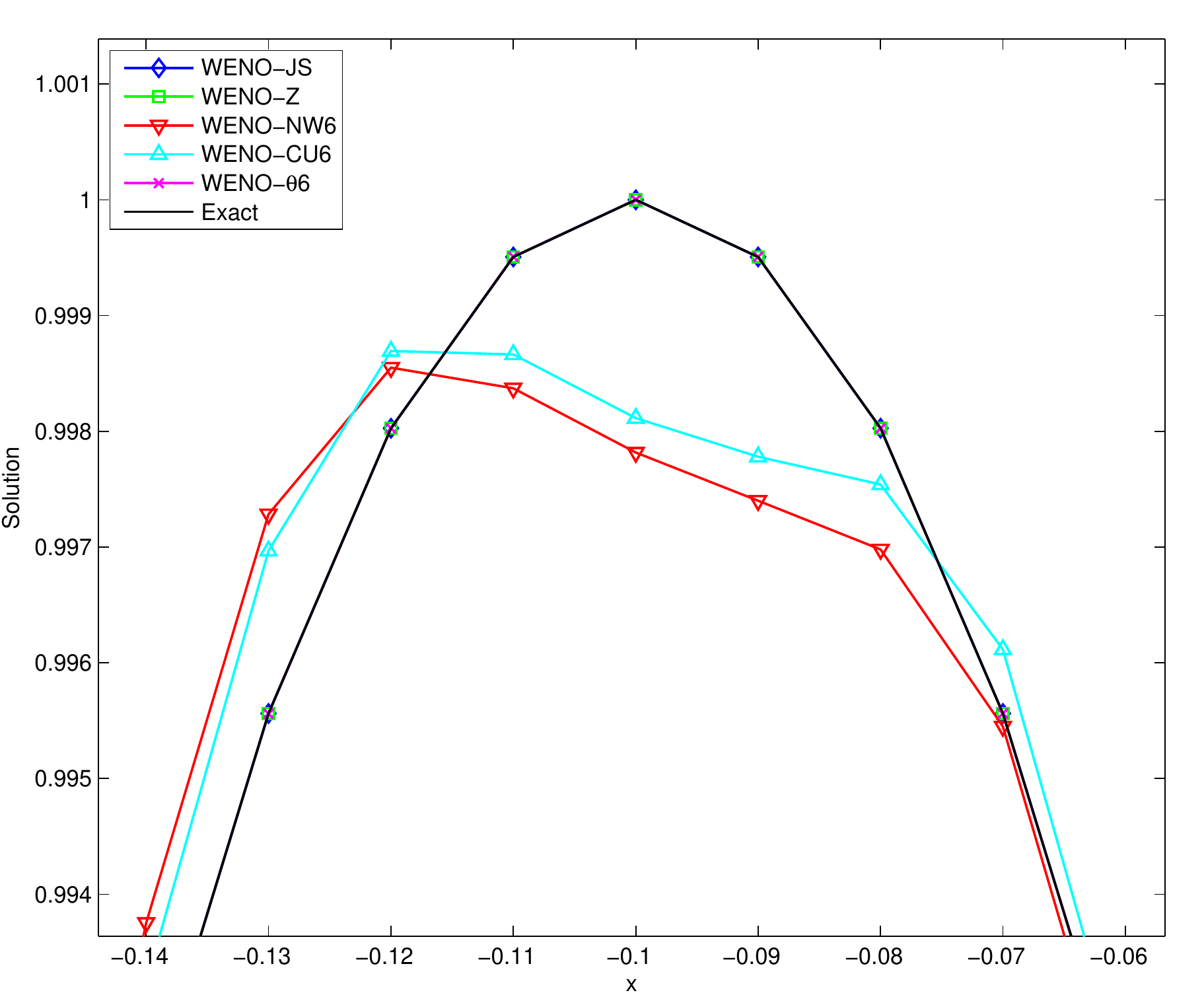}}&
      \resizebox{52mm}{!}{\includegraphics{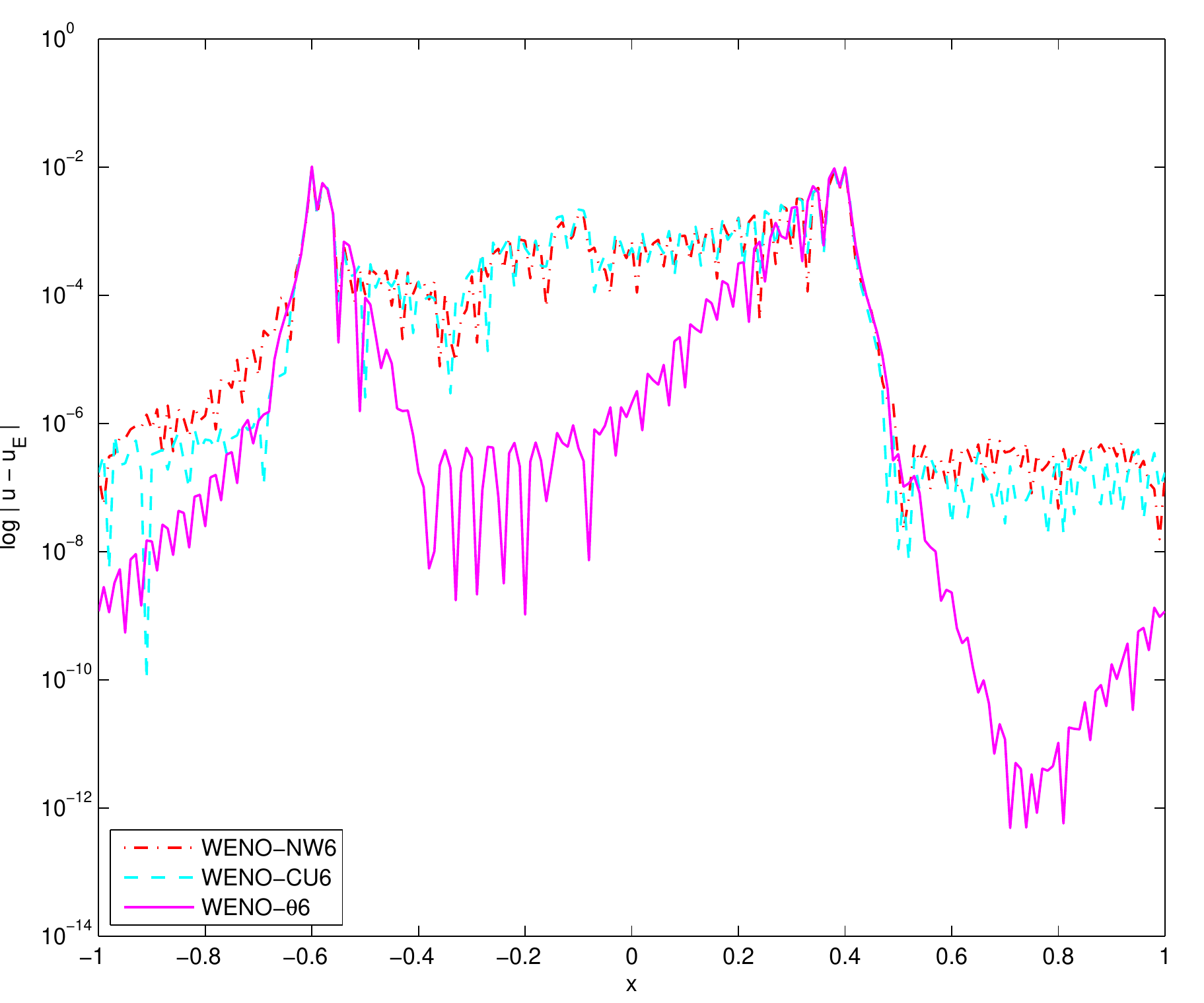}}
    \end{tabular}
    \caption{\small{Left: Numerical solutions of Eq. (\ref{mod_con_scalar}) with initial condition
    (\ref{sca_ini_3}) at time $t=2.4$. Grid $200$.
    Middle: Zoom near the critical point.
    Right: Pointwise errors in log scale.}}
    \label{fig_soln_ini3}
  \end{center}
\end{figure}



The nonlinear weights $\omega_k$'s of these schemes are plotted in
Fig. \ref{fig_om_ini3}. We observe that around the critical point,
the $\omega_k$'s of WENO-Z and WENO-$\theta6$ are stable and
converge to their optimal values $\gamma_k$'s. We note that for the
latter scheme, the nonlinear weights keep fluctuating between the
optimal weights of the 5th-order upwind and 6th-order central linear
schemes. This clearly shows the effect of switching mechanism
(\ref{tau_we6}) in improving the accuracy of the scheme near the
critical region. We also notify the non-convergence of $\omega_k$'s
of the WENO-CU6 and WENO-NW6 schemes around the critical region.
This shows the improvement of our new scheme over the other 6th-order ones.

\begin{figure}
  \begin{center}
    \begin{tabular}{cc}
      \resizebox{70mm}{!}{\includegraphics{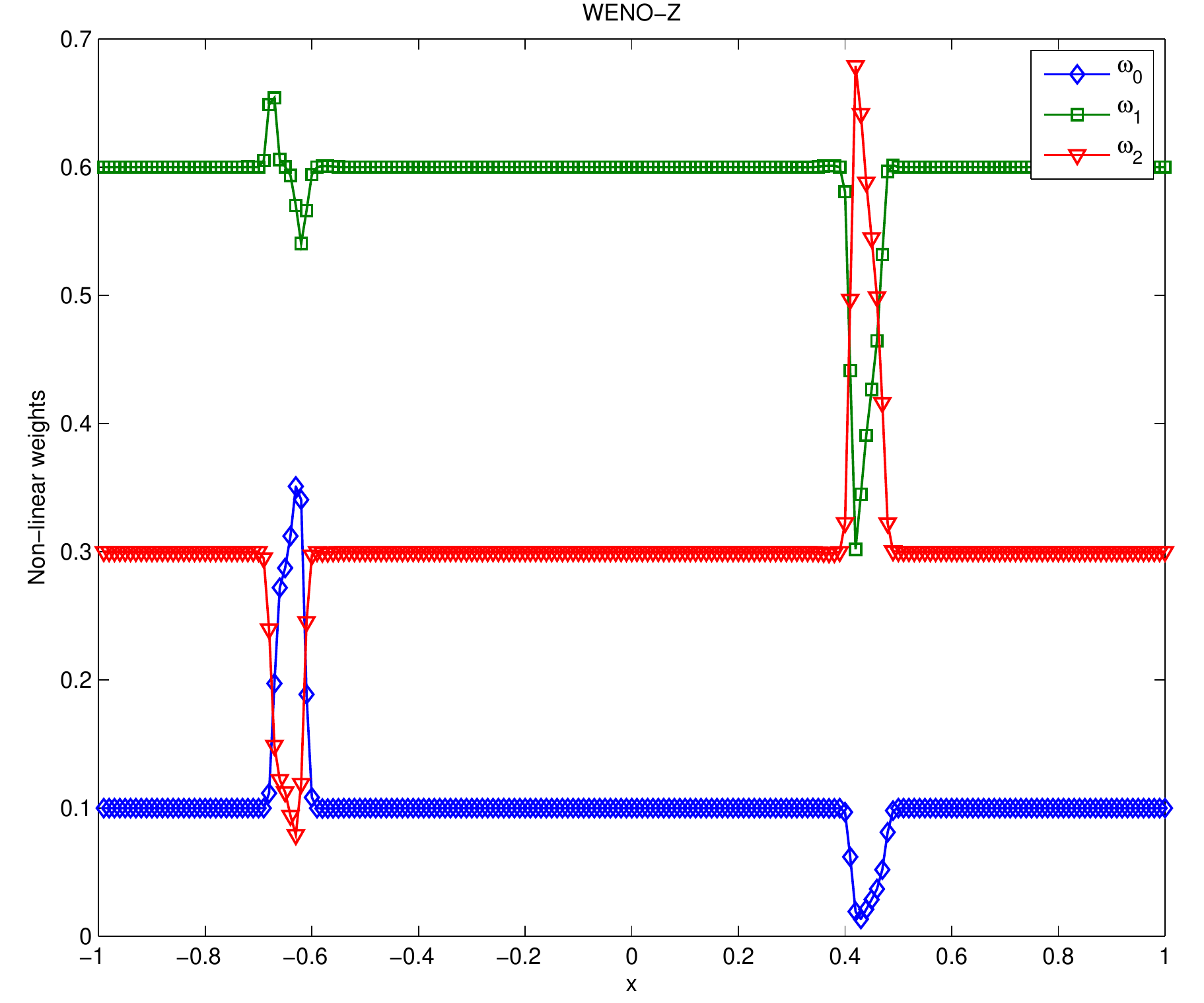}} &
      \resizebox{70mm}{!}{\includegraphics{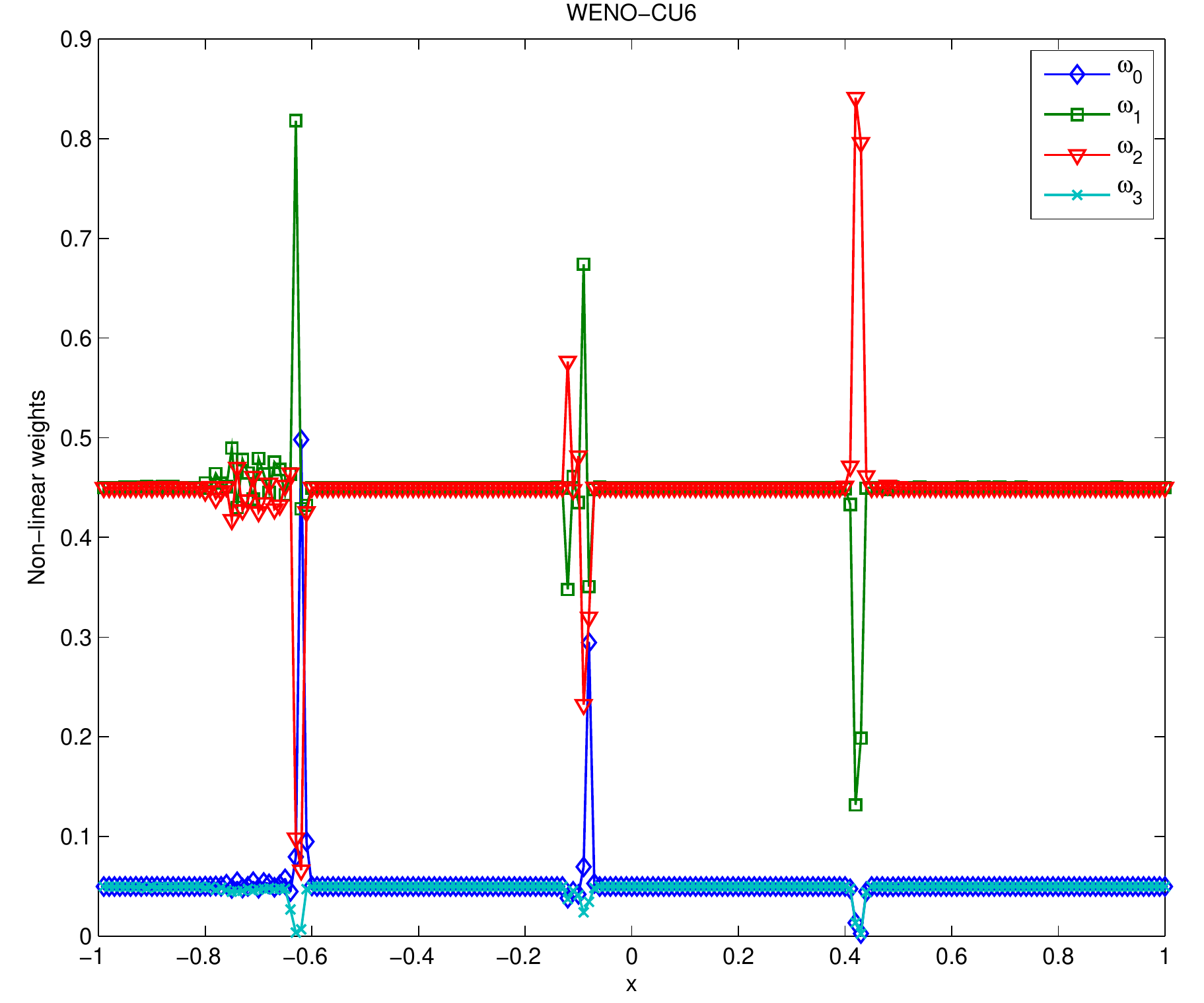}}\\
      \resizebox{70mm}{!}{\includegraphics{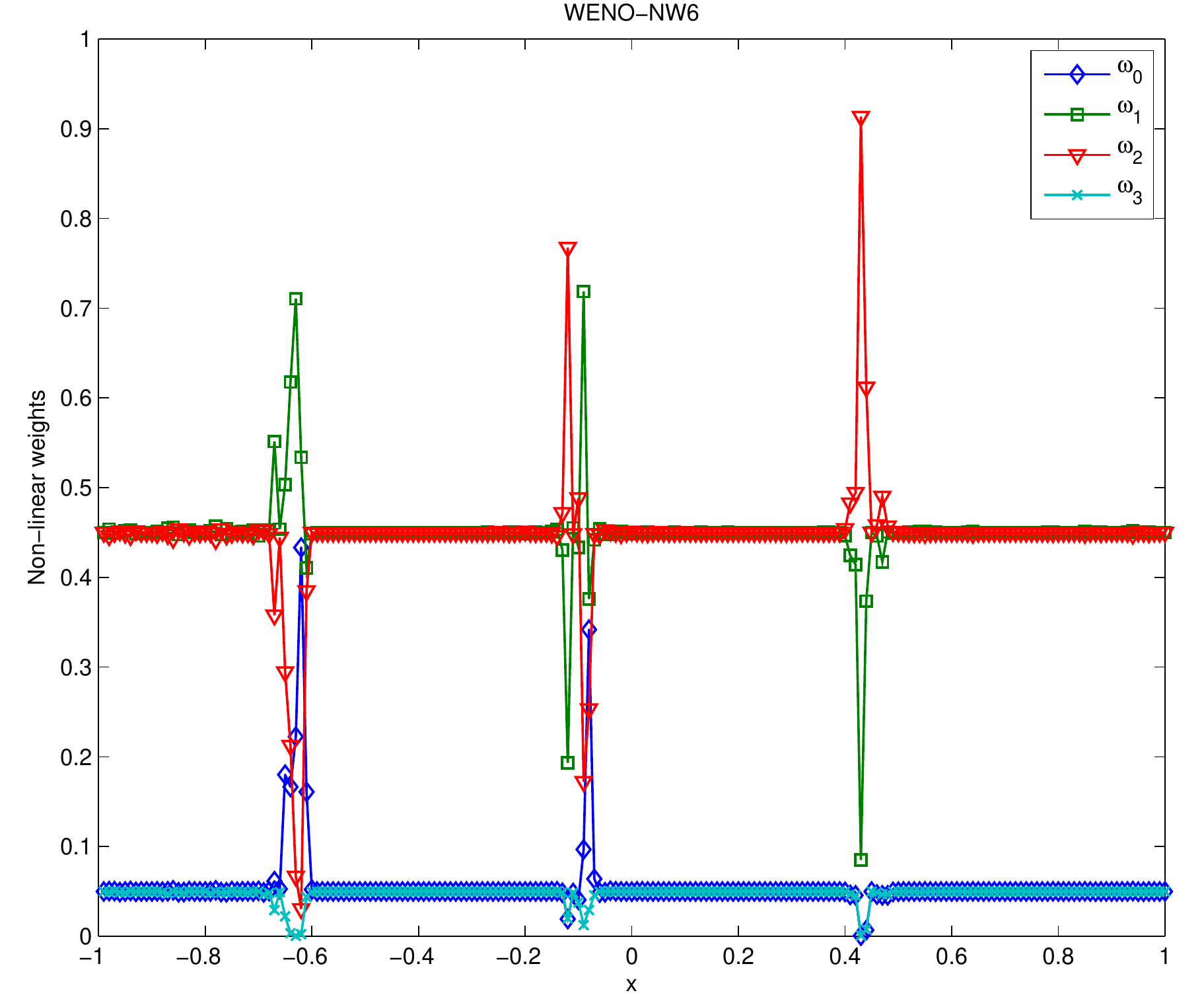}} &
      \resizebox{70mm}{!}{\includegraphics{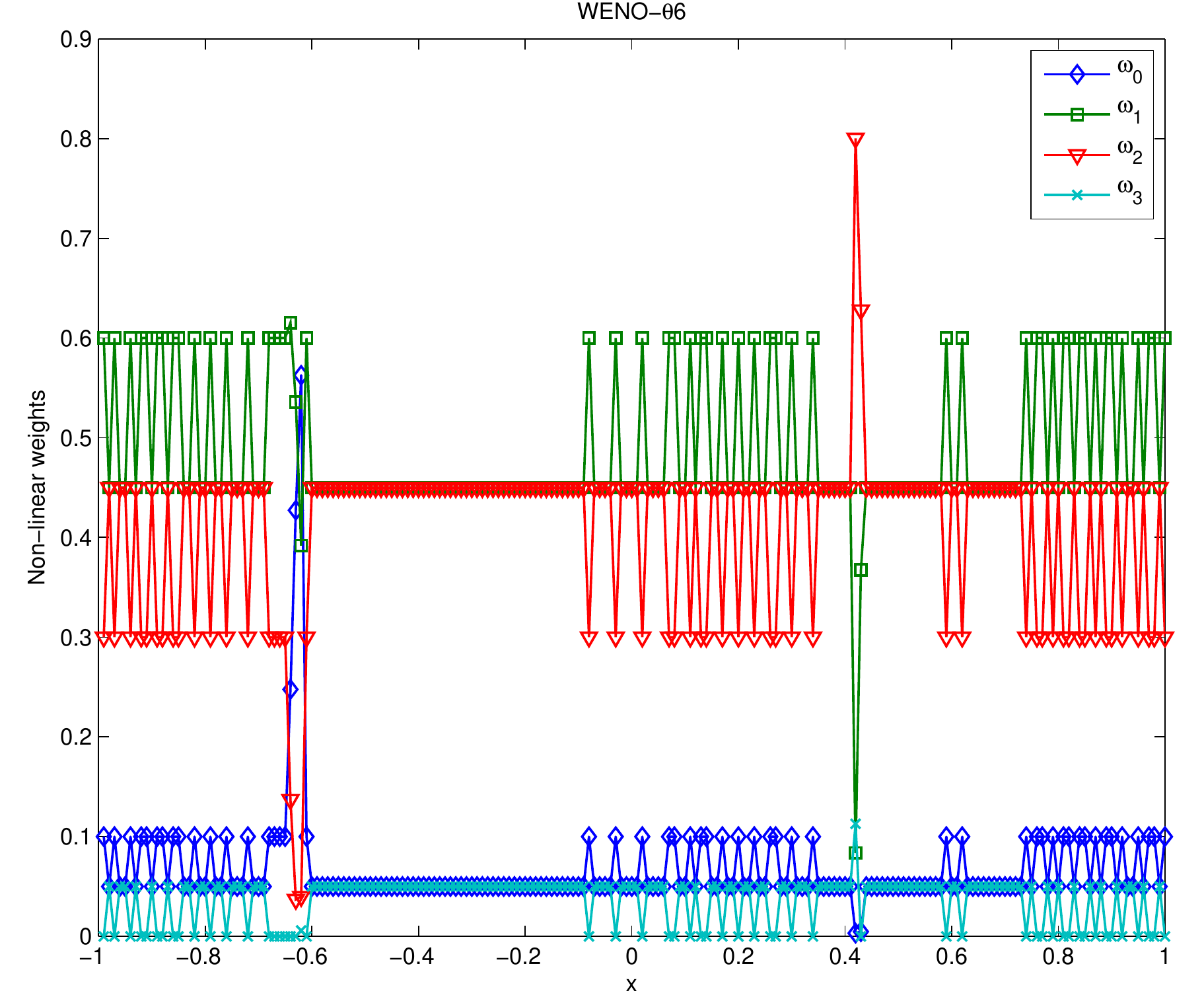}}
    \end{tabular}
    \caption{\small{Distribution of the non-linear weights for the initial data (\ref{sca_ini_3}).
    From top to bottom, left to right: WENO-Z, WENO-CU6, WENO-NW6, and
    WENO-$\theta$6.}}
    \label{fig_om_ini3}
  \end{center}
\end{figure}

\section{Numerical Results}

In this section, we perform a number of tests to compare the results
of our new scheme with those obtained from the other WENO schemes,
including the 5th-order upwind WENO-JS, WENO-Z, and the 6th-order
central WENO-CU6, WENO-NW6.

\subsection{Scalar Conservation Laws}

\subsubsection{TEST 1: Linear Case}

We solve the one-dimensional linear advection equation
(\ref{lin_eq}) with the following initial condition $u_0(x)$ which
contains a $C^{\infty}$ Gaussian, a square wave, a triangle, and a
semi-ellipse (see \cite{JS}),
\begin{align}\label{sca_ini_4}
u_0(x)=
\begin{cases}
&\frac{1}{6}[G(x,\beta,z-\delta)+4G(x,\beta,z)+G(x,\beta,z+\delta)],\quad
-0.8\le x\le-0.6,\\
&1,\quad\quad\quad\quad\quad\quad\quad\quad\quad\quad\quad\quad\quad\quad\quad\quad\quad\quad\quad\quad\
\
-0.4\le x\le-0.2,\\
&1-|10(x-0.1)|,\quad\quad\quad\quad\quad\quad\quad\quad\quad\quad\quad\quad\quad\quad
0\le x\le0.2,\\
&\frac{1}{6}[F(x,\alpha,a-\delta)+4F(x,\alpha,a)+F(x,\alpha,a+\delta)],\quad
0.4\le x\le0.6,\\
&0,\quad\quad\quad\quad\quad\quad\quad\quad\quad\quad\quad\quad\quad\quad\quad\quad\quad\quad\quad\quad\
\ \text{otherwise},
\end{cases}
\end{align}
where
\begin{align}
G(x,\beta,z)=\exp(-\beta(x-z)^2),
\end{align}
\begin{align}
F(x,\alpha,a)=\sqrt{\max(1-\alpha^2(x-a)^2,0)};
\end{align}
the constants are $z=-0.7$, $\delta=0.005$,
$\beta=\frac{\log2}{36\delta^2}$, $a=0.5$, and $\alpha=10$.

We compute the solution up to time $t=6.3$ with $N=400$ and periodic
boundary conditions. The results obtained from the WENO-CU6,
WENO-NW6, and WENO-$\theta6$ schemes are plotted in Fig.
\ref{fig_soln_lin}. We choose $\alpha_R=50$. Zooms around the shocks
and top of the semi-ellipse are also shown in the same figure. It is
observed that WENO-$\theta6$ is comparable to WENO-NW6 in capturing
the shocks, but the former is much better than the latter and
WENO-CU6 in approximating top of the semi-ellipse.

\begin{figure}
  \begin{center}
    \begin{tabular}{ccc}
      \resizebox{52mm}{!}{\includegraphics{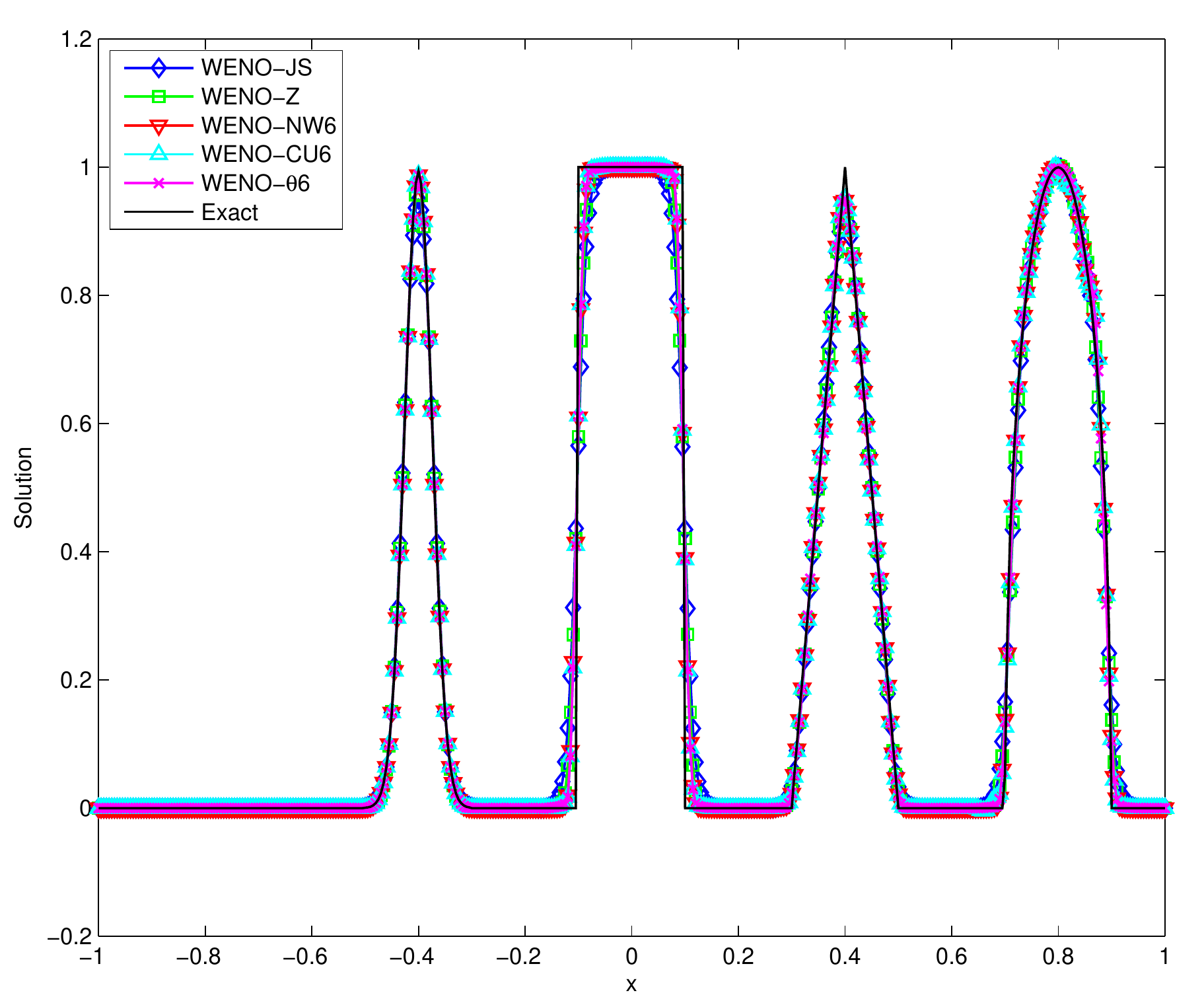}} &
      \resizebox{52mm}{!}{\includegraphics{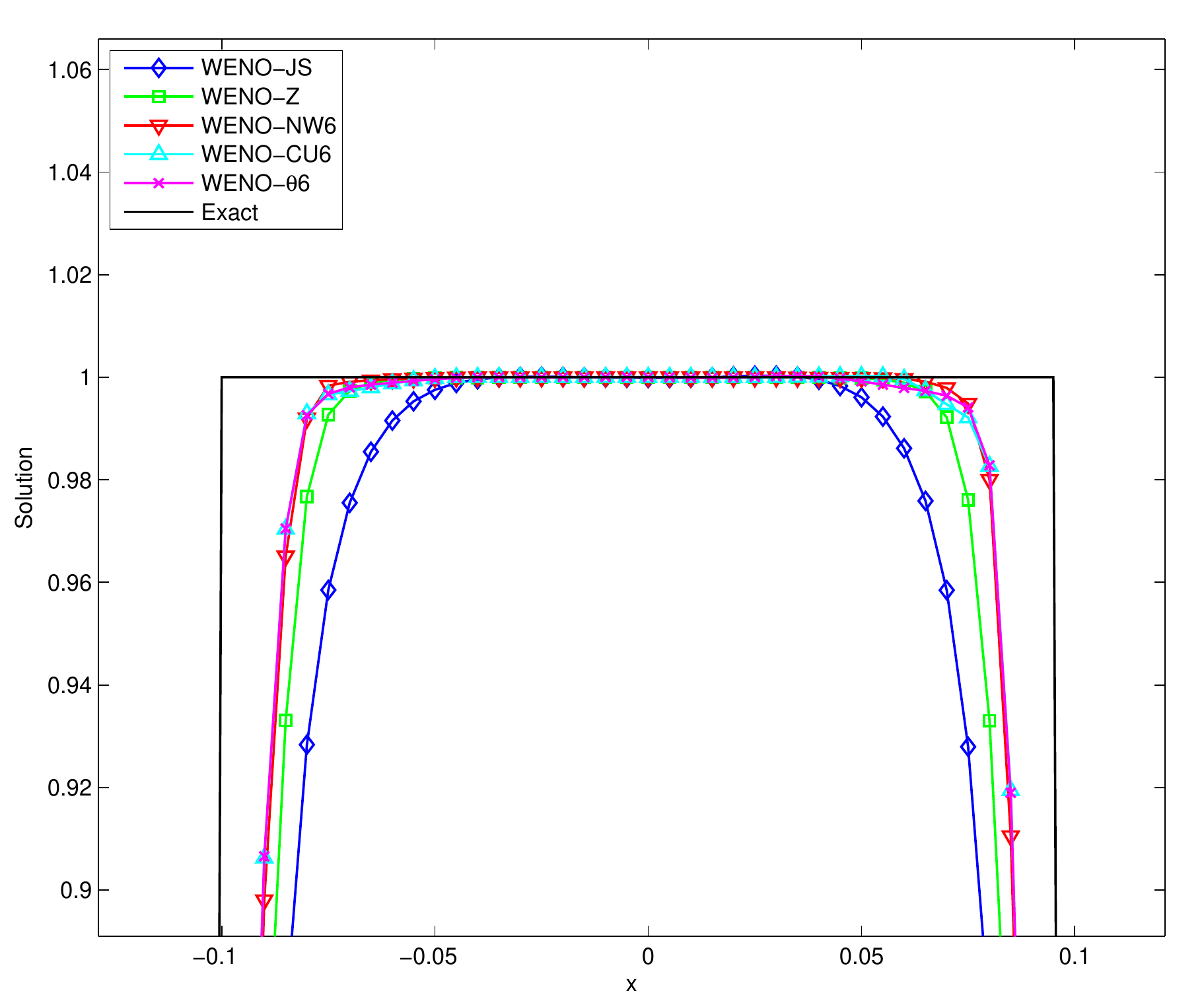}}&
      \resizebox{52mm}{!}{\includegraphics{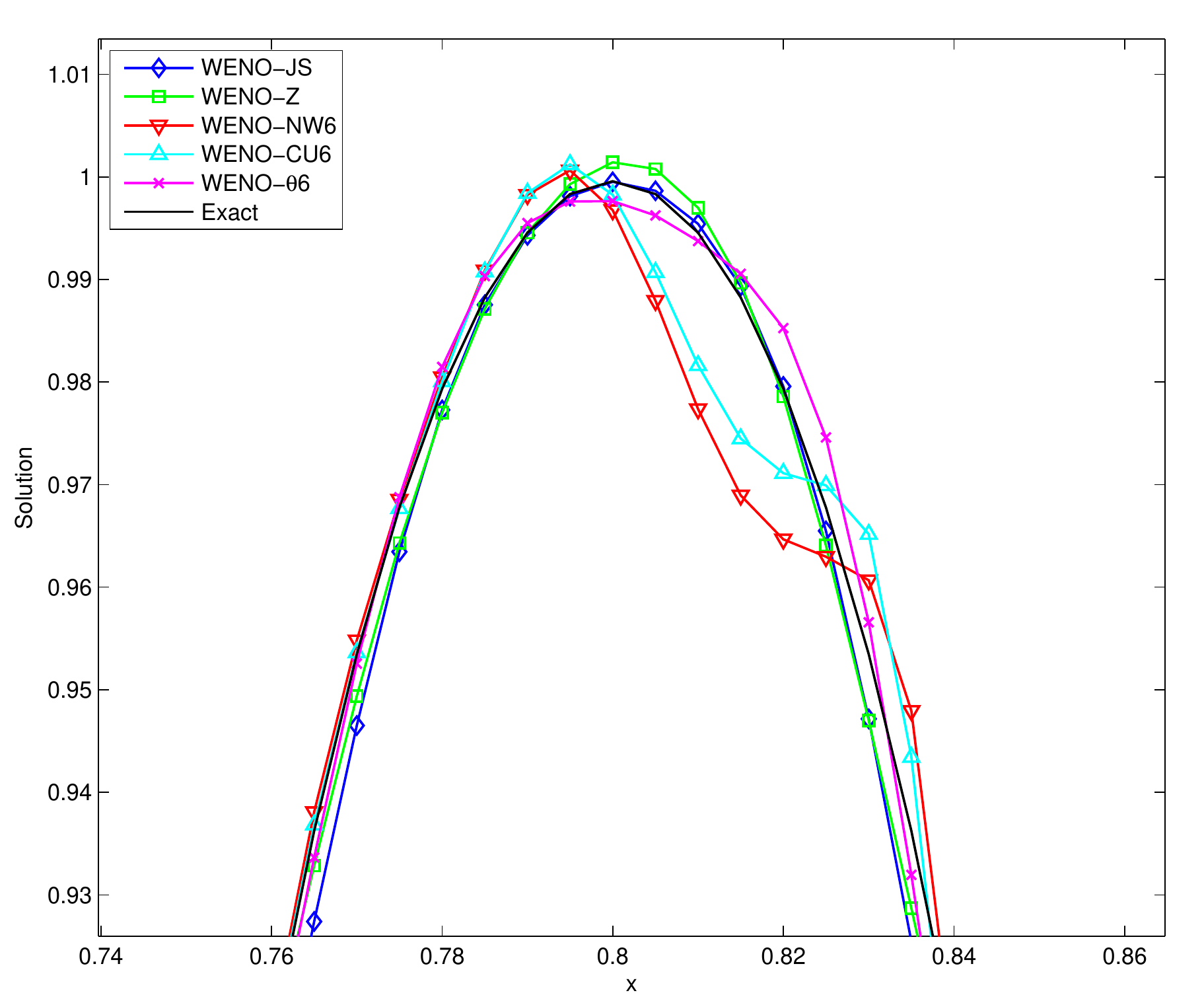}}
    \end{tabular}
    \caption{\small{Left: Linear advection Eq. (\ref{lin_eq}) with initial condition
    (\ref{sca_ini_4}). Time $t=6.3$. Grid $400$. The others: zooms at critical regions.}}
    \label{fig_soln_lin}
  \end{center}
\end{figure}

\subsubsection{TEST 2: Nonlinear Case}

For this, we choose the Burgers equation
\begin{align}\label{bur_eq}
\begin{cases}
u_t+\left(\dfrac{u^2}{2}\right)_x=0,\quad x\in(-1,1),\\
u(x,0)=u_0(x),
\end{cases}
\end{align}
subject to periodic boundary conditions.

In Fig. \ref{fig_soln_bur}, we show the numerical results of the
6th-order WENO schemes for the initial condition
\begin{align}\label{bur_ini_1}
u_0(x)=-\sin(\pi x);
\end{align}
at time $t=1.5$ and
\begin{align}\label{bur_ini_2}
u_0(x)=\dfrac{1}{2}+\sin(\pi x)
\end{align}
at $t=0.55$. We choose a grid of $N=200$ grid intervals and
$\alpha_R=50$. It is shown
that the shocks are very well captured by all schemes. 

\begin{figure}
  \begin{center}
    \begin{tabular}{cc}
      \resizebox{70mm}{!}{\includegraphics{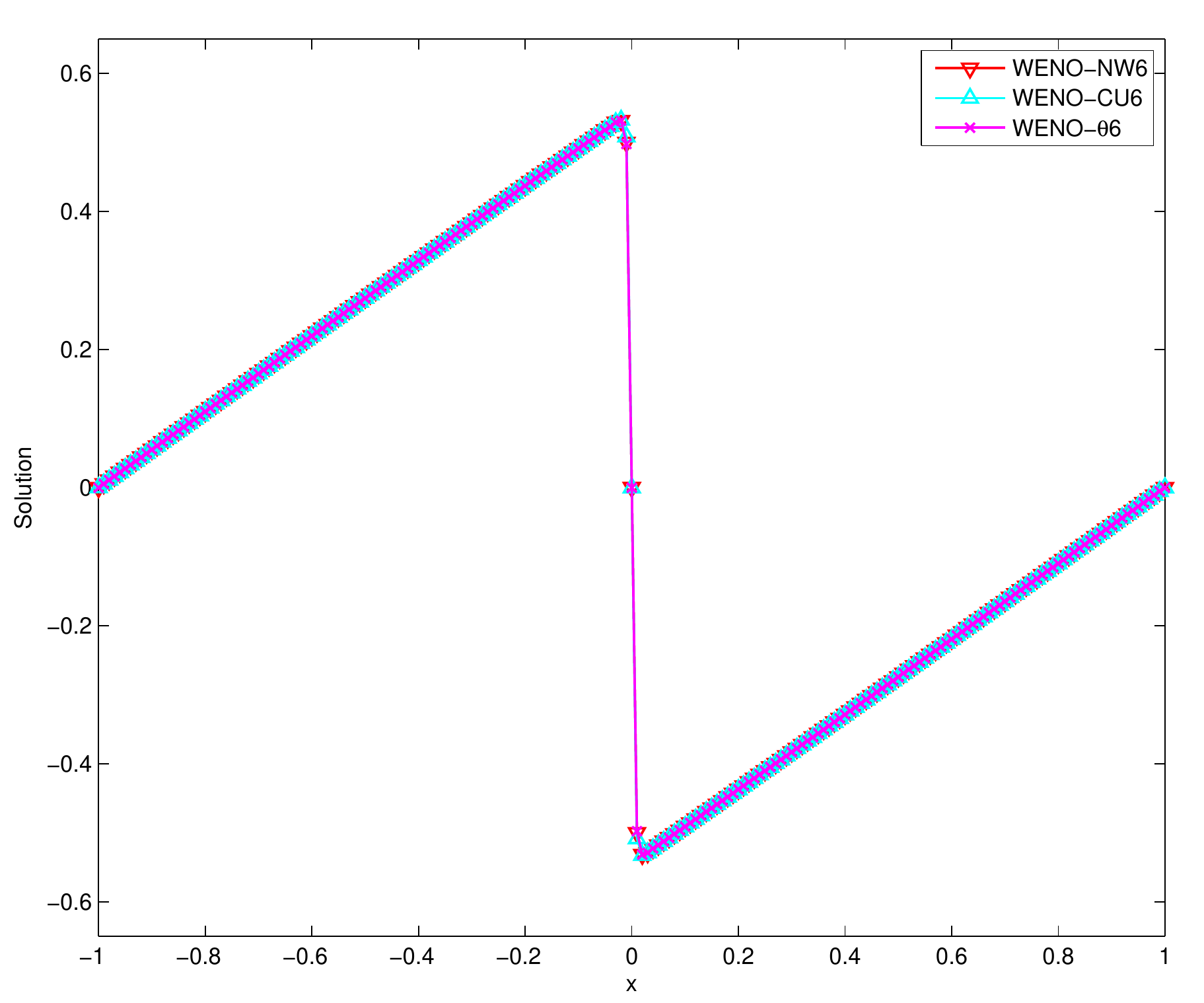}} &
      \resizebox{70mm}{!}{\includegraphics{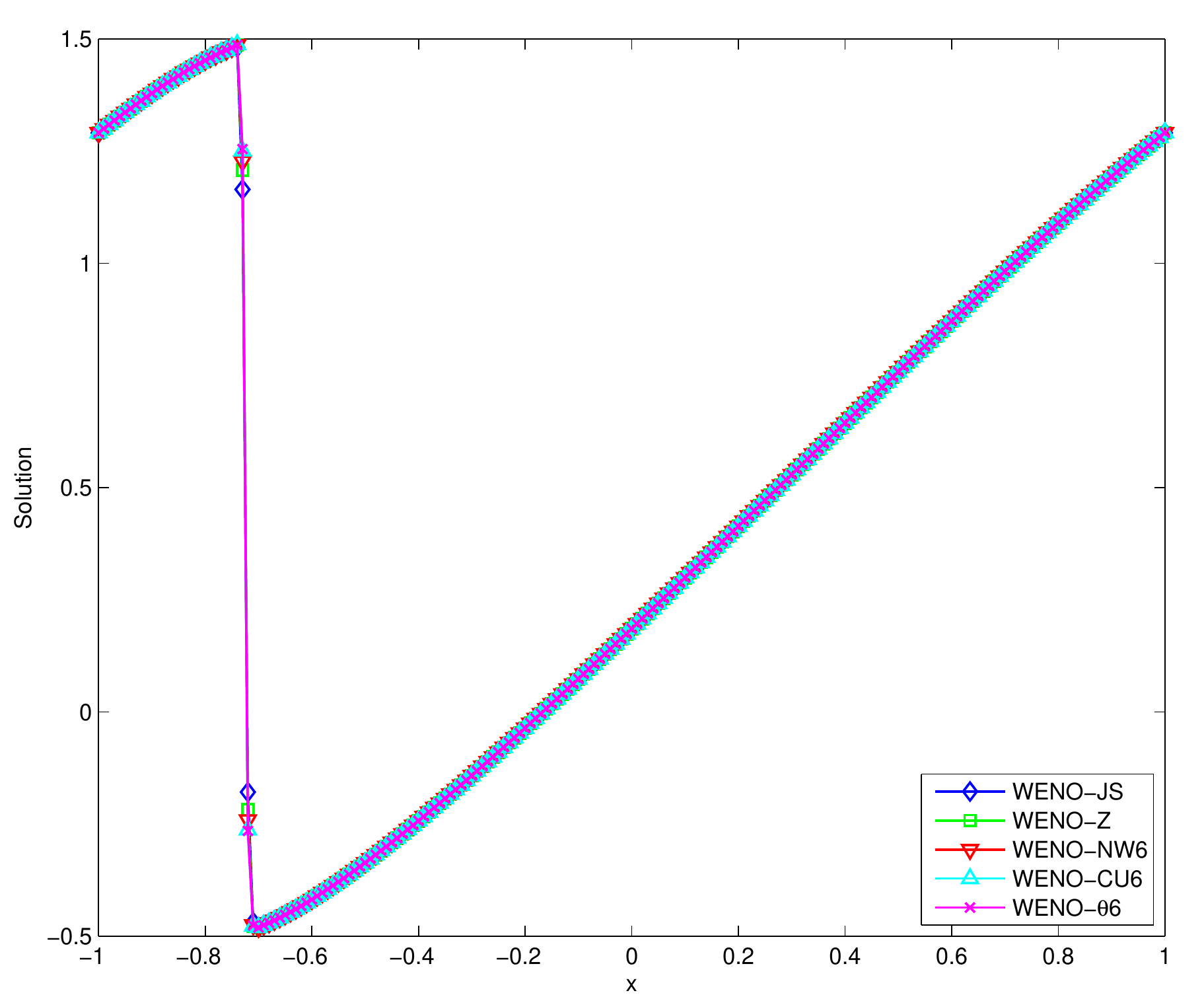}}
    \end{tabular}
    \caption{\small{Burgers' Eq. (\ref{bur_eq}). Grid $200$.
    Left: initial condition (\ref{bur_ini_1}) at time $t=1.5$;
    Right: initial condition (\ref{bur_ini_2}) at time $t=0.55$.}}
    \label{fig_soln_bur}
  \end{center}
\end{figure}

\subsection{Euler Equations of Gas Dynamics}

In this subsection, we consider the one-dimensional Euler equations
of gas dynamics given in the following,
\begin{align}\label{eu_eq}
\mathbf{u}_t+\mathbf{f}(\mathbf{u})_x=0,
\end{align}
where
\begin{align}
\mathbf{u}=(\rho,\rho u,E)^T,\quad\mathbf{f}(\mathbf{u})=(\rho
u,p+\rho u^2,(E+p)u)^T,
\end{align}
where $\rho$, $u$, $p$, $E$ are density, velocity, pressure, and
total energy, respectively. The equation of state is as follows,
\begin{align}
p=(\gamma-1)\left(E-\dfrac{1}{2}\rho u^2\right),\quad \gamma=1.4,
\end{align}
where $\gamma$ is the ratio of specific heats. More details on the
Euler equations can be found in, e.g., \cite{LV}, \cite{To}.

For all below numerical simulations, we apply the WENO schemes in
characteristic fields of the flux $\mathbf{f}(\mathbf{u})$. That is,
we first find an average of the Jacobian $A_{j+\frac{1}{2}}$ at the
interface $x=x_{j+\frac{1}{2}}$. For this, we apply the Roe's mean
matrix (see \cite{Ro}).  Then the eigenvalues $\lambda_s$'s,
$L=[\mathbf{l}_s]_{s=1}^m$, $R=[\mathbf{r}_s]_{s=1}^m$ the complete
sets of the left and right eigenvectors, respectively, of
$A_{j+\frac{1}{2}}$ are determined. We next project the flux
$\mathbf{f}(\mathbf{u})$ into the characteristic fields by left
multiplying it with $L$. WENO schemes with a global Lax-Friedrichs
flux splitting are applied to approximate the components of the
flux. After that, the approximation in each characteristic field is
projected back to the component space by a right multiplying with
the matrix $R$. 

\subsubsection{TEST 3: Riemann Problems}

We consider the shock-tube problems which are Eq. (\ref{eu_eq}) with
Riemann initial data. In particular, the Sod problem, the Lax
problem, and the $123$ problem are given below.

$\bullet$ Sod's problem:
\begin{align}\label{sod_prob}
(\rho,u,p)=
\begin{cases}
&(0.125,\ \ 0,\ \ 0.1),\quad -5<x<0,\\
&(1,\ \ 0,\ \ 1),\quad 0<x<5;
\end{cases}
\end{align}
and the final time $t=1.7$.

$\bullet$ Lax's problem:
\begin{align}\label{lax_prob}
(\rho,u,p)=
\begin{cases}
&(0.445,\ \ 0.698,\ \ 3.528),\quad -5<x<0,\\
&(0.5,\ \ 0,\ \ 0.571),\quad 0<x<5;
\end{cases}
\end{align}
and the final time $t=1.3$.

$\bullet$ 123 problem:
\begin{align}\label{123_prob}
(\rho,u,p)=
\begin{cases}
&(1,\ \ -2,\ \ 0.4),\quad -5<x<0,\\
&(1,\ \ 2,\ \ 0.4),\quad 0<x<5;
\end{cases}
\end{align}
and the final time $t=1$.

We apply a transmissive condition at both boundaries. The exact
solution of these shock-tube problems can be found in, for example,
\cite{To}. For Lax's and the blast waves problems, we choose
$\alpha_R=10$; and $\alpha_R=50$ for other problems.

Numerical results of the density obtained from all WENO schemes with
a grid of $N=300$ are shown in Figs. \ref{fig_soln_sod} -
\ref{fig_soln_123}, respectively. We observe that for Sod's and
Lax's problems, there are overshoots at the contact discontinuities
for WENO-CU6, whereas WENO-$\theta6$ gives the sharpest capturing
without generating oscillations. For the $123$ problem, WENO-CU6
shows the most wiggling behavior around the trivial contact
discontinuity.


\begin{figure}
  \begin{center}
    \begin{tabular}{cc}
      \resizebox{70mm}{!}{\includegraphics{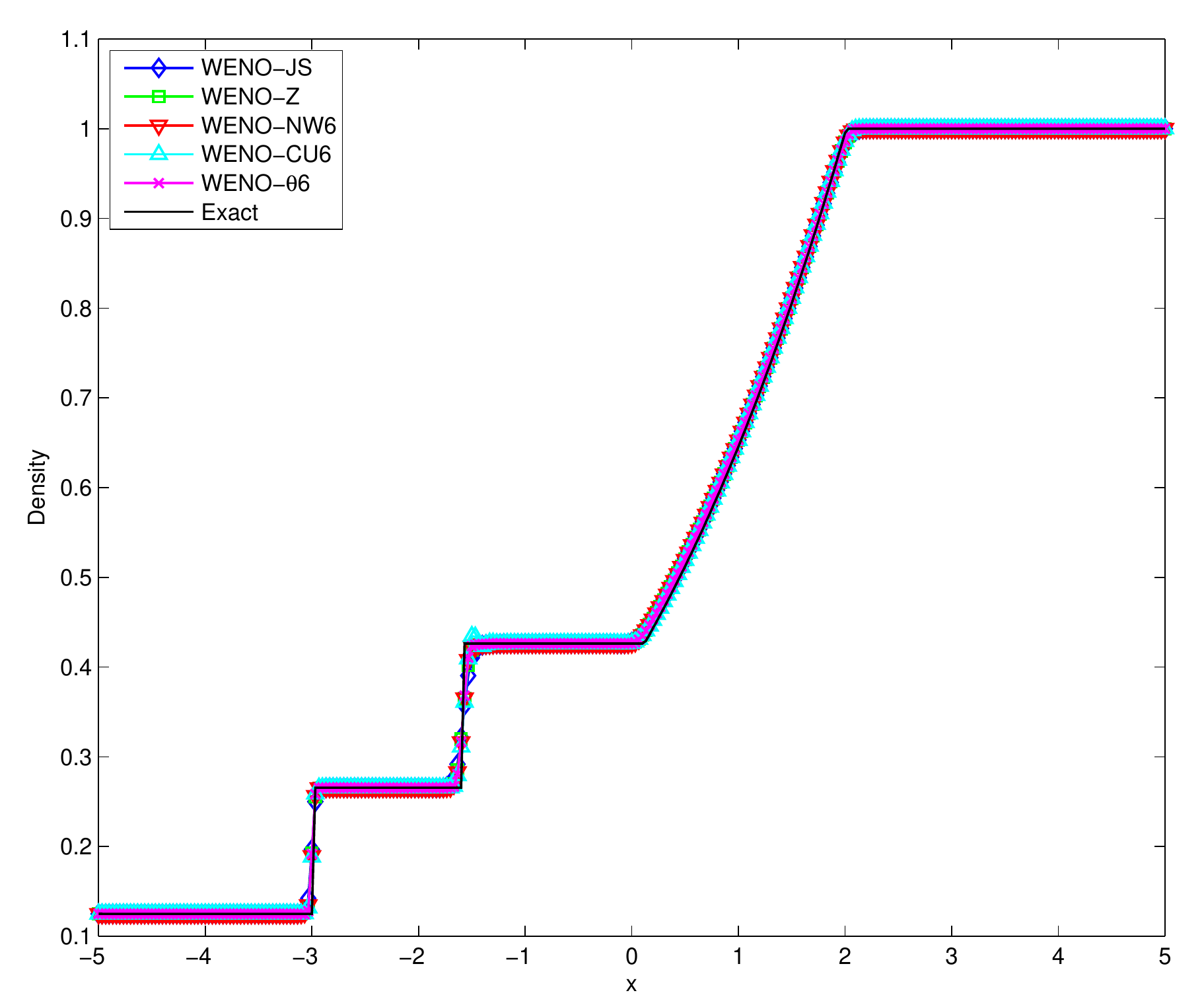}} &
      \resizebox{70mm}{!}{\includegraphics{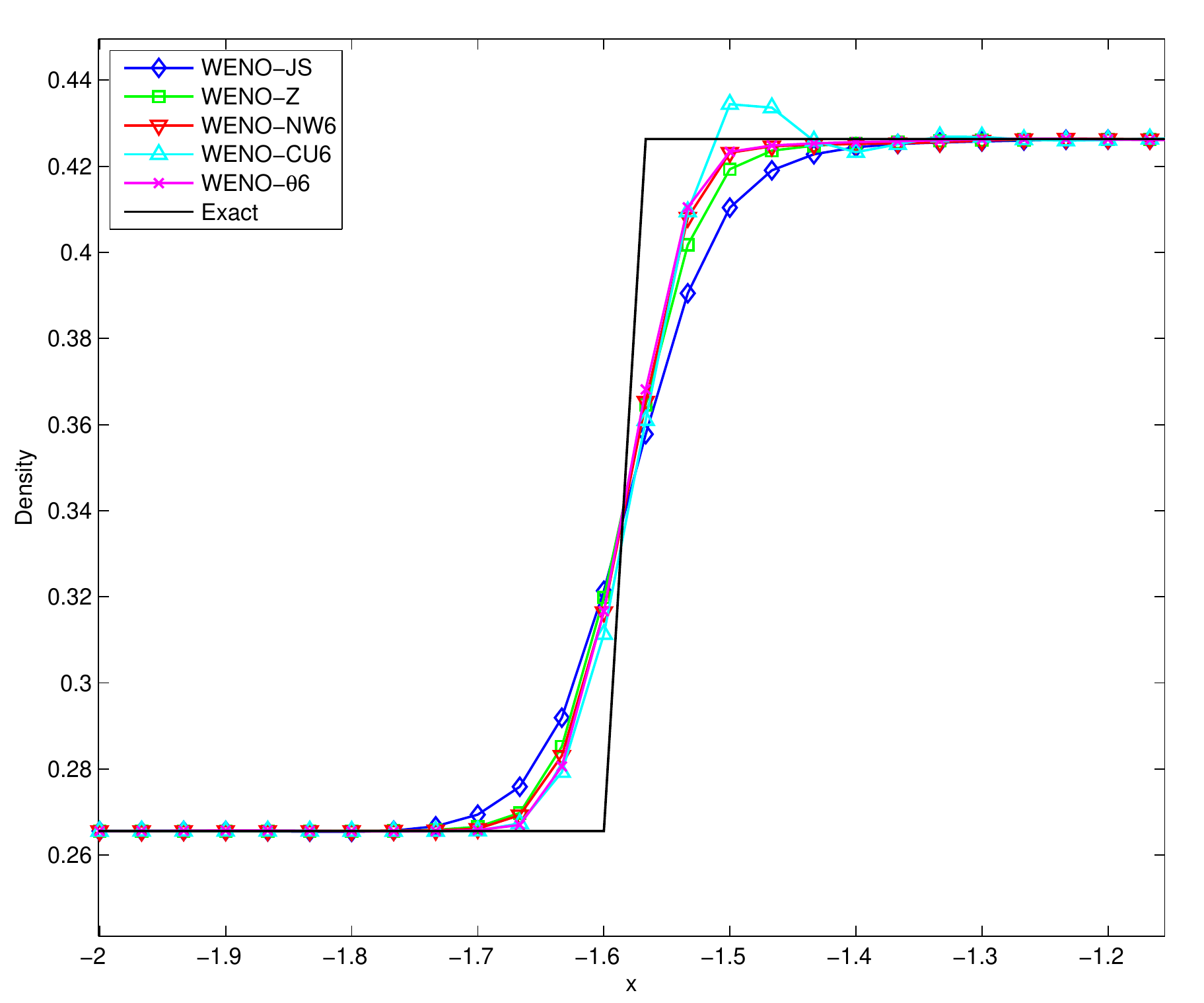}}
    \end{tabular}
    \caption{\small{Left: Sod's problem with initial data (\ref{sod_prob}). Time $t=1.7$.
    Grid $300$. Right: zoom at the contact discontinuity.}}
    \label{fig_soln_sod}
  \end{center}
\end{figure}

\begin{figure}
  \begin{center}
    \begin{tabular}{cc}
      \resizebox{70mm}{!}{\includegraphics{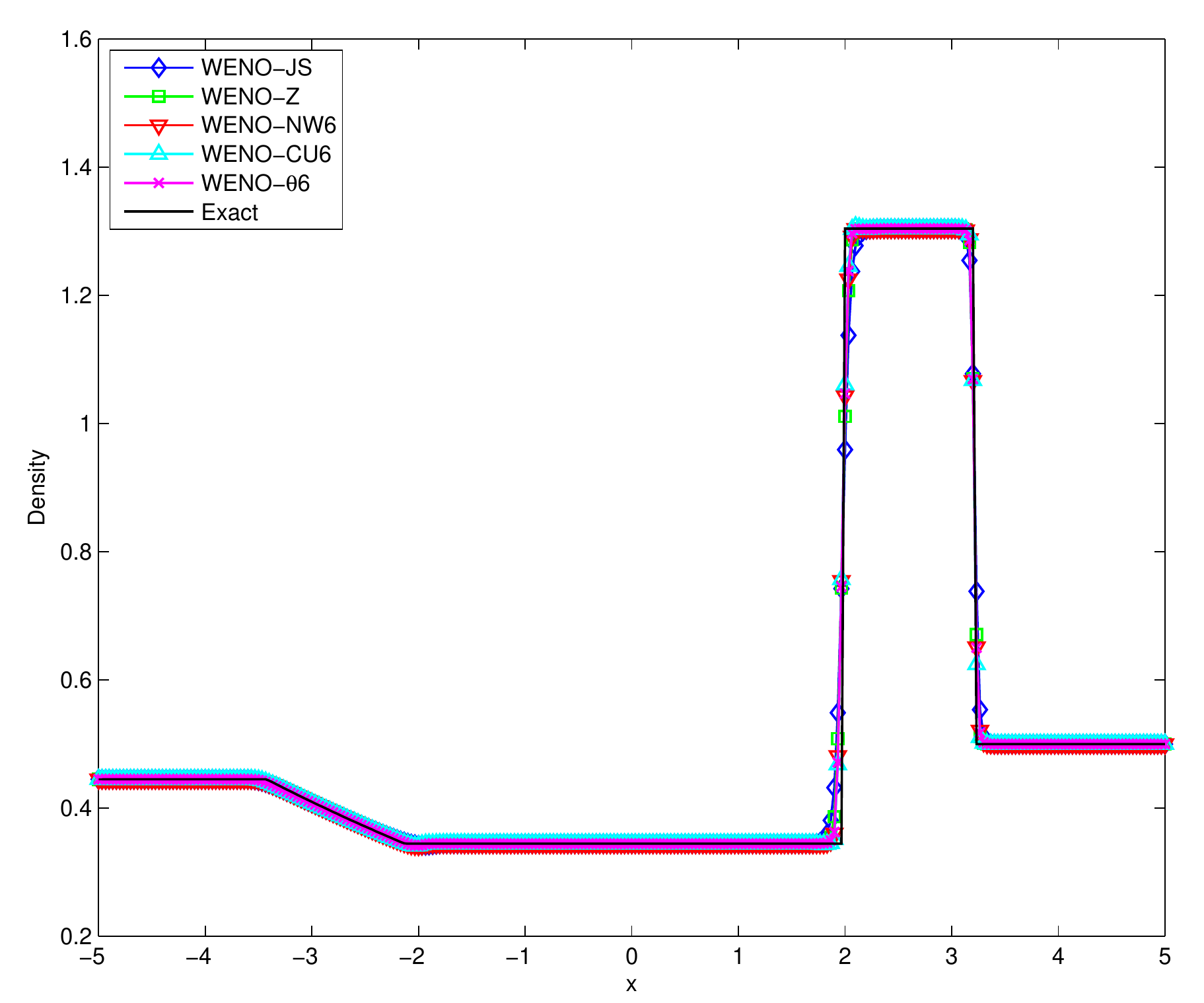}} &
      \resizebox{70mm}{!}{\includegraphics{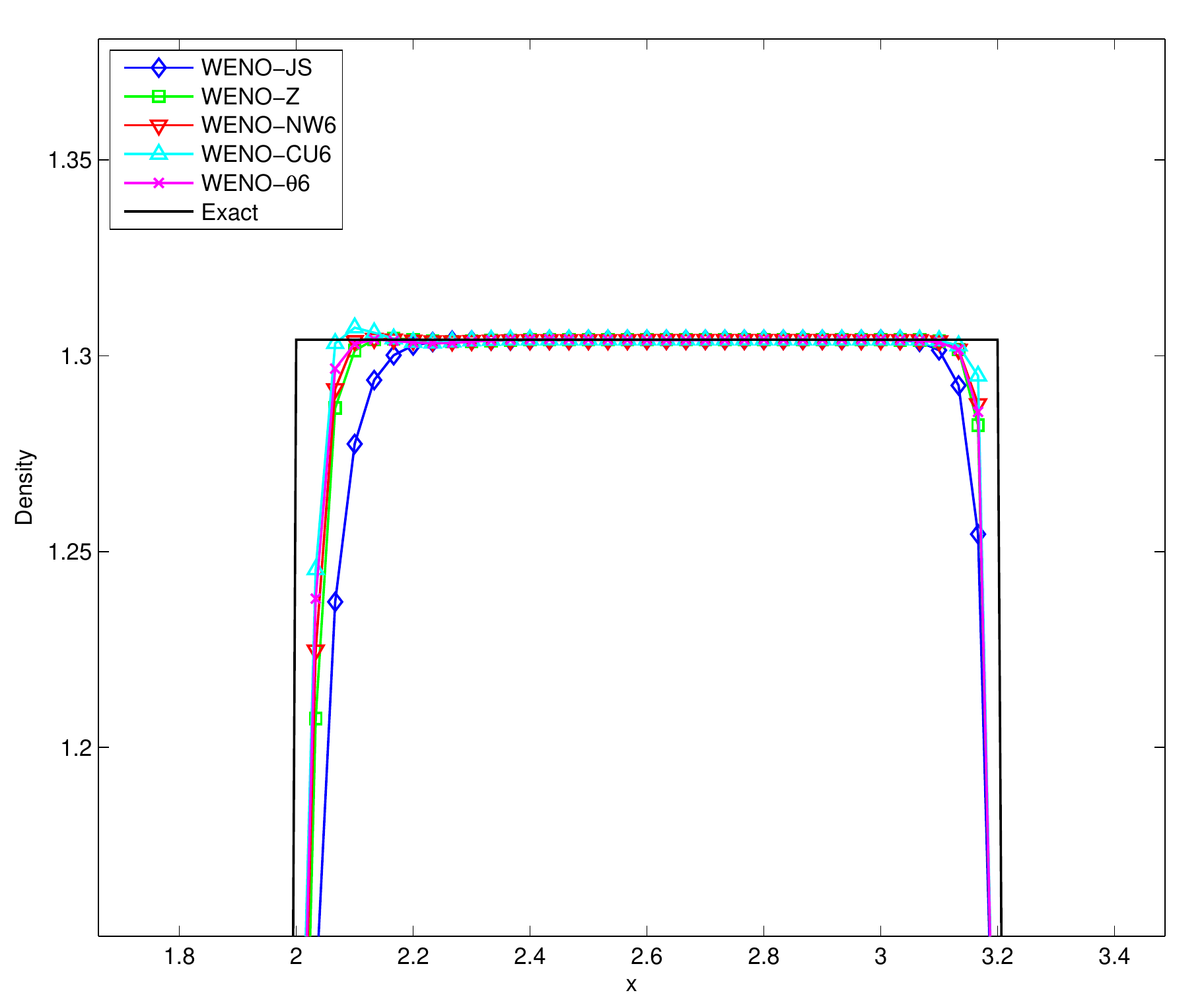}}
    \end{tabular}
    \caption{\small{Left: Lax's problem with initial data (\ref{lax_prob}). Time $t=1.3$.
    Grid $300$. Right: zoom at the contact discontinuity.}}
    \label{fig_soln_lax}
  \end{center}
\end{figure}

\begin{figure}
  \begin{center}
    \begin{tabular}{cc}
      \resizebox{70mm}{!}{\includegraphics{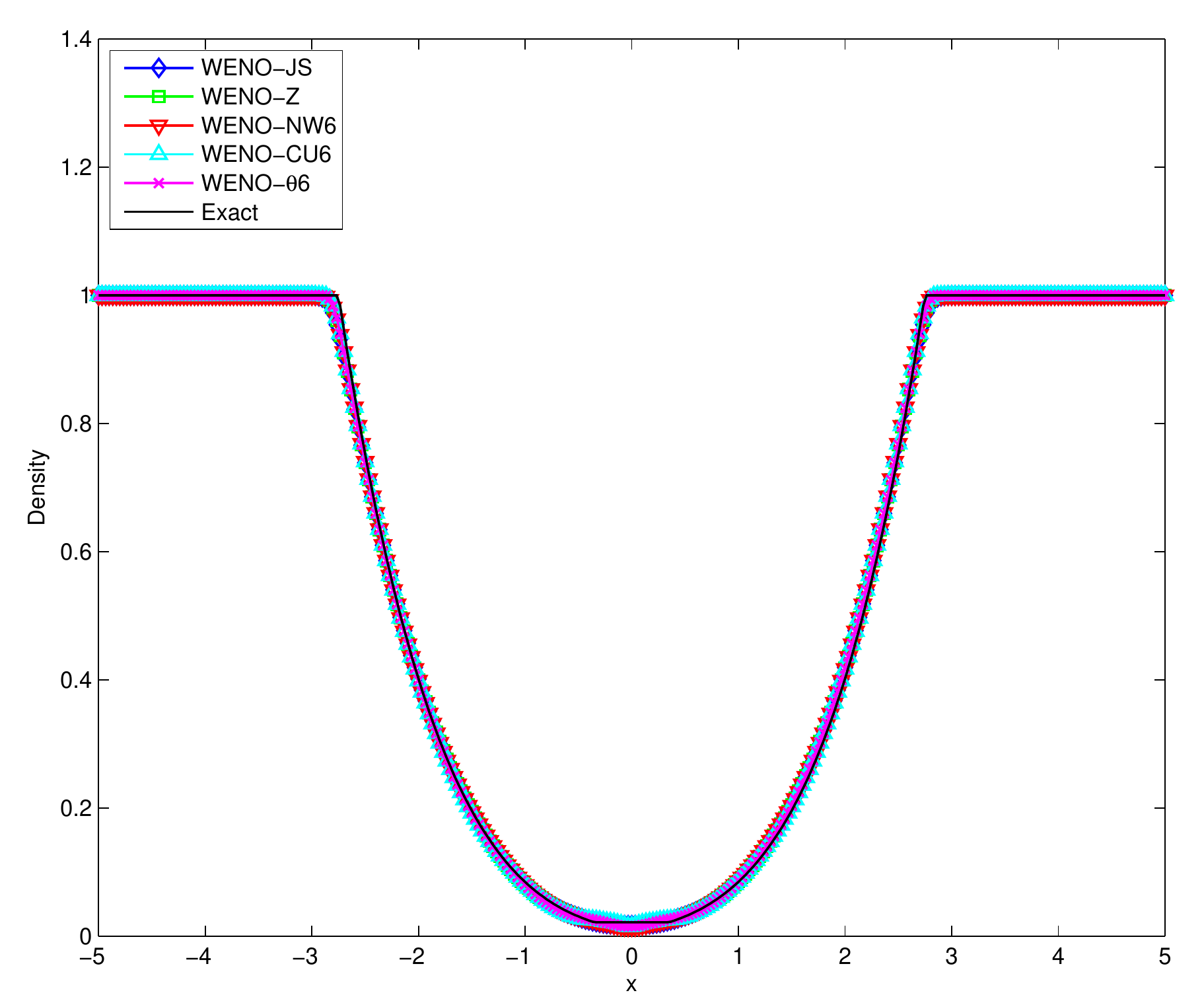}} &
      \resizebox{70mm}{!}{\includegraphics{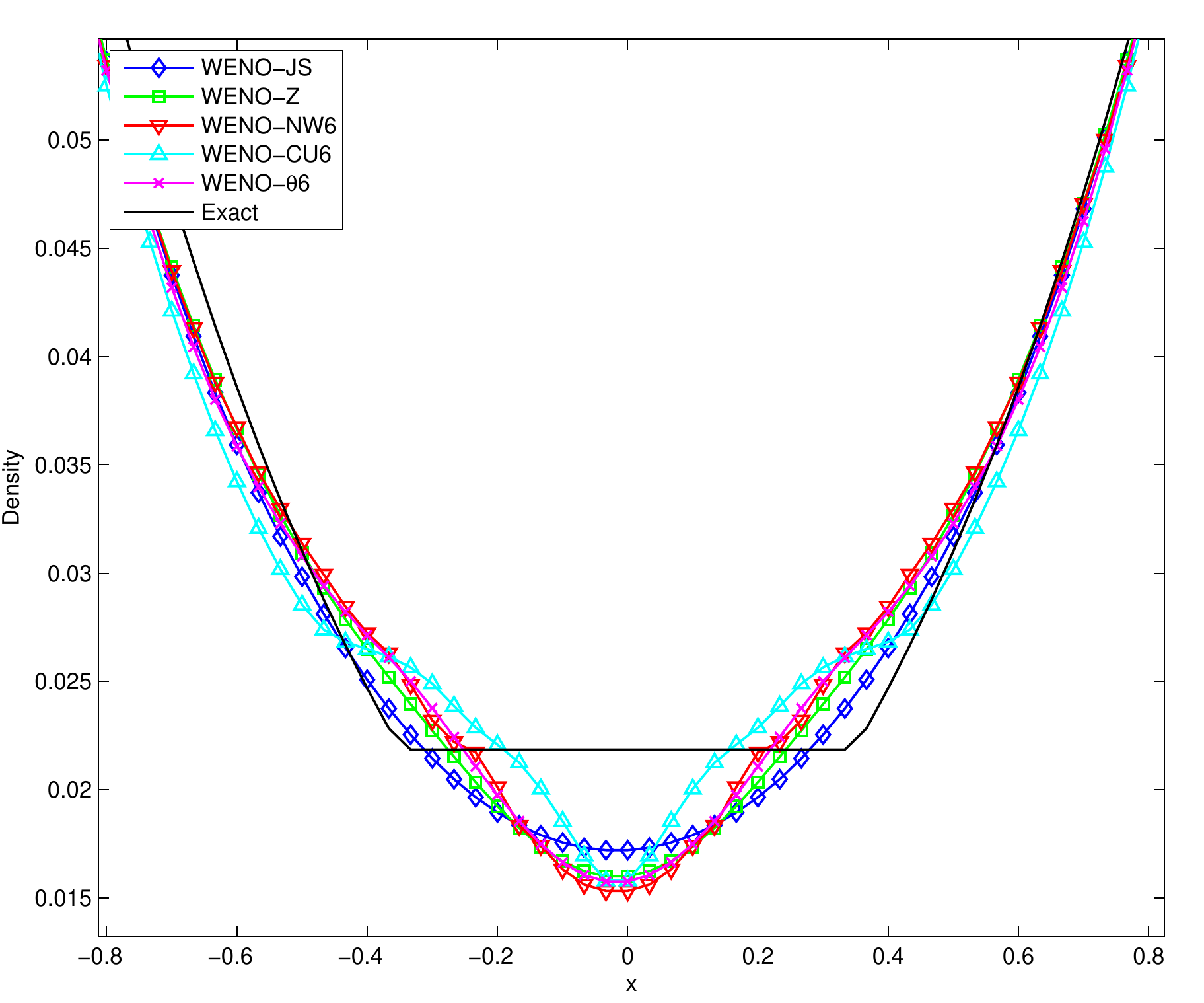}}
    \end{tabular}
    \caption{\small{Left: The $123$ problem with initial data (\ref{123_prob}).
    Time $t=1.0$. Grid $300$. Right: zoom at the trivial contact discontinuity.}}
    \label{fig_soln_123}
  \end{center}
\end{figure}

\subsubsection{TEST 4: Shock Density Wave Interaction, Shu-Osher's test}

We consider the following initial data,
\begin{align}\label{shock_den_prob}
(\rho,u,p)=
\begin{cases}
&(3.857143,\ \ 2.629369,\ \ 31/3),\quad -5<x<-4,\\
&(1+0.2\sin(5x),\ \ 0,\ \ 1),\quad -4<x<5,
\end{cases}
\end{align}
with zero-gradient boundary conditions.

The problem simulates the interaction of a right-moving Mach 3 shock
with a wavelike perturbed density whose magnitude is much smaller
than the shock. As a result, a flow field of compressed and
amplified wave trails is created right behind the shock. For more
details, see \cite{JS}. In Fig. \ref{fig_shock_den}, we show the
numerical results of the 6th-order WENO schemes at time $t=1.8$ with
grids of $N=200$ and $N=400$ intervals. The ``exact'' solution is
computed by WENO-JS with a fine grid $N=4000$. It is shown that all
schemes give satisfactory approximations of the compressed wavelike
structures behind the shock. A careful observation reveals that
WENO-$\theta6$ resolves the wave package as well as WENO-CU6,
whereas WENO-NW6 is more dissipative for both grid levels.

\begin{figure}
  \begin{center}
    \begin{tabular}{cc}
      \resizebox{70mm}{!}{\includegraphics{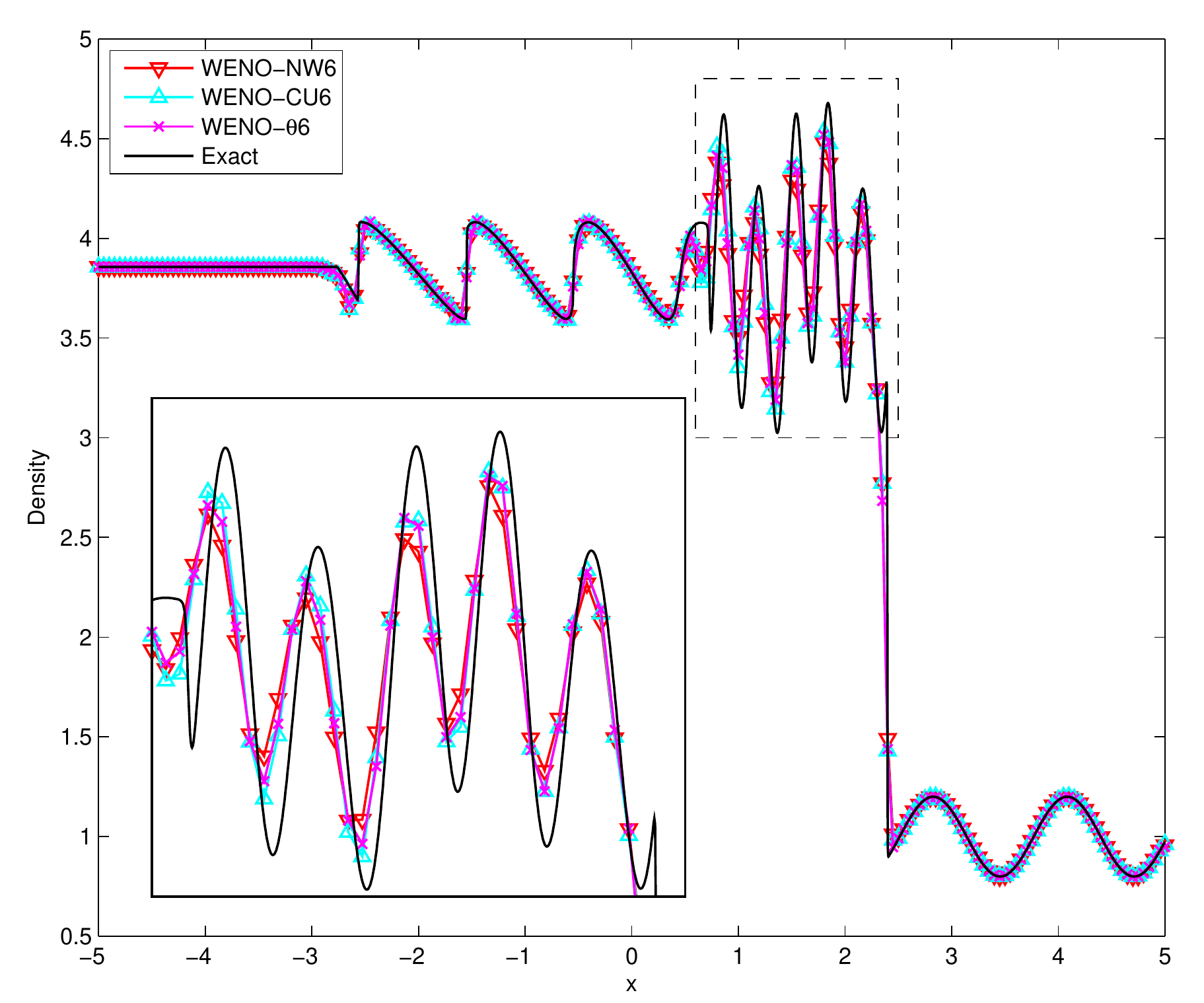}} &
      \resizebox{70mm}{!}{\includegraphics{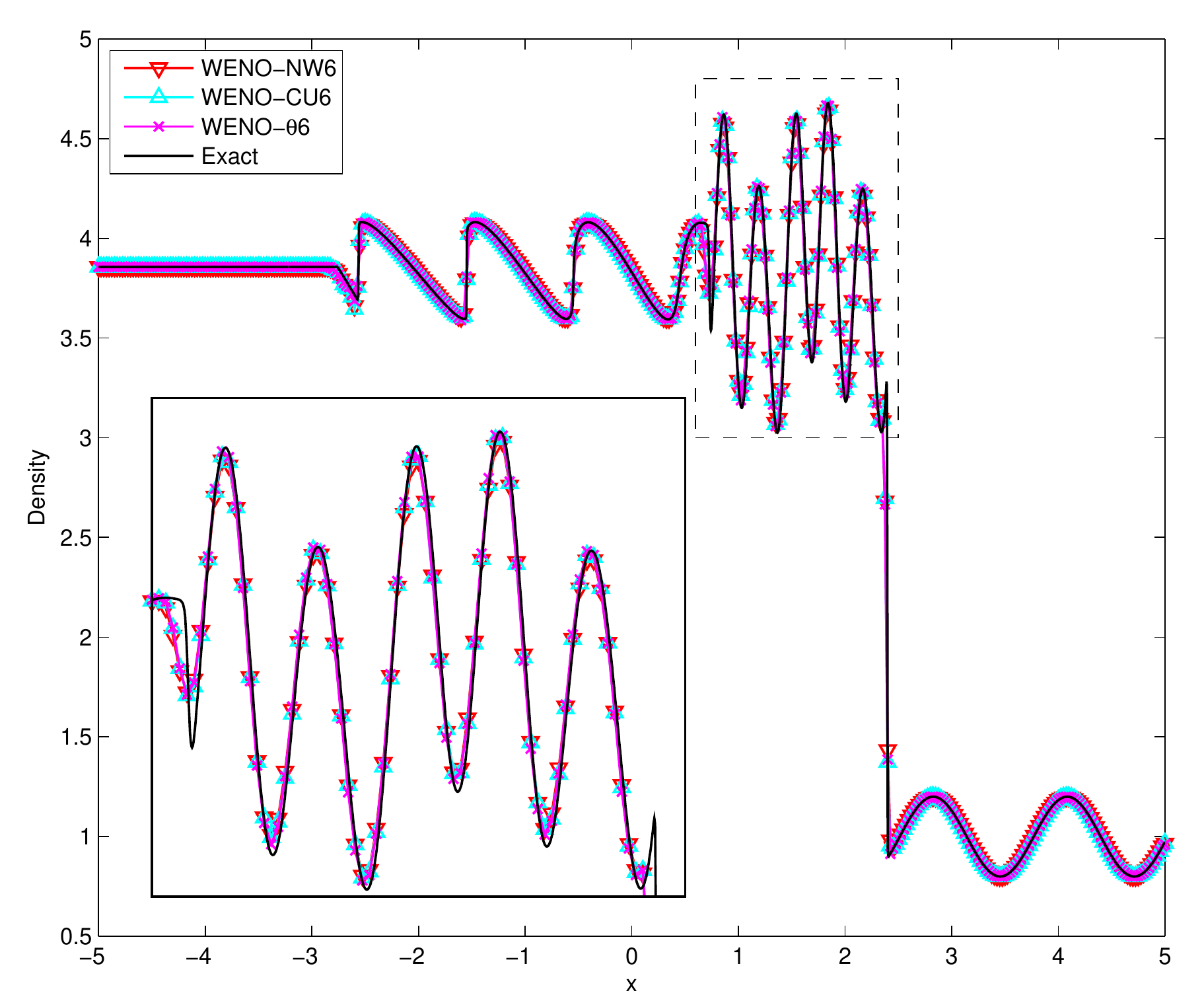}}
    \end{tabular}
    \caption{\small{Shu-Osher's problem with initial data (\ref{shock_den_prob}).
    Density. Time $t=1.8$. Left: medium grid $200$; Right: fine grid $400$.}}
    \label{fig_shock_den}
  \end{center}
\end{figure}

\subsubsection{TEST 5: Two Interacting Blast Waves}

In this test, we show that our new scheme WENO-$\theta6$ passes the
tough test of two interacting blast waves which the initial data are
given as follows,
\begin{align}\label{blast_prob}
(\rho,u,p)=
\begin{cases}
&(1,\ \ 0,\ \ 1000),\quad 0<x<0.1,\\
&(1,\ \ 0,\ \ 0.01),\quad 0.1<x<0.9,\\
&(1,\ \ 0,\ \ 100),\quad 0.9<x<1,
\end{cases}
\end{align}
and a reflective condition is applied at both boundaries. This
problem is used to test the robustness of shock-capturing methods
since many interactions are observed in a small area. A detailed
discussion of this problem can be found in \cite{WC84}.

Numerical results of 6th-order WENO schemes are computed up to time
$t=0.038$ with a grid of $N=801$ and plotted in Fig.
\ref{fig_blast_prob} for the density. The exact solution is
approximated by WENO-JS with a much fine grid $N=4001$. It is shown
that all schemes well capture the shocks as well as contact
discontinuities. A zoom near $x=0.745$ indicates that WENO-$\theta6$
gives better resolution than WENO-NW6 and WENO-CU6. We also
emphasize that there exists a stair-casing phenomenon in the
solutions of the latter methods in this region, which is similar to
that at the top of the semi-ellipse in test $1$ (see Fig.
\ref{fig_soln_lin}), and the $123$ problem (see Fig.
\ref{fig_soln_123}).

\begin{figure}
  \begin{center}
    \begin{tabular}{c}
      \resizebox{110mm}{!}{\includegraphics{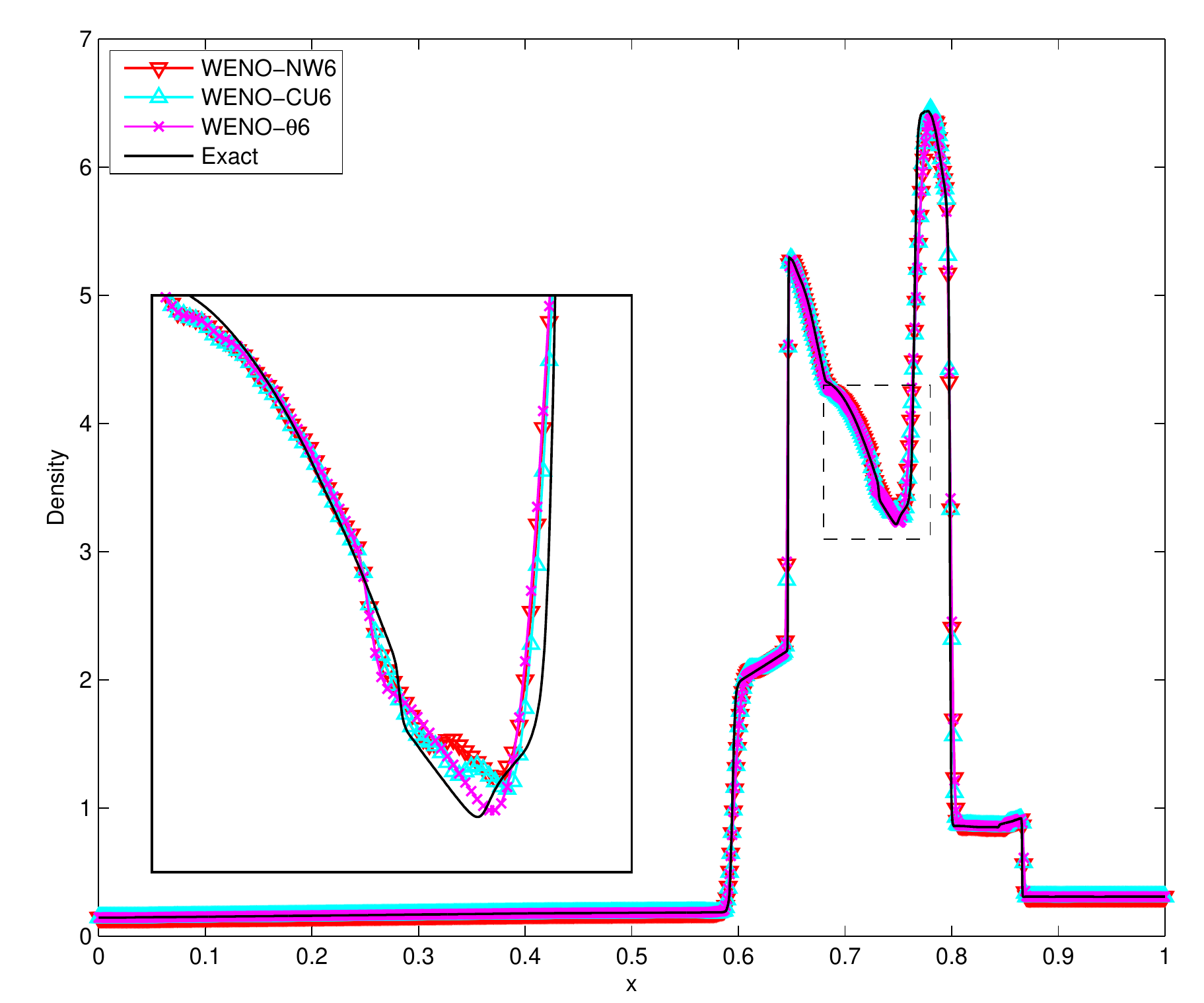}}
    \end{tabular}
    \caption{\small{Two interacting blast waves with initial data
    (\ref{blast_prob}). Density. Time $t=0.038$. Grid $801$.}}
    \label{fig_blast_prob}
  \end{center}
\end{figure}

\subsection{Two-dimensional Euler's Equations}

In this subsection, we extend the problem to two-dimensional cases.
We choose the 2D Euler equations which are as follows,
\begin{align}\label{eu2D}
\begin{cases}
&\mathbf{u}_t+\mathbf{f}(\mathbf{u})_x+\mathbf{g}(\mathbf{u})_y=0,\\
&\mathbf{u}(x,y,0)=\mathbf{u}_0(x,y),
\end{cases}
\end{align}
where $\mathbf{u}=(\rho,\rho u,\rho v, E)^T$,
$\mathbf{f}(\mathbf{u})=(\rho u,p+\rho u^2,\rho uv,u(E+p))^T$,
$\mathbf{g}(\mathbf{u})=(\rho v,\rho uv,p+\rho v^2,v(E+p))^T$. The
relation of pressure and conservative quantities is through the
equation of state
\begin{align}
p=(\gamma-1)\left(E-\dfrac{1}{2}(u^2+v^2)\right).
\end{align}
Here, we choose the ratio of specific heats $\gamma=1.4$.

\subsubsection{TEST 6: Rayleigh-Taylor Instability}

In the following tests, we show numerical evidence that
WENO-$\theta6$ maintains symmetry in the solutions much better than
the other 6th-order schemes, and outperforms 5th-order schemes in
resolving small-scaled structures occurring in flow configurations.
We first simulate the Rayleigh-Taylor instability. The instability
occurs where there is a heavy fluid falling into a light fluid (see
\cite{GGLO}, \cite{At}, \cite{FSTY}). Following \cite{At}, we set up
the problem as follows. The domain is
$(x,y)=(-0.25,0.25)\times(-0.75,0.75)$. Initial density has a
discontinuity at the interface, i.e., $\rho=2$ for $y\ge0$ and
$\rho=1$ for $y<0$. The pressure is set at hydrostatic equilibrium
initially $p=2.5-\rho gy$ where $g=0.1$ is the gravitational
acceleration. The $x$-component velocity $u=0$, while the
$y$-component is perturbed with $v=\frac{0.01}{4}(1+\cos(4\pi
x))(1+\cos(\frac{4}{3}\pi y))$ for a single mode perturbation.
Boundary conditions are set periodic in $x$-direction, and
reflective in $y$-direction. The ratio of specific heats
$\gamma=1.4$. We add $-g\rho$ and $-g\rho v$ in the $y$-momentum and
energy equations of (\ref{eu2D}) as source terms.

\begin{figure}
  \begin{center}
    \begin{tabular}{c}
      \resizebox{150mm}{!}{\includegraphics{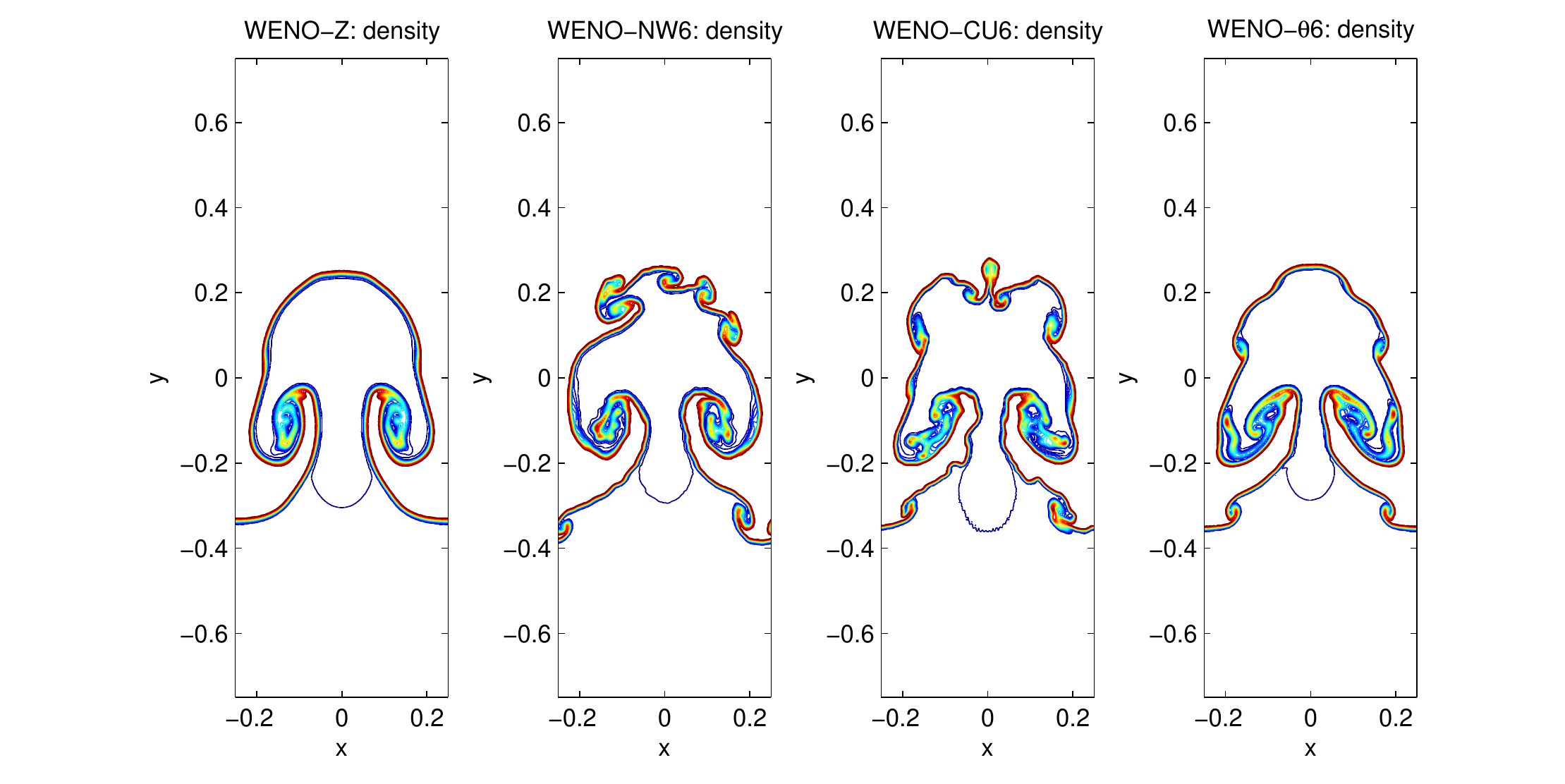}}
    \end{tabular}
    \caption{\small{The Rayleigh-Taylor instability.
    Density at time $t=9.5$. Grid $120\times360$. $CFL=0.5$.
    From left to right: WENO-Z, WENO-NW6, WENO-CU6, and WENO-$\theta6$.}}
    \label{fig_RT}
  \end{center}
\end{figure}

In Fig. \ref{fig_RT}, we plot the density with 20 equally spaced
contours obtained from 5th- and 6th-order WENO schemes at time
$t=9.5$ with a $120\times360$ grid. It is shown that the 6th-order
schemes have much better numerical resolution comparing with the
5th-order ones. We notice that WENO-$\theta6$ preserves the symmetry
of the solution; whereas WENO-NW6 and WENO-CU6 do not. We conjecture
the lack of symmetry of WENO-NW6 is due to the loss of accuracy
around critical regions which is shown in previous numerical tests.
The test also shows that the discontinuous switching of
$\tau^{\theta}$ in Eq. (\ref{tau_we6}) does not affect the
robustness of the new WENO-$\theta$ scheme, even for problem with
highly unstable fluid flows as the Rayleigh-Taylor instability.

\subsubsection{TEST 7: Implosion problem}

The next numerical test is the implosion problem (see \cite{At},
\cite{LW}) with initial data as follows,
\begin{align}\label{eu2D_im}
(\rho,p)=\begin{cases} &(1,1)\quad\text{for }x+y>\frac{1}{2},\\
&(0.125,0.14)\quad\text{otherwise},
\end{cases}
\end{align}
and zero velocity everywhere initially. We choose reflecting
conditions for all boundaries.

Symmetry is important for this test. For such a scheme, due to the
interactions of shock waves and reflecting boundaries, jets along
the diagonal are created. Longer and narrower jets are produced for
less dissipative schemes.

In Fig. \ref{fig_eu2D_im}, we show the results obtained from
different schemes on computational domain $(0,1)\times(0,1)$ at
final time $t=5$. We choose a grid of $400\times400$. For
WENO-$\theta$, we choose $\alpha_R=1$. It is shown that only WENO-Z
and WENO-$\theta$ well preserve the symmetry of the problem; whereas
the other 6th-order schemes do not. The jets created by WENO-NW6 and
WENO-CU6 tend to diverge from the main diagonal $x=y$. We also note
that the jets produced by WENO-$\theta$ is much longer and narrower
than those of WENO-Z, which means that the former scheme is less
dissipative than the latter one.

\begin{figure}
  \begin{center}
    \begin{tabular}{cc}
      \resizebox{70mm}{!}{\includegraphics{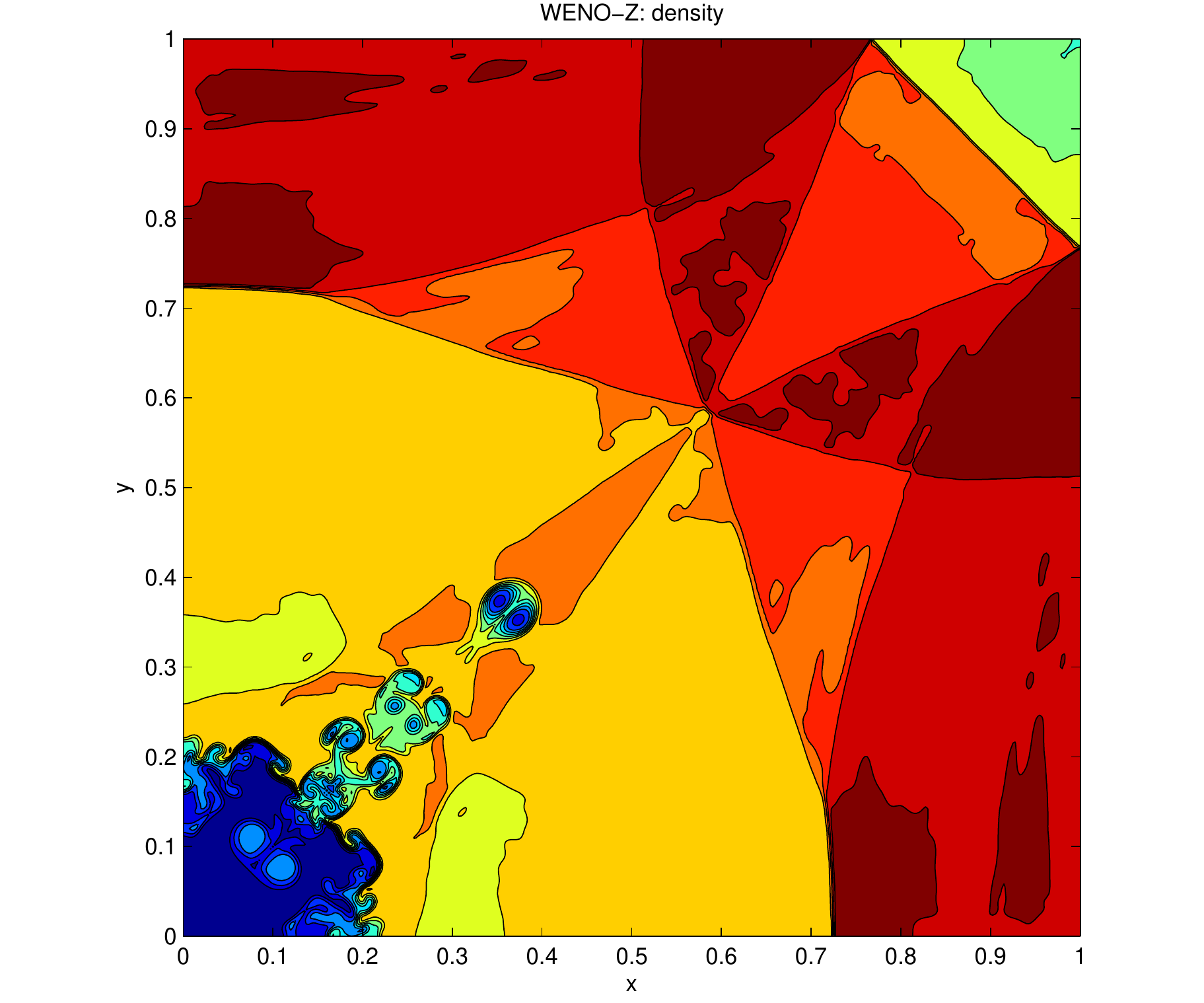}}&
      \resizebox{70mm}{!}{\includegraphics{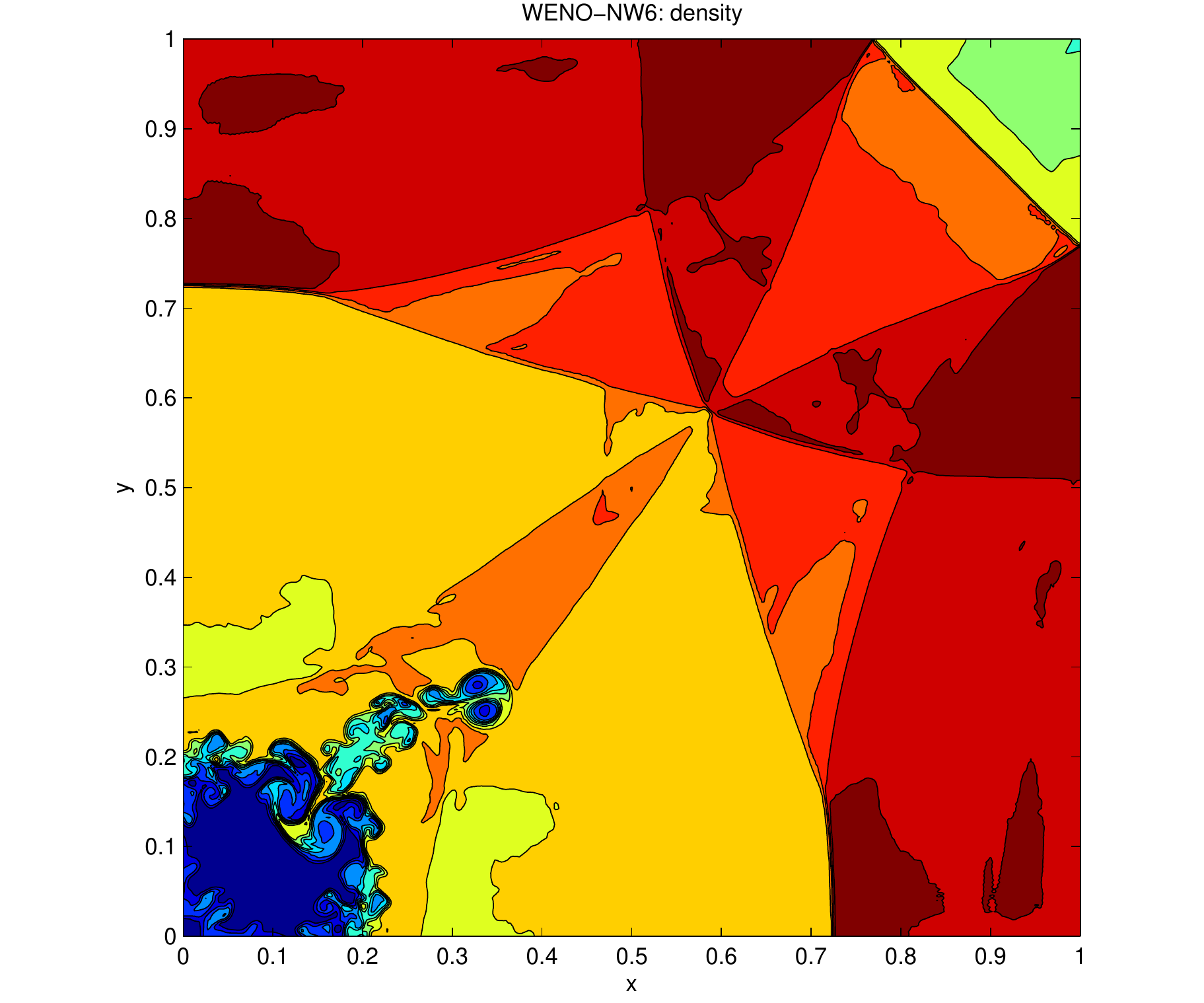}}\\
      \resizebox{70mm}{!}{\includegraphics{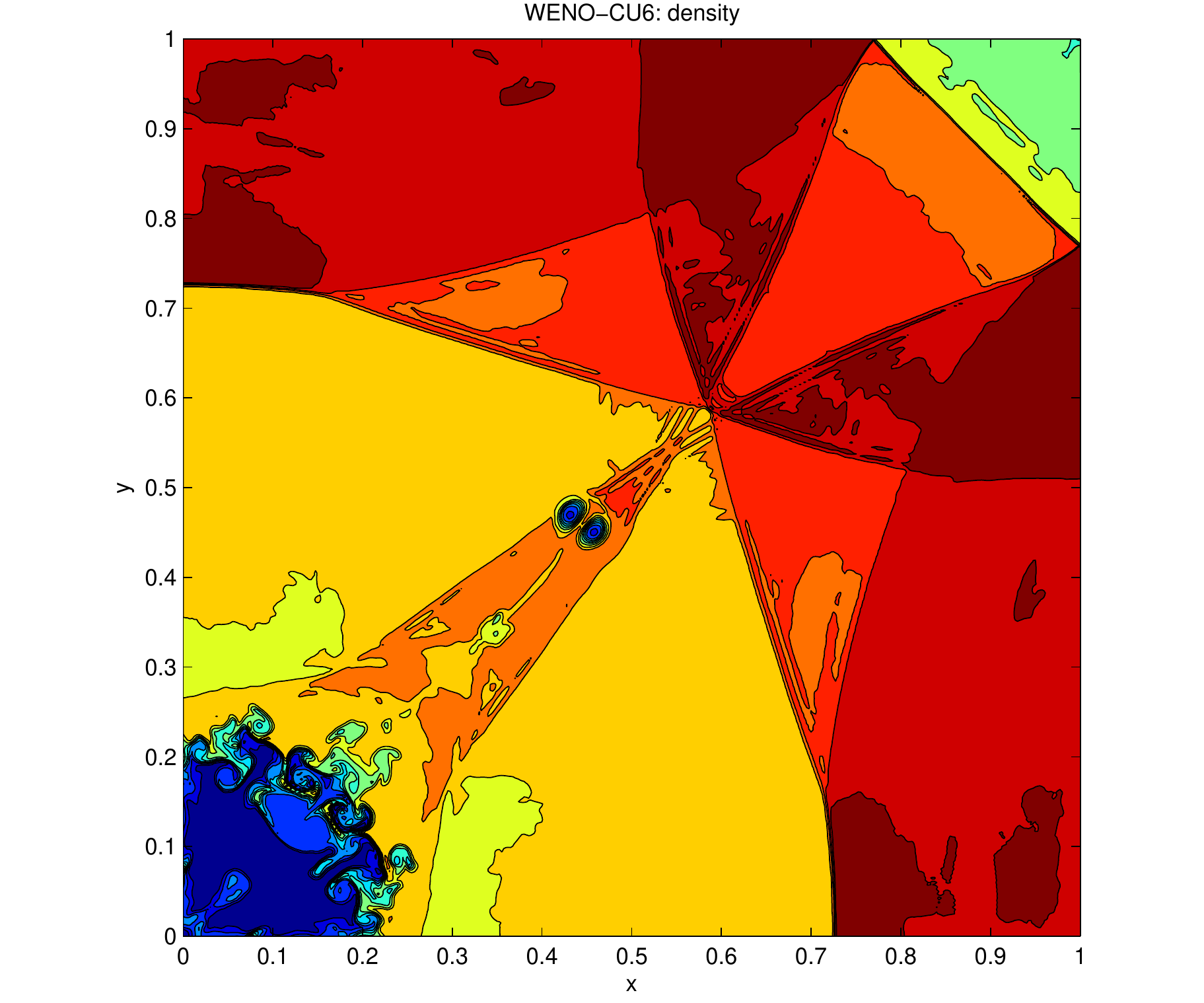}}&
      \resizebox{70mm}{!}{\includegraphics{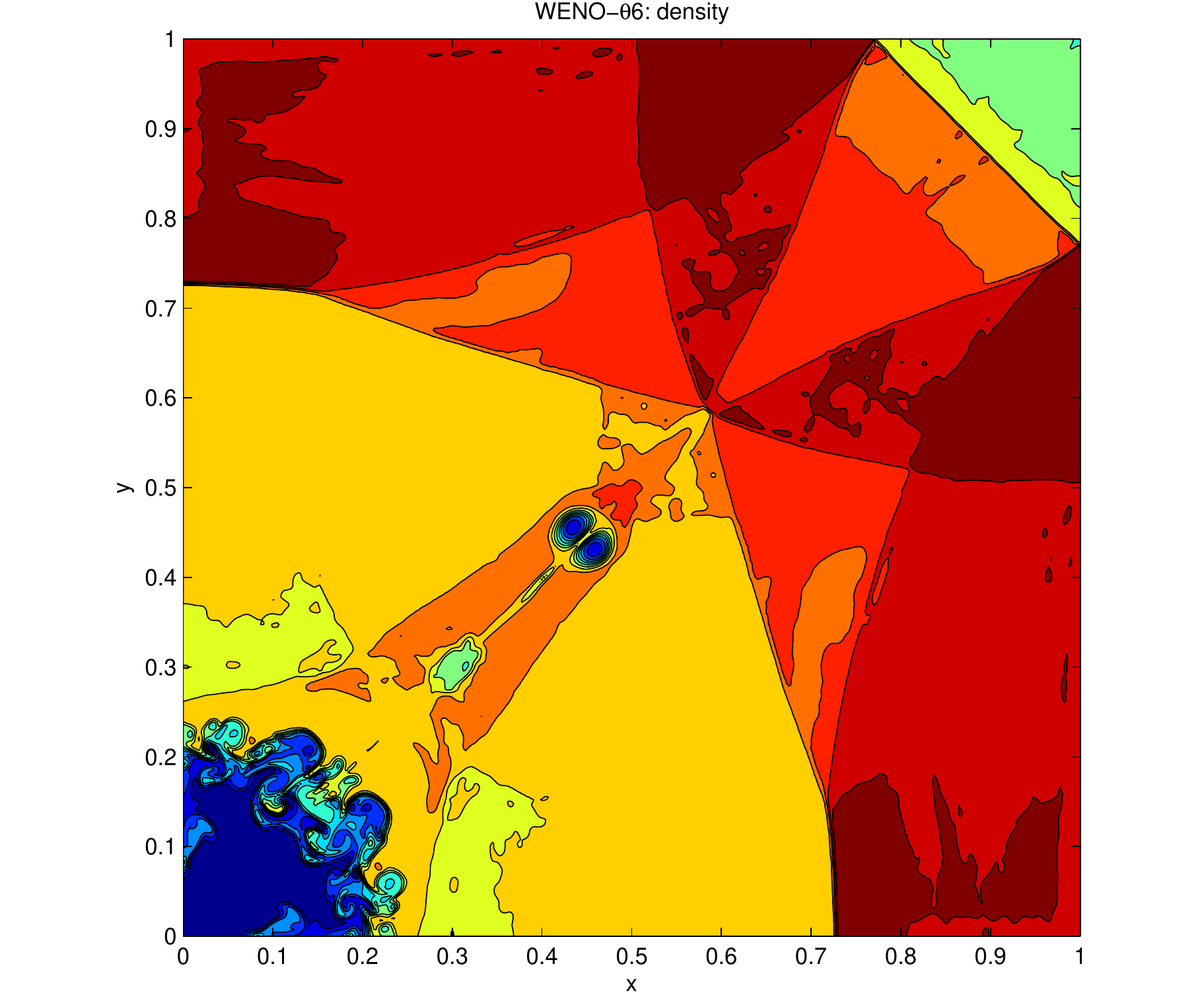}}
    \end{tabular}
    \caption{\small{ The implosion problem. Density with 20 contours uniformly distributing from $0$ to 1. Grid $400\times400$. Final
    time $t=5$.
    Left to right, top to bottom: WENO-Z, WENO-NW6, WENO-CU6, WENO-$\theta$6.}}
    \label{fig_eu2D_im}
  \end{center}
\end{figure}

\subsubsection{TEST 8: 2D Riemann Initial Data}

The 2D Riemann problem is set up by assigning different constant
states of $(\rho_k,u_k,v_k,p_k)$, $k=1,2,3,4$, to four quadrants of
the computational domain $\Omega=(0,1)\times(0,1)$. The constant
states are chosen so that there is only a single elementary wave,
namely, shock-, rarefaction-, and contact-wave, connecting two
neighboring quadrants (see \cite{SCG}). For our test, we choose the
following configuration for the initial data, respectively, for
quadrants $1$, $2$, $3$, $4$,
\begin{align}\label{eu2D_Rie_ini12}
(\rho,u,v,p)=
\begin{cases}
&(0.5313,\ 0,\ 0,\ 0.4),\ \quad x>0.5,\ y>0.5,\\
&(1,\ 0.7276,\ 0,\ 1),\quad\quad x<0.5,\ y>0.5,\\
&(0.8,\ 0,\ 0,\ 1),\ \ \quad\quad\quad x<0.5,\ y<0.5,\\
&(1,\ 0,\ 0.7276,\ 1),\ \ \ \quad x>0.5,\ y<0.5,
\end{cases}
\end{align}
which has shocks through quadrants $1$ - $2$ and $1$ - $4$, and
contact discontinuities through quadrants $2$ - $3$ and $3$ - $4$.
Transmissive boundary conditions are imposed on all boundaries for
these two cases.

\begin{figure}
  \begin{center}
    \begin{tabular}{cc}
      \resizebox{70mm}{!}{\includegraphics{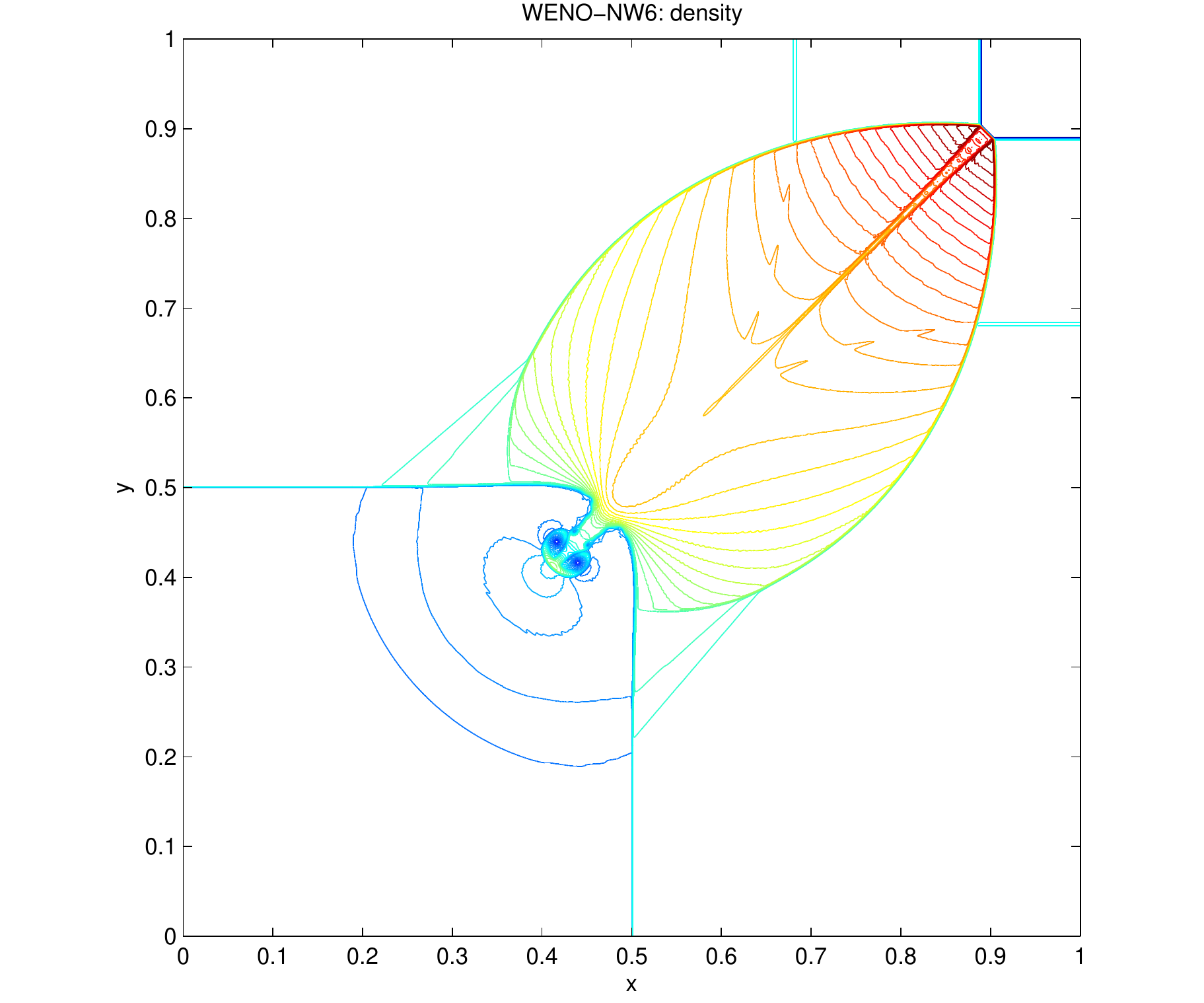}}
     &\resizebox{70mm}{!}{\includegraphics{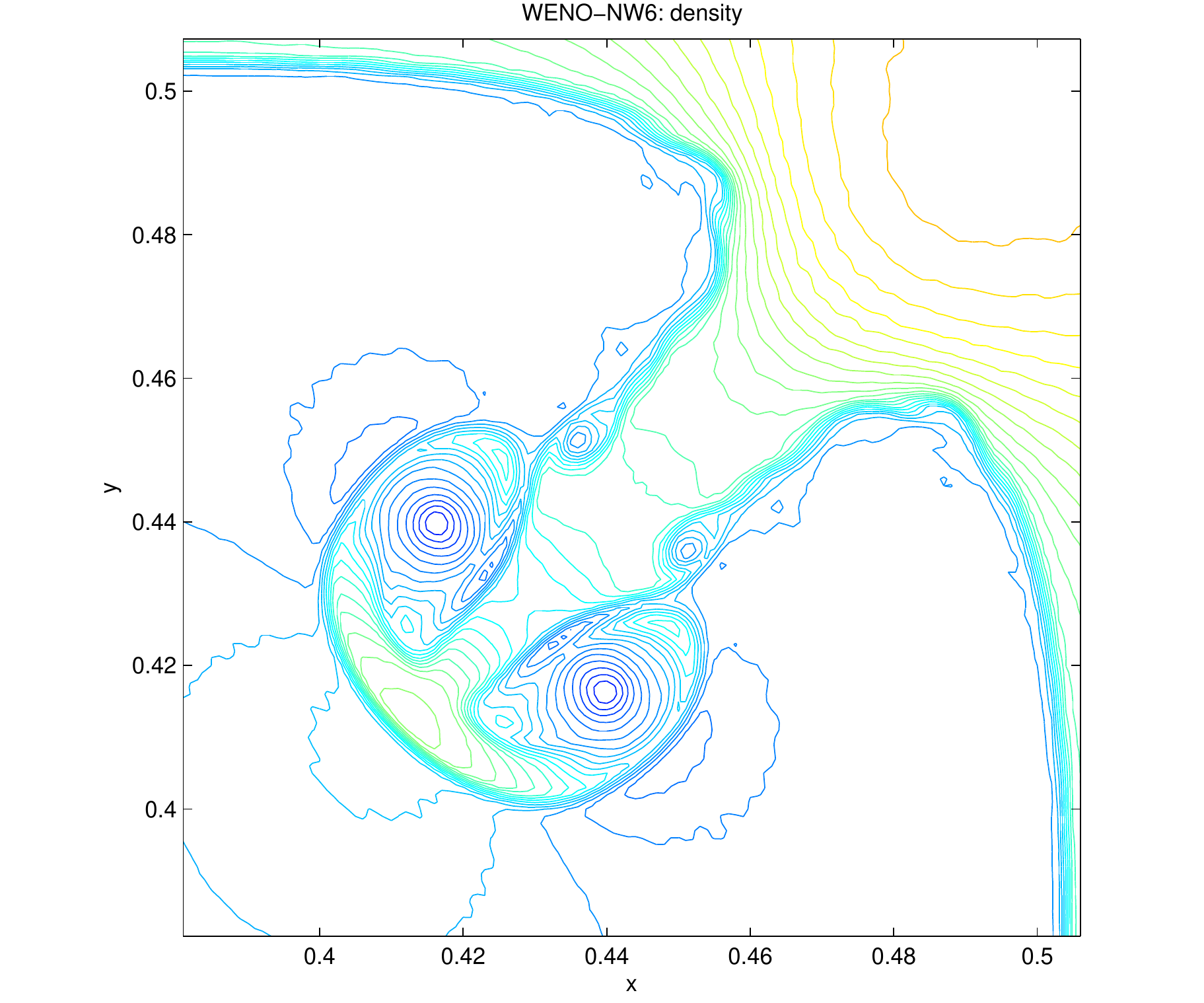}}\\
      \resizebox{70mm}{!}{\includegraphics{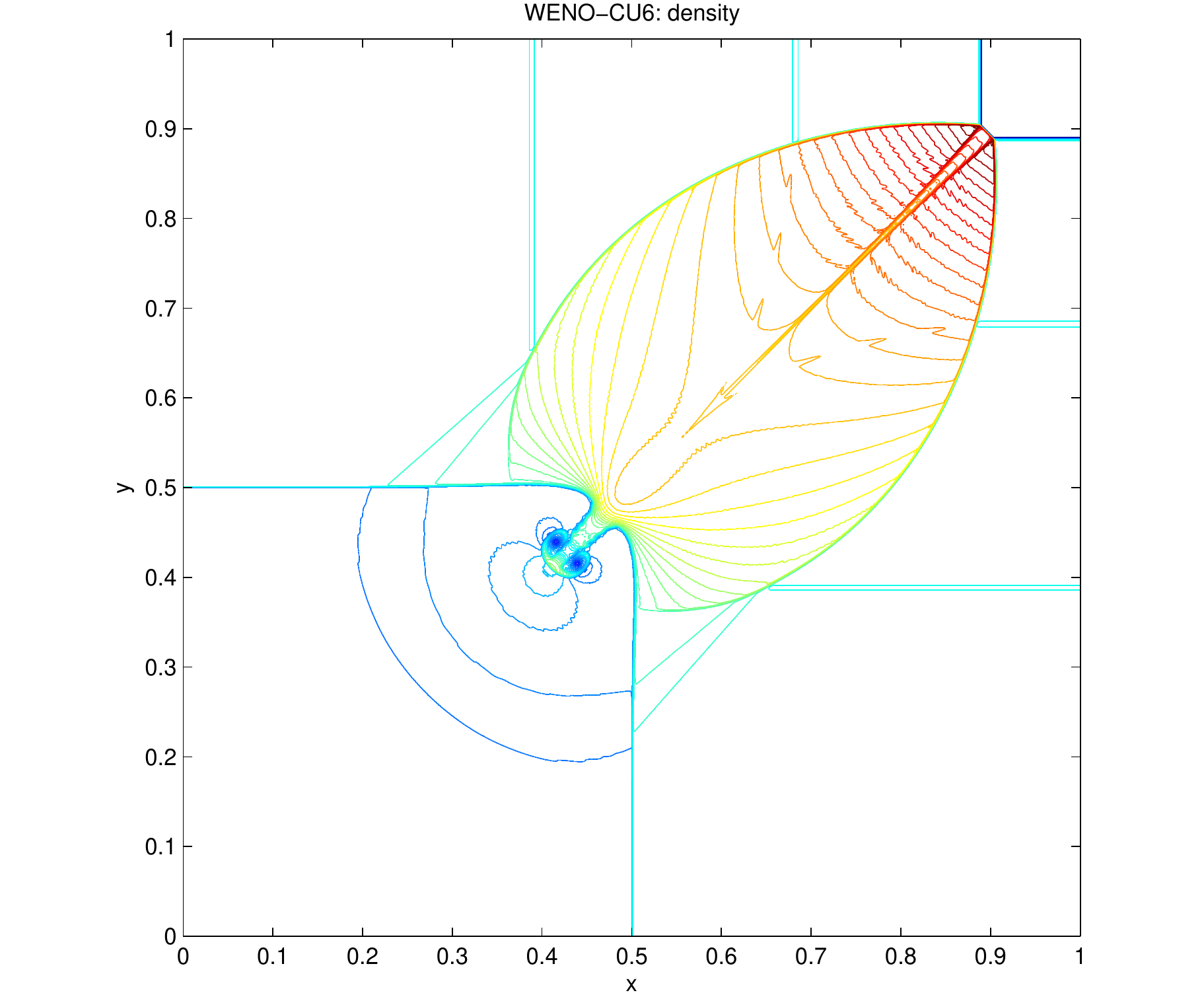}}
     &\resizebox{70mm}{!}{\includegraphics{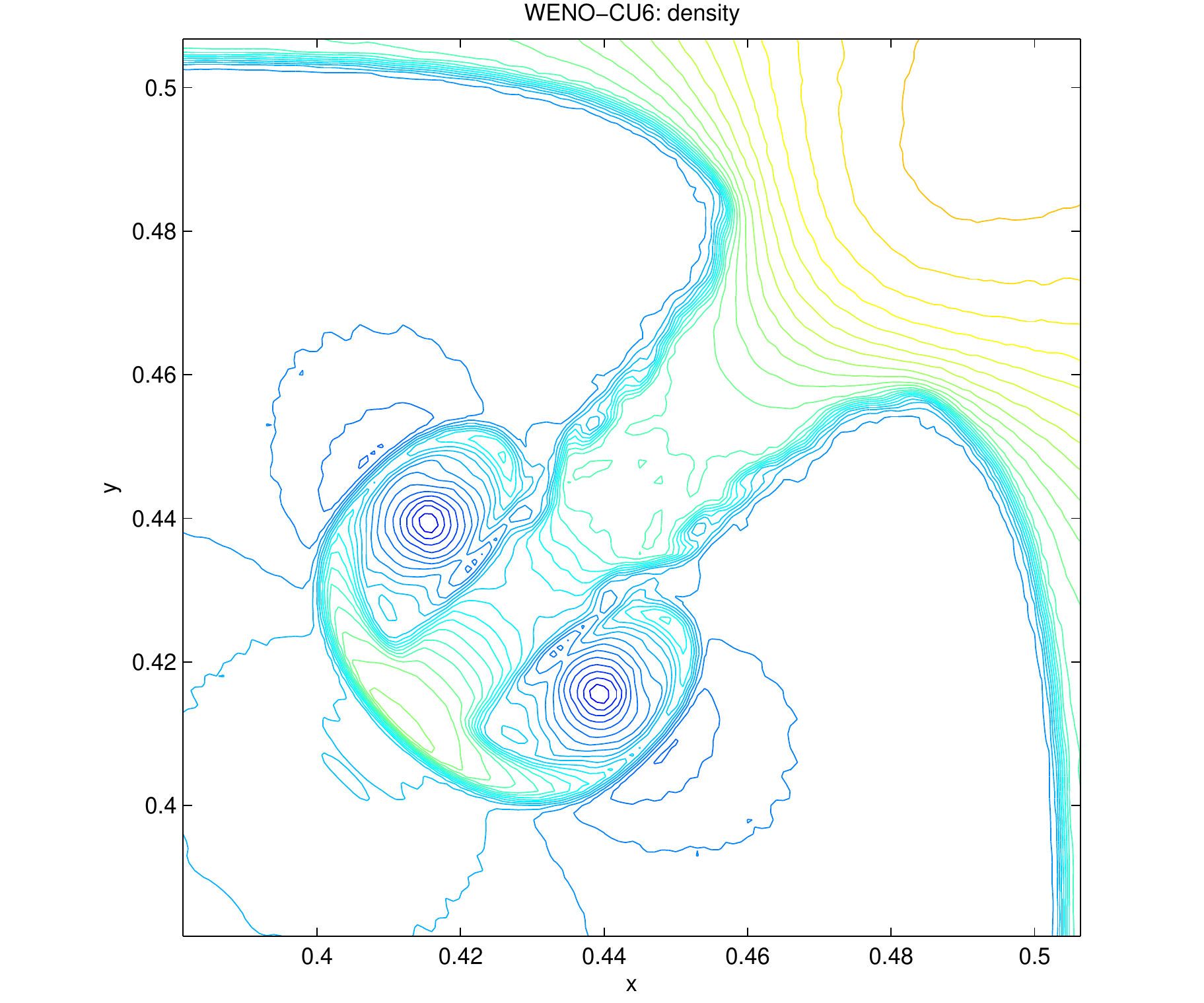}}\\
      \resizebox{70mm}{!}{\includegraphics{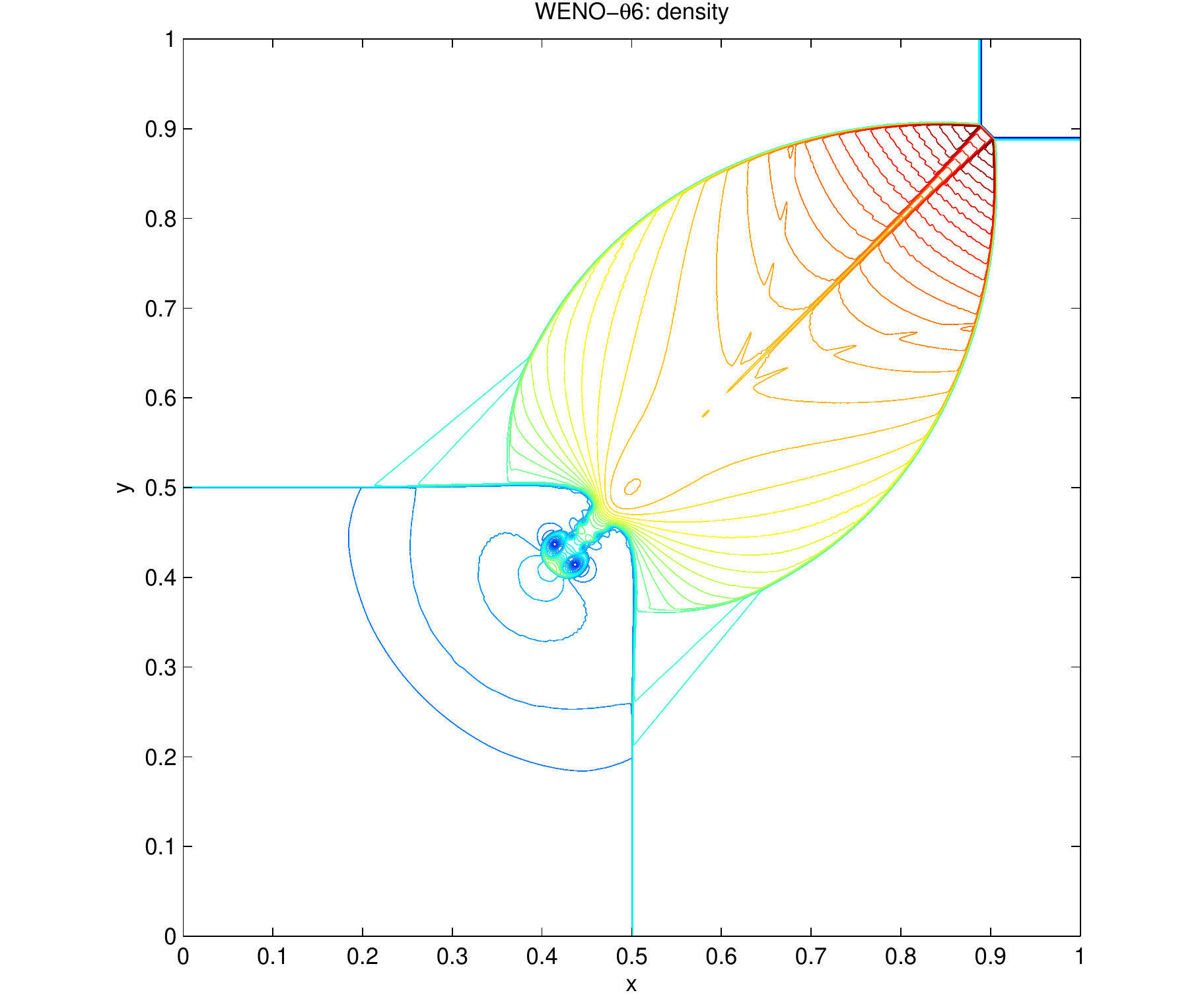}}
     &\resizebox{70mm}{!}{\includegraphics{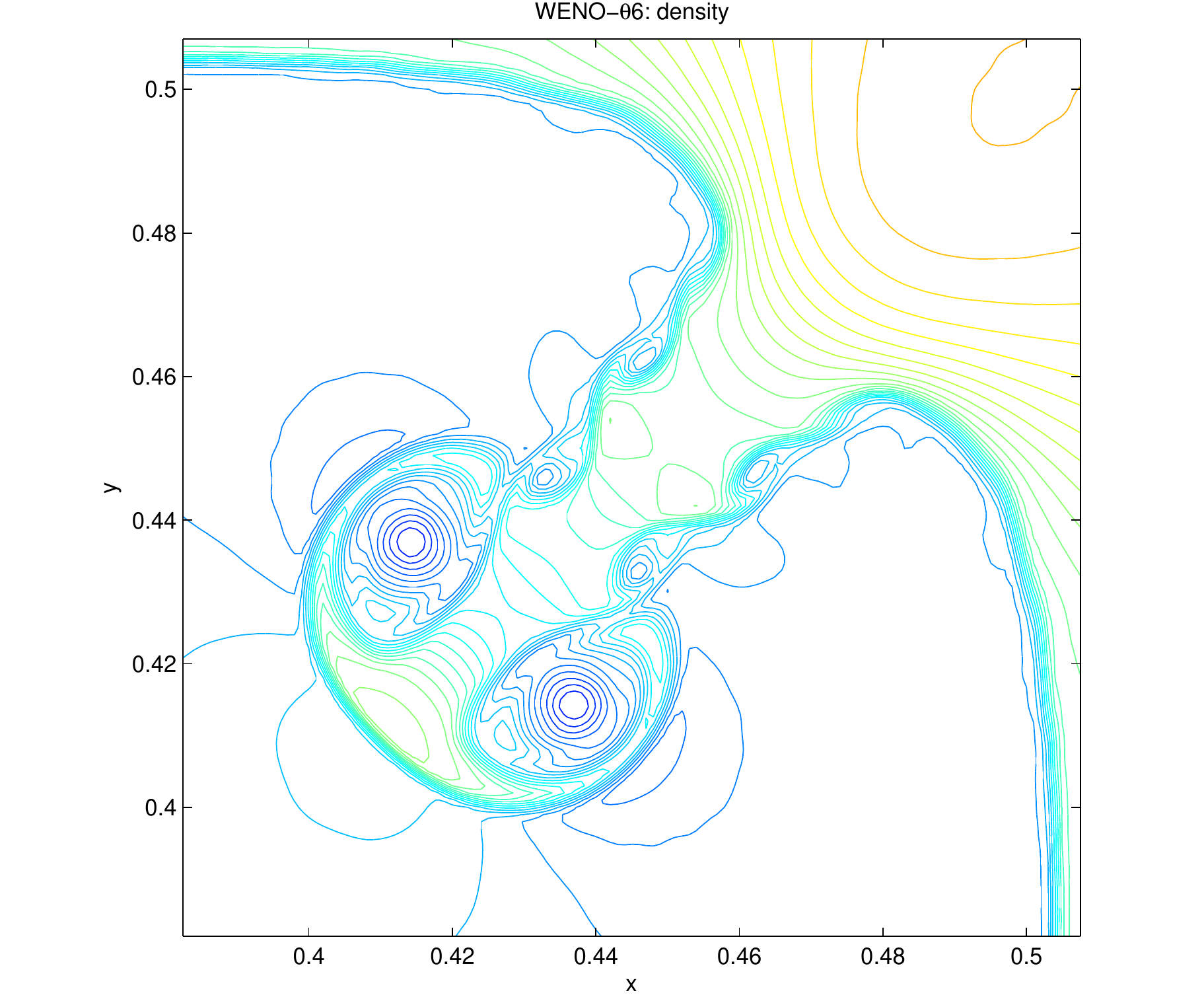}}
    \end{tabular}
    \caption{\small{Left: The 2D Riemann problem with initial data (\ref{eu2D_Rie_ini12}).
    Density with $50$ contours. Time $t=0.25$. Grid $1000\times1000$.
    Right: Zoom at the spirals.
    From top to bottom, respectively: WENO-NW6, WENO-CU6, WENO-$\theta$6.}}
    \label{fig_eu2D_ini12}
  \end{center}
\end{figure}

The approximations of the density with initial data
(\ref{eu2D_Rie_ini12}) at time $t=0.25$ are plotted in Fig.
\ref{fig_eu2D_ini12} with $50$ contours for WENO-NW6 and
WENO-$\theta$6. Here, we use a fine grid with $1000\times1000$
intervals for the capturing of the small vortices along the
contacts. Zooms near the spirals region are also shown on the right
column in the same figure. Again, we observe a better performance of
the WENO-$\theta$6 over WENO-NW6 and WENO-CU6 schemes over these
small-scaled structures without oscillations on the contours.

\subsubsection{TEST 9: Double Mach Reflection of a Strong Shock}

Finally, we investigate the double Mach reflection of a strong shock
which is a typical benchmark test for shock-capturing methods. The
problem simulates the reflection occurring when a simple planar
shock interacts with a wedge making with the $x$-axis an angle
$\alpha$. The strength of the moving shock is characterized by the
Mach number $M_s$. For a double Mach reflection problem, $M_s=10$
and the wedge angle is chosen as $\alpha=30^{\circ}$. Detailed
discussions on this type of problems can be found in \cite{WC84} and
the references therein. For numerical purpose, we choose the
computational domain $\Omega=(0,4)\times(0,1)$. Initially the shock
is located at $x_0=1/6$, inclined with the $x$-axis by the angle
$90^{\circ}-\alpha$. Inflow and zero gradients conditions are
imposed on the left and right boundaries, respectively. On the
bottom one, a reflective condition is applied to the interval
$[x_0,4]$ representing the wedge, and the exact post-shock state is
imposed over $[0,x_0]$. The top boundary is treated in a way that
there are no interactions of the shock with this boundary. That is,
the exact post- and pre-shock states are employed over the intervals
$[0,x_s(t)]$ and $[x_s(t),4]$, respectively, on the top boundary.
Here, $x_s(t)=x_0+\dfrac{1}{\tan
60^{\circ}}+\dfrac{M_sa_{pre}}{\cos30^{\circ}}t$, where $a_{pre}$ is
the sound speed of the pre-shock state, is the location of the shock
in time. These states can be computed exactly when one of them is
pre-described (see, e.g., \cite{To}). In particular, for our problem
the initial conditions are given as follows,
\begin{align}\label{eu2D_mach}
(\rho,u,v,p)=
\begin{cases}
&(8,\ 8.25\cos 30^{\circ},\ -8.25\sin 30^{\circ},\
116.5),\quad x<x_0+\frac{y}{\tan 60^{\circ}},\\
&(1.4,\ 0,\ 0,\
1),\quad\quad\quad\quad\quad\quad\quad\quad\quad\quad\quad\ x\ge
x_0+\frac{y}{\tan 60^{\circ}}.
\end{cases}
\end{align}

\begin{figure}
  \begin{center}
    \begin{tabular}{c}
      \resizebox{110mm}{!}{\includegraphics{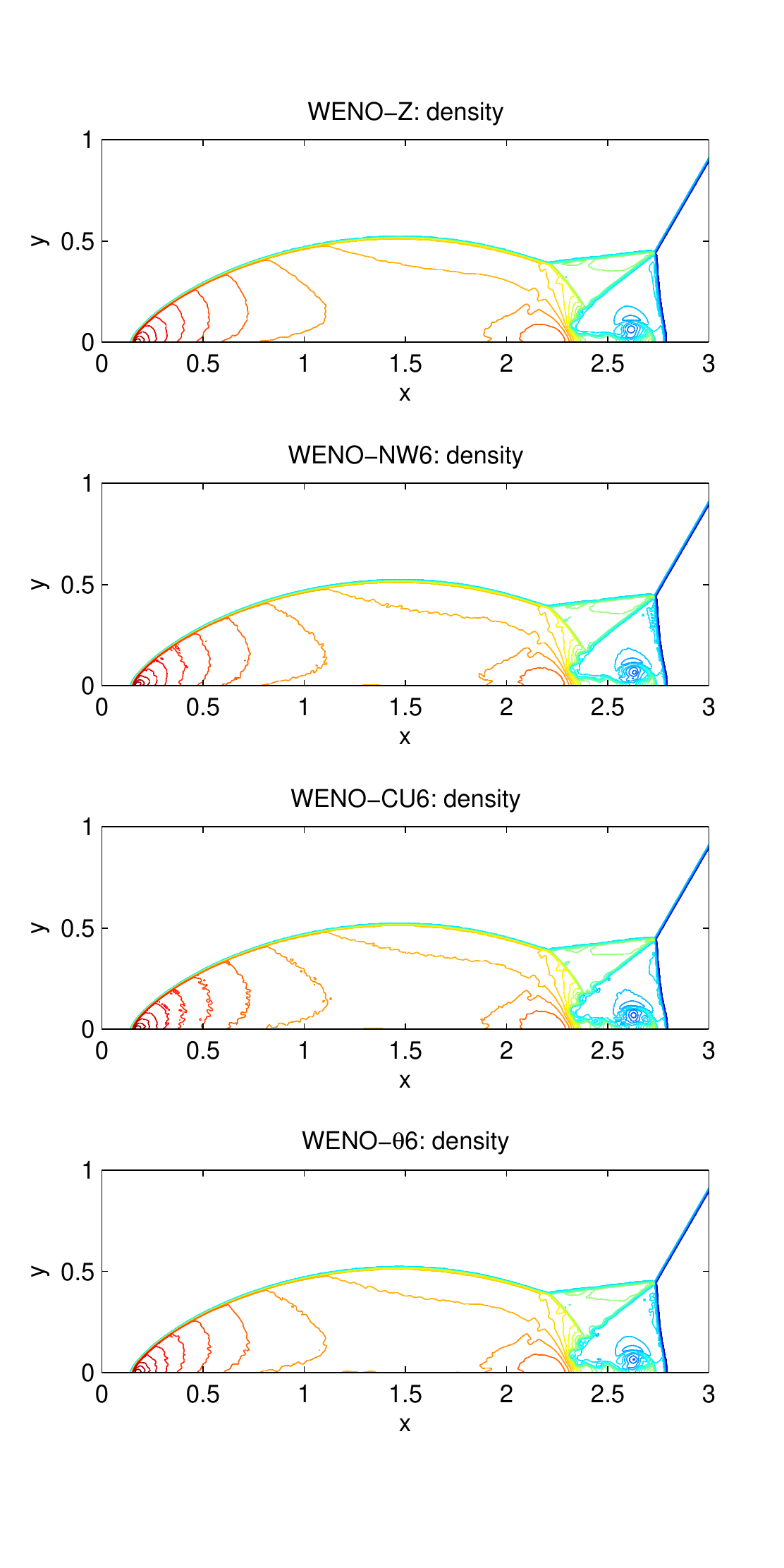}}
    \end{tabular}
    \caption{\small{The double-Mach reflection problem with initial
    data (\ref{eu2D_mach}). Density with $30$ contours. Time $t=0.2$. Grid $800\times200$.
    From top to bottom, respectively: WENO-Z, WENO-NW6, WENO-CU6, and WENO-$\theta$6.}}
    \label{fig_eu2D_mach}
  \end{center}
\end{figure}

\begin{figure}
  \begin{center}
    \begin{tabular}{cc}
      \resizebox{70mm}{!}{\includegraphics{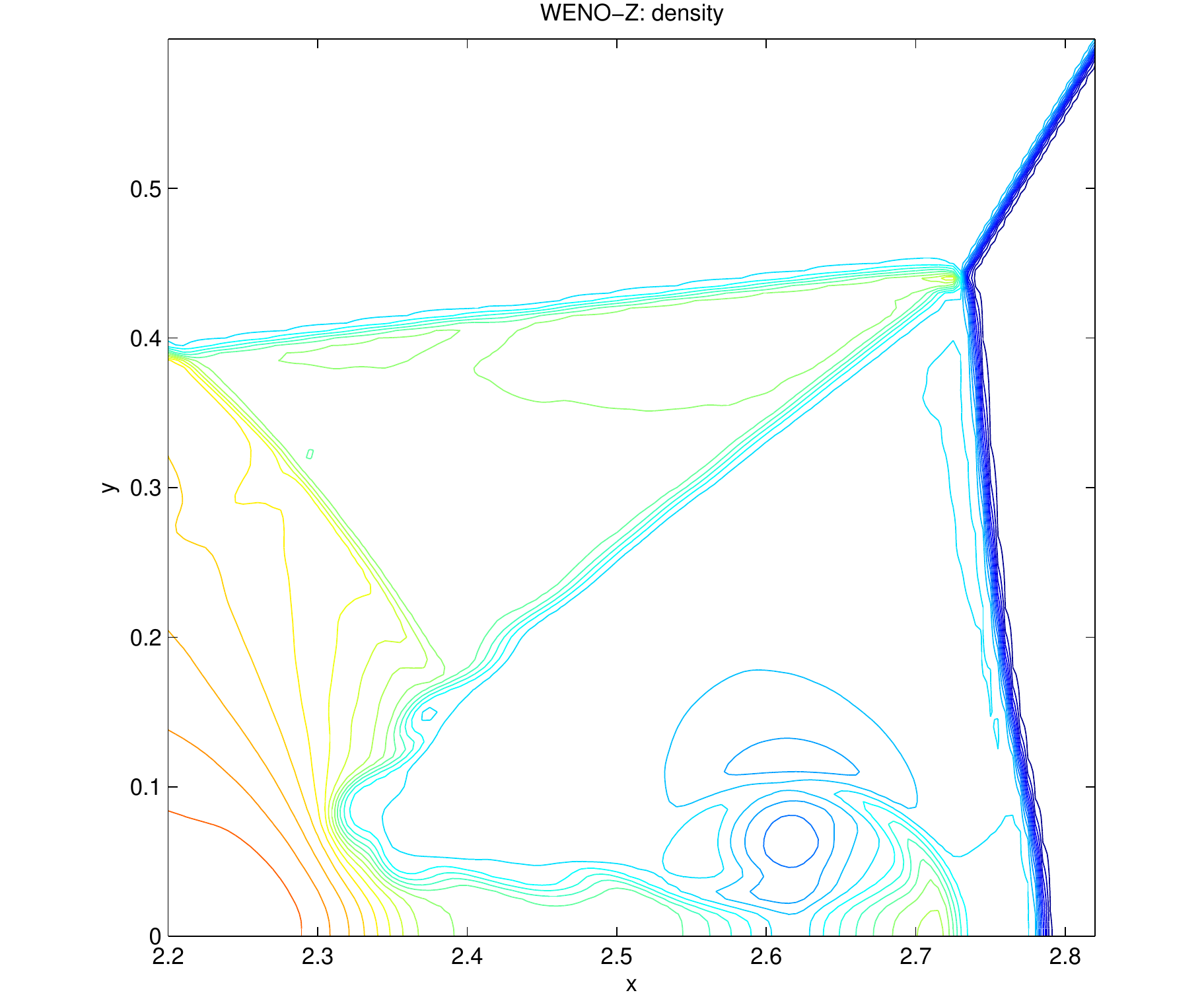}}&
      \resizebox{70mm}{!}{\includegraphics{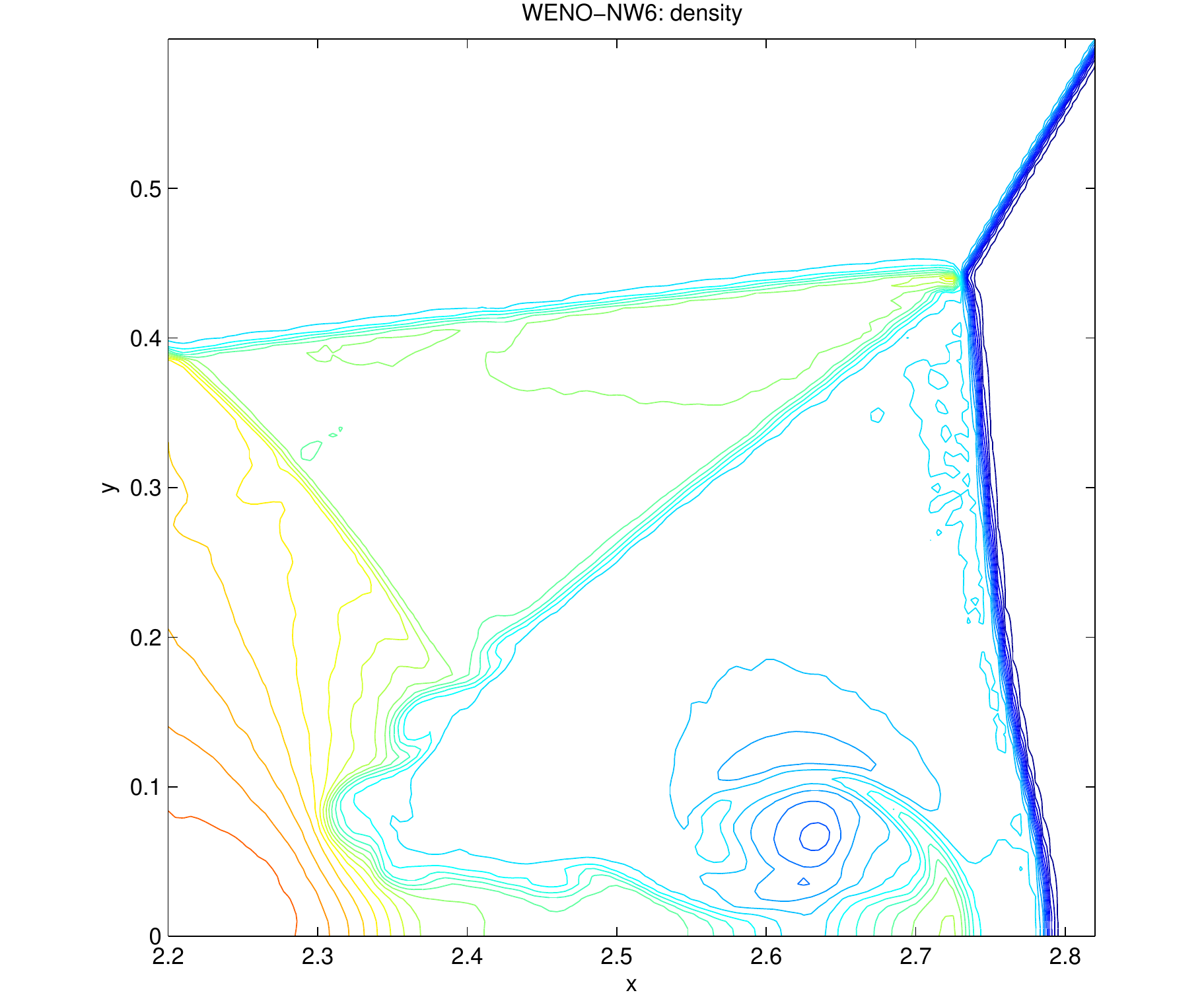}}\\
      \resizebox{70mm}{!}{\includegraphics{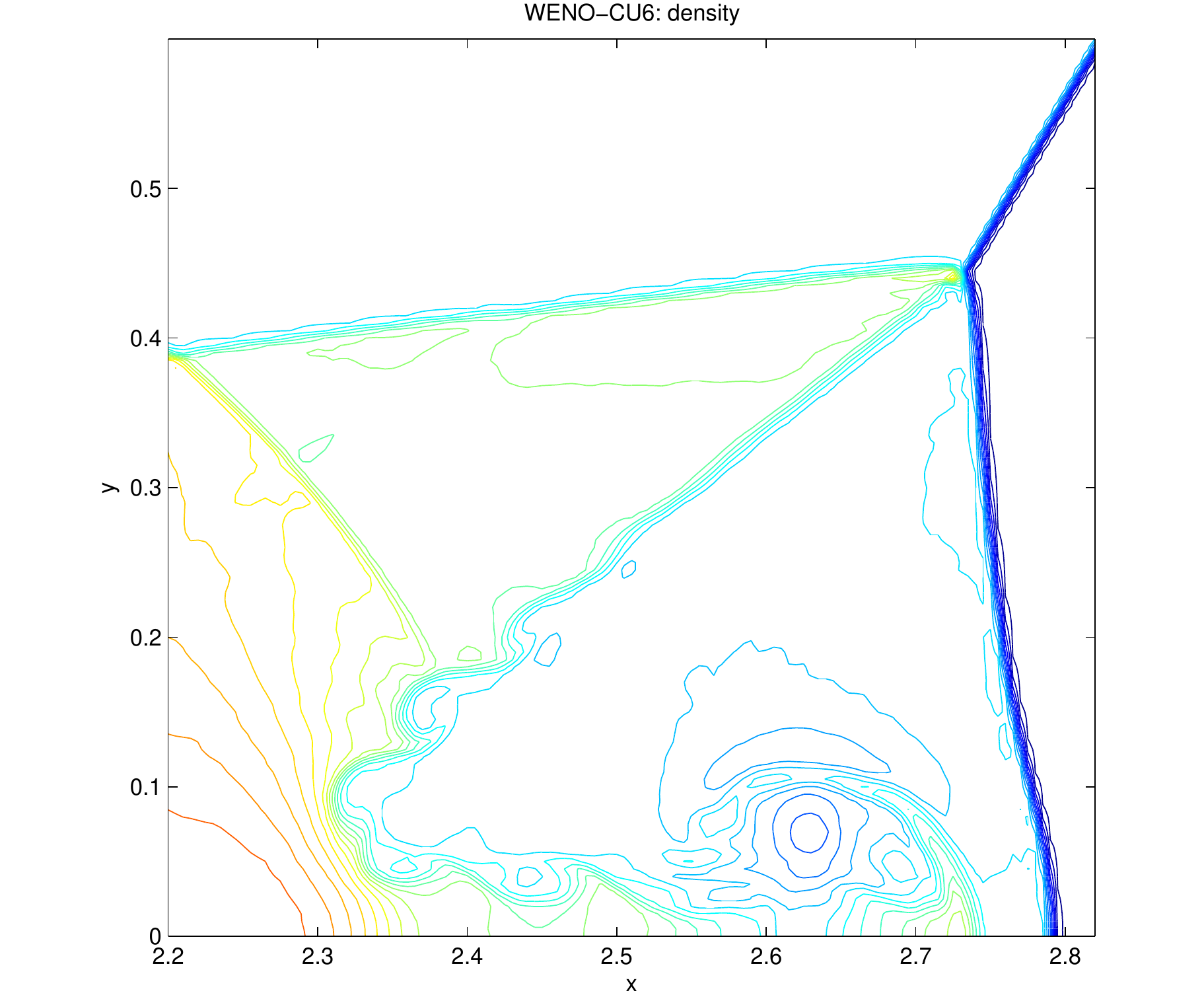}}&
      \resizebox{70mm}{!}{\includegraphics{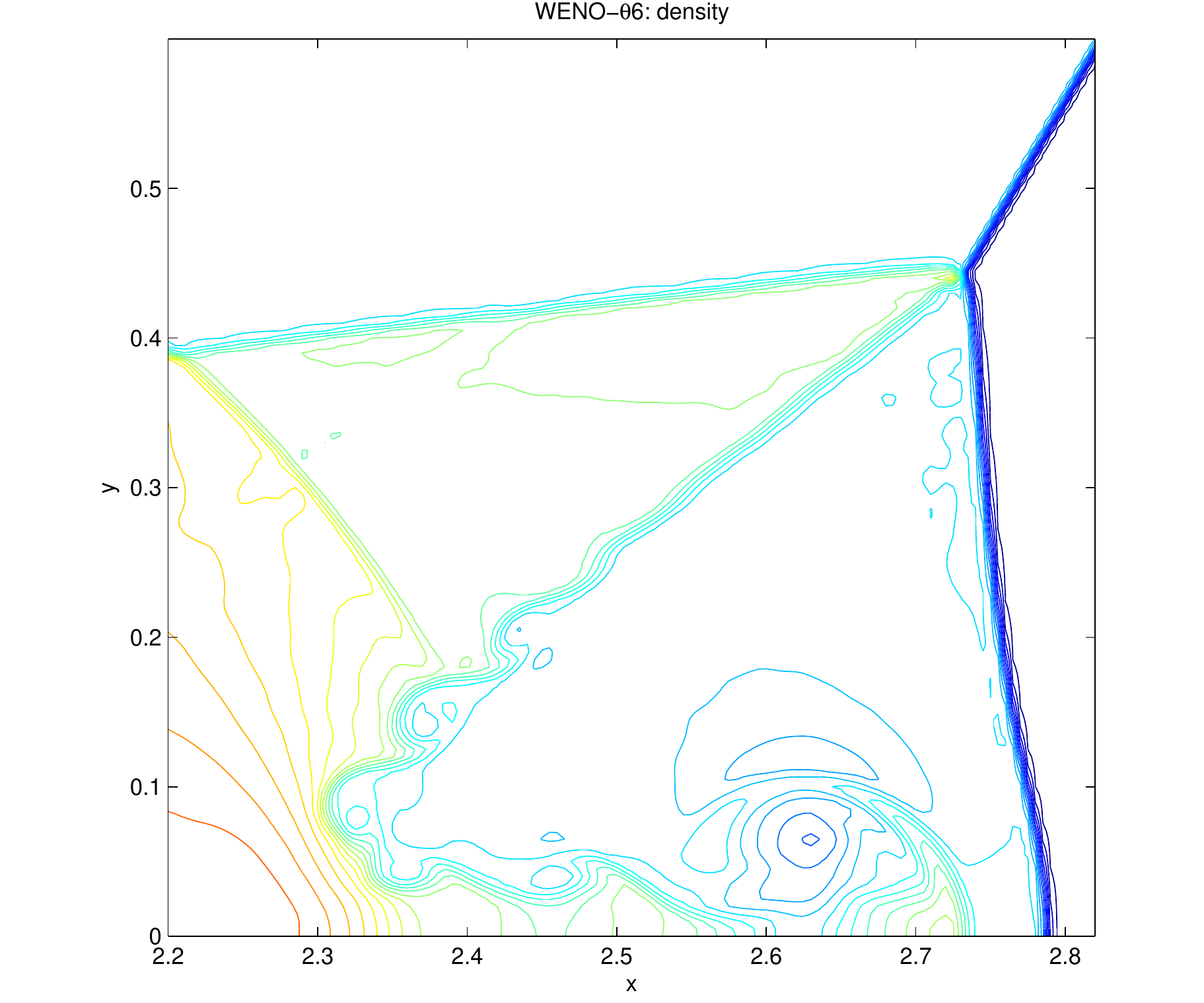}}
    \end{tabular}
    \caption{\small{The double-Mach reflection problem with initial data
    (\ref{eu2D_mach}). Zoom at the double Mach stems region.
    From top to bottom, respectively: WENO-Z, WENO-NW6, WENO-CU6, and WENO-$\theta$6.}}
    \label{fig_eu2D_mach_zoom}
  \end{center}
\end{figure}

Numerical results of the density obtained from the 5th-order WENO-Z,
the 6th-order WENO-NW6, WENO-CU6, and WENO-$\theta$6 schemes at time
$t=0.2$ are plotted in Fig. \ref{fig_eu2D_mach} with $30$ contours.
For this case, we choose a fine grid of $800\times200$ points for
all schemes. We notice the rendering of small vortices at the end of
the slip line and the wall jet, starting from WENO-Z and becoming
clearer for the 6th-order schemes. The zoom-in on the Mach stems
region shown in Fig. \ref{fig_eu2D_mach_zoom} reveals that the
WENO-$\theta$6 scheme gives more satisfactory resolution than the
WENO-NW6 and WENO-CU6 ones.

\section{Conclusion}

In this work, we have presented a new WENO-$\theta$ scheme which
adaptively switches between a 5th-order upwind and 6th-order central
scheme, depending on the smoothness of not only the sub-stencils but
also the large one. Unlike the other 6th-order WENO methods in which
this switch depends solely on the smoothness of the most downwind
sub-stencil, it is that of the large stencil which decides this
mechanism in our new scheme. \emph{Main features of our new scheme
are that the new scheme combines good properties of both 5th-order
upwind and 6th-order central schemes. That is, the new scheme is
more dispersive than the 5th-order ones in terms of better
resolution of small-scaled structures and capturing discontinuities.
Moreover, the scheme overcomes the loss of accuracy around some
critical regions and has the ability to maintain symmetry in the
solutions which are drawbacks of other comparing 6th-order WENO
schemes.}

We have also developed the new smoothness indicators of the
sub-stencils $\tilde{\beta}_k$'s which are symmetric in terms of
Taylor expansions around the point $x_{j}$ and a new $\tau^{\theta}$
for the large stencil. The latter is chosen as the smoother one
among two candidates which are computed based on the possible
highest-order variations of the reconstruction polynomials in $L^2$
sense. From then, value of the parameter $\theta$ is determined to
decide if the scheme is 5th-order upwind or 6th-order central.

A number of numerical tests for both scalar cases, linear and
nonlinear, and system case with the Euler equations of gas dynamics
are carried out to check the accuracy, resolution, and robustness of
our new scheme. It is shown that our new method is more accurate
than WENO-JS, WENO-Z, and WENO-CU6; more robust than WENO-CU6 and
WENO-NW6; and outperforms comparing schemes in capturing
small-scaled structures, and around critical regions.

Numerical simulations of higher dimensional problems will be
investigated in a subsequent work. Since the new smoothness
indicators $\tilde{\beta}_3$ and $\tau^{\theta}$ are constructed in
a systematic manner, we expect a development of the scheme to a
higher order of accuracy. This will be considered in our future
work.

\section*{Acknowledgments}

This work was supported by the National Research Foundation of Korea
(NRF) grant funded by the Korea government (MSIP)
(2012R1A1B3001167). The authors would like to thank Prof. Chi-Wang
Shu for generously giving us the WENO-JS codes for the 2D Euler
system case, and Prof. Xiangyu Hu for pointing out the typo in the
WENO-CU6 scheme.

\pagebreak

\end{document}